\documentclass[oneside, a4paper,11pt,reqno]{amsart}
\usepackage{color}
\usepackage{amsthm}
\usepackage{graphics,dsfont}
\usepackage[T1]{fontenc}
\usepackage{tikz}
\usepackage{euscript}
\usepackage{amsmath,amssymb}
\usepackage{dsfont}
\usepackage{graphicx,psfrag}%,qsymbols}
\usepackage{multicol}

\newcommand{\mysection}{\setcounter{equation}{0} \section}
\newtheorem{theo}{Theorem}[section]

\newtheorem{lemma}[theo]{Lemma}

\newtheorem{fact}[theo]{Fact}
\newtheorem{deff}[theo]{Definition}
\newtheorem{ex}{Example}

\def\un{ \mathds{1}}
\def\a{\alpha}

\def\d{\delta}
\def\e{\varepsilon}
\def\E{\mathbb{E}}
\def\o{\omega}
\def\F{\mathcal{F}}

\def\N{\mathbb{N}}
\def\P{\mathbb{P}}

\def\R{\mathbb{R}}

\def\t{\theta}
\def\tn{\overline{\theta}}
\def\tt{\overline{\theta}}
\def\Z{\mathbb{Z}}

\def\VV{\widehat{V}}

\def\po{P_{\omega}}
\def\eo{E_{\omega}}

\def\dd{\textnormal{d}}

\title[Persistence of some additive functionals of Sinai's walk]
{Persistence of some additive functionals of Sinai's walk}

\date{\today}

\author{Alexis Devulder}
\address{Universit\'e de Versailles Saint-Quentin-en-Yvelines, Laboratoire de Math\'ematiques de
Versailles, CNRS UMR 8100, Bât. Fermat,
45 avenue des Etats-Unis,
78035 Versailles Cedex, France.}
\email{devulder@math.uvsq.fr }

 %60F10
\keywords{Random walk in random
environment, Sinai's walk, integrated random walk, one-sided exit problem, persistence, survival exponent.\\
This research was partially supported by the french ANR project MEMEMO2 2010 BLAN 0125.}

\begin{document}

%\date{}

\begin{abstract}
%Sinai proved in 1992 that for a simple random walk $(R_n)_{n\in\N}$, $\P(\min_{1\leq k\leq N}\sum_{k=1}^n %R_k>0)\asymp N^{-1/4}$ as $N\to+\infty$.
%We are interested in the corresponding probability for Sinai's walk $(S_n)_{n\in\N}$. More generally we prove that %the annealed probability that
We are interested in Sinai's walk $(S_n)_{n\in\N}$. We prove that the annealed probability that
$\sum_{k=0}^n f(S_k)$ is strictly positive for all $n\in[1,N]$ is equal to $1/(\log N)^{\frac{3-\sqrt{5}}{2}+o(1)}$, for a large class of functions $f$, and in particular for $f(x)=x$.
The persistence exponent $\frac{3-\sqrt{5}}{2}$ first appears in a non-rigorous paper of Le Doussal, Monthus and Fischer, with motivations coming from physics.
The proof relies on techniques of localization for Sinai's walk and uses results of Cheliotis
about the sign changes of the bottom of valleys of a two-sided Brownian motion.
\end{abstract}

\maketitle

%We are interested in lower tail problem for additive functionals of Sinai's walk. In particular, we obtain the %persistence of the temporally averaged particle (studied by Le Doussal et al. \cite{LMF}). To this aim we use %techniques of localization of Sinai's walk and results of Cheliotis about the sign changes of the bottom of the %valleys of a two-sided Brownian motion.
%\end{abstract}

%\noindent\textbf{\sc{Key Words:}}\emph{ Random walk in random
%environment, Sinai's walk, integrated random walk, one-sided exit problem, persistence, survival exponent.}

%\bigskip
%\noindent\textbf{AMS $(2010)$ Classification:} 60K37,%60F10,60J55.

%%%%%%%%%%%%%%%%%%%%%%%%%%%%%%%%%%%%%%%%%%%%%%%%%%%%%%%%%%%%%%%%
%                                                              %
%                           SECTION 1                          %
%                                                              %
%%%%%%%%%%%%%%%%%%%%%%%%%%%%%%%%%%%%%%%%%%%%%%%%%%%%%%%%%%%%%%%%

\mysection{Introduction}\label{intro}

In this paper we consider random walks in random environments in $\Z$.
 Let $\o:=(\o_i)_{i\in\Z}$ be a collection of independent and identically distributed random variables
 taking values in $(0,1)$, with joint law
$\eta$. A realization of $\o$ is called an {\it environment}. Conditionally on $\o$, we define a Markov chain $(S_n)_{n\in\N}$ by $S_0=0$ and for $n\in\N$, $k\in\Z$ and $i\in\Z$,
\begin{equation*}
\po(S_{n+1}=k|S_n=i)=\left\{
\begin{array}{ll}
\o_i & \text{ if } k=i+1,\\
1-\o_i & \text{ if } k=i-1,\\
0      & \text{ otherwise}.
\end{array}
\right.
\end{equation*}
We say that $(S_n)_{n\in\N}$ is a {\it random walk in random environment} (RWRE).
This model has many applications in physics (see e.g. Hughes \cite{Hug}) and in biology (see e.g. Cocco and Monasson \cite{CoccoMonasson} about DNA reconstruction), and has unusual properties. Moreover, its properties are used to study several other mathematical models, see e.g. Zindy \cite{Zindy}, Enriquez, Lucas and Simenhaus \cite{ELS} and Devulder \cite{Devulder_SPL}.

The probability $\po$ is called the {\it quenched law}. We denote by $\po^x$ the quenched law for a RWRE starting at $x\in\Z$ instead of $0$. We also consider the {\it annealed law}, which is defined by
$$\P(.)=\int\po(.)\eta(\textnormal{d}\o).$$
Notice in particular that $(S_n)_{n\in\N}$ is not Markovian under $\P$. We also denote by $\E$, $\eo$ and $\eo^x$ the expectations under $\P$, $\po$ and $\po^x$ respectively. We assume that the following ellipticity condition holds:
\begin{equation}\label{elliptic}
\exists \e_0\in(0,1/2), \qquad \eta(\e_0\leq \o_0\leq 1-\e_0)=1.
\end{equation}
This ensures that $|\log(\frac{1-\o_0}{\o_0})|$ is $\eta$-a.s. bounded by $\log(\frac{1-\e_0}{\e_0})$.
Solomon \cite{S1} proved that $(S_n)_{n\in\N}$ is recurrent for almost every environment $\o$ if and only if
\begin{equation}\label{eqrecA}
\int\log\left(\frac{1-\o_0}{\o_0}\right)\eta(\dd \o)  =  0.\\
\end{equation}
We assume that this condition is satisfied throughout the paper.
Moreover, in order to avoid the degenerate case of simple random walks, we suppose in the following that
\begin{equation}\label{variance}
\sigma:=\left(\int\log^2\left(\frac{1-\o_0}{\o_0}\right)\eta(\dd \o)\right)^{1/2} >  0.
\end{equation}
A RWRE $(S_n)_{n\in\N}$ satisfying conditions \eqref{elliptic}, \eqref{eqrecA} and \eqref{variance}
 is referred to as {\it Sinai's walk}. Sinai (\cite{S3}, see also Andreoletti \cite{Papier_4_A} for extensions) proved that in this (recurrent) case,
$$
\sigma^2\frac{S_n}{\log^2 n}\to_{\text{law}}b_{\infty}
$$
as $n\to+\infty$, where $b_{\infty}$ is a non degenerate and non gaussian random variable and $\to_{\text{law}}$ denotes convergence in law under $\P$. We refer to  Hughes \cite{Hug}, R{\'e}v{\'e}sz \cite{Revesz}  and Zeitouni \cite{Z1} for more properties of RWRE.

Sinai \cite{SinaiSRW} also showed in 1992 that for a symmetric simple random walk $(R_n)_{n\in\N}$, we have
$\P(\sum_{k=1}^n R_k>0\ \forall 1\leq n\leq N)\asymp N^{-1/4}$ as $N\to+\infty$.
In this paper, we are interested in the corresponding probability for Sinai's walk $(S_n)_{n\in\N}$, and more generally in the one-sided exit problem for some additive functionals of Sinai's walk under the annealed law $\P$.
%In this paper, we are interested in the one-sided exit problem for some additive functionals of Sinai's walk under %the annealed law $\P$.
We say that $g(x)=o(1)$ as $x\to+\infty$ (resp $-\infty$) if $g(x)\to 0$ as $x\to+\infty$ (resp $-\infty$).
Our main result is the following.

\begin{theo}\label{TheoremePersistenceAdditive}
Let $f$ be a function $\Z\to\R$, such that
$f(0)=0$;
%$\forall x>0, f(x)\geq 1$;
$f(x)\geq 1$ for all $x>0$;
%$\forall x<0, f(x)\leq -1$
$f(x)\leq -1$ for all $x<0$;
and $|f(x)|\leq\exp(|x|^{o(1)})$ as $x\to\pm\infty$.
%and $f(x)\geq -\exp(|x|^{o(1)})$ as $x\to-\infty$.
We consider a RWRE $(S_n)_{n\in\N}$ satisfying conditions \eqref{elliptic}, \eqref{eqrecA} and \eqref{variance}, and a real number $u\leq 0$. We have as $N\to+\infty$,
$$\P\left(\sum_{k=0}^n f(S_k)>u\quad  \forall 1\leq n\leq N\right)
=\frac{1}{(\log N)^{\frac{3-\sqrt{5}}{2}+o(1)}}.$$
\end{theo}

\medskip

Let $(A_t)_{t\in D}$ be a real valued stochastic process starting from $0$, where $D=\R_+$ or $D=\N$.
The asymptotic study of the survival function $\P(\forall t\in(0,T]\cap D,  A_t\geq x)$ for $x\leq 0$, when $T\to+\infty$, is called {\it one sided exit problem} or {\it persistence probability}. This problem is equivalent to the study of $\P(T_x>T)$,
where $T_x$ is the first passage time of the process $(-A_t)_t$ strictly above the level $y=-x\geq 0$.
In many cases with physical relevance, the survival function  behaves asymptotically like $1/T^{\alpha+o(1)}$ as $T\to+\infty$, with $\a>0$. The exponent $\alpha$ is called the {\it persistence} or {\it survival} exponent. This problem, which is well known for random walks or L\'evy processes, is less understood for the integrals of these processes, in particular in the discrete case. We refer to Aurzada and Simon \cite{AurzadaSimon} for a recent review on this subject from the mathematical point of view. Persistence properties have also received a considerable attention in physics, see e.g. Bray, Majumdar and Schehr \cite{Bray} for an up-to-date survey.

In our case, the probability we obtain in Theorem \ref{TheoremePersistenceAdditive} for the integrals of $(f(S_n))_n$ is a power of $\log N$ instead of $N$, which is quite unusual
and contrasts with all the cases presented in the review paper \cite{AurzadaSimon}.
%The explicit value of the exponent $\alpha$ is known in few cases, especially for non-Markov processes, see e.g. %Simon \cite{Simon} (verifier).
The value of the survival exponent is $\frac{3-\sqrt{5}}{2}$; it does not depend on the function $f$ for a wide class of functions, and it also
does not depend on the law $\eta$ of the environment, as long as \eqref{elliptic}, \eqref{eqrecA} and \eqref{variance} are satisfied.
It is derived from the results of Cheliotis \cite{Cheliotis} about the number of sign changes of the bottom of valleys of Brownian motion,
and was first stated in a non rigorous paper of Le Doussal, Monthus and Fisher \cite{LMF}, with motivations coming from physics.
%Notice also that the strict inequality $>u$ can be replaced by a large one in Theorem \ref{TheoremePersistenceAdditive}.

Before giving some examples, we introduce some more notation.
We denote by
%$\Z$ the set of integers. $\N$ is the set of nonnegative integers,
$\N^*$ the set of positive integers, and $\Z_-^*$ is the set of negative integers.
We define the local time of the RWRE $(S_n)_{n\in\N}$ at time $n\in\N$ as follows:
\begin{equation*}
L(A, n)  :=  \sum_{k=0}^n \mathds{1}_{\{S_k\in A\}},\qquad  L(x,n):=L(\{x\},n)
\end{equation*}
for $A\subset \Z$ and $x\in\Z$. In words $L(A,n)$ is the number of visits of the random walk $S$ to the set $A$ in the first $n$ steps. This quantity will be useful in the proof of Theorem \ref{TheoremePersistenceAdditive}, because
\begin{equation}\label{eqFormuleOccupation}
\sum_{k=0}^n g(S_k)=\sum_{x\in\Z} g(x) L(x,n),\qquad n\in\N,
\end{equation}
for every function $g$.

It can be interesting to keep in mind the first example:

%For a set $A\subset \Z$ and $n\in\N$, we define $L(A, n):=\sum_{k=0}^n \mathds{1}_{\{S_k\in A\}}$, and call it the
%occupation time of $A$ at time $n\in\N$. If $x\in\Z$, we also define $L(x,n):=L(\{x\},n)$ the local time of $x$ at time $n$.

\begin{ex}
For $f(x)=\mathds{1}_{\{x>0\}}-\mathds{1}_{\{x<0\}}$, Theorem \ref{TheoremePersistenceAdditive} gives
$$\P\left[L(\N^*, n)>L(\Z_-^*,n)\quad  \forall 1\leq n\leq N\right]
=\frac{1}{(\log N)^{\frac{3-\sqrt{5}}{2}+o(1)}}.$$
\end{ex}
The following example gives for $\alpha=1$ the persistence of the {\it temporal average}  or {\it running average} of Sinai's walk,
that is $\frac{1}{n}\sum_{k=0}^n S_k$, with the terminology of
Le Doussal et al. (\cite{LMF} Section IV):

\begin{ex}\label{ex2}
Let $\alpha>0$, $\textnormal{sgn}(x):=\mathds{1}_{\{x> 0\}}-\mathds{1}_{\{x<0\}}$ for $x\in\R$,
and $f(x)=\textnormal{sgn}(x)|x|^\alpha$ for $x\in\Z$. We get for $u\leq 0$,
$$
    \P\left(
        %\frac{1}{n}
        \sum_{k=0}^n \textnormal{sgn}(S_k)|S_k|^\alpha>u\quad  \forall 1\leq n\leq N
    \right)
=
    \frac{1}{(\log N)^{\frac{3-\sqrt{5}}{2}+o(1)}}.$$
\end{ex}
We recall that the corresponding probability for $\alpha=1$ for a simple random walk is of order $N^{-1/4}$ (see Sinai \cite{SinaiSRW}; see also Vysotsky \cite{Vito} and Dembo, Ding and Gao \cite{DemboGao} for recent extensions). Example \ref{ex2} is also, for $\alpha>0$ arbitrary, the analogue for Sinai's walk of the results obtained by Simon \cite{Simon} for some additive functionals of stable processes with no negative jumps. We can also consider functions increasing more rapidly, such as $f(x)=\textnormal{sgn}(x)|x|^{|\log (2+|x|)|^\alpha}$, $x\in\Z$ for $\alpha>0$.

The rest of the paper is organized as follows.
We introduce some notation and basic facts in Section \ref{SectNotation}.
In Section \ref{SectionUpper} we build a set $\mathcal{B}(N)$ of {\it bad environments}, such that in a bad environment,
$\sum_{k=0}^n f(S_k)$ is less than $u$ for at least one integer $n\in[1,N]$ with a great quenched probability.
To this aim, we approximate the potential of the environment by a two-sided Brownian motion, and
we define {\it strong changes of sign} for the valleys of this Brownian motion. We prove that in a bad environment, the existence of such a strong change of sign  forces the walk to stay a long time in $\Z_-^*$ with a large quenched probability, leading to the upper bound of Theorem \ref{TheoremePersistenceAdditive}.
A sketch of this proof is provided in Subsection \ref{SubSectSketchProofPart3}.
In Section \ref{good} we build a set $\mathcal{G}(N)$ of {\it good environments}. We prove, using a mathematical induction,
%Section \ref{SectionLower}
that in such a good environment $\sum_{k=0}^n f(S_k)$ is strictly positive for all $1\leq n\leq N$
with a large quenched probability, which leads to the lower bound of Theorem \ref{TheoremePersistenceAdditive}.
A sketch of this proof is given in Subsection \ref{SubSectSketchProofPart4}.
Finally, Section \ref{SectionPreuveduLemme21} is devoted to the proof of two technical lemmas.

Throughout the paper, $c_i, i\in\N$, denote positive constants, and $\log$ denotes the natural logarithm.

%%%%%%%%%%%%%%%%%%%%%%%%%%%%%%%%%%%%%%%%%%%%%%%%%%%%%%%%%%%
%                                                         %
%                SOME NOTATION                            %
%                                                         %
%%%%%%%%%%%%%%%%%%%%%%%%%%%%%%%%%%%%%%%%%%%%%%%%%%%%%%%%%%%

\mysection{Preliminaries}\label{SectNotation}

\subsection{Potential}
We recall that the potential $V$ is a function of the environment $\o$, which is defined on $\Z$ as follows:
\begin{equation*}
  V(n):=\left\{
  \begin{array}{lr}
  \sum_{i=1}^n \log\frac{1-\o_i}{\o_i} &
  \textnormal{if } n>0,\\
  0 & \textnormal{if }n=0,\\
  -\sum_{i=n+1}^0 \log\frac{1-\o_i}{\o_i} &  \textnormal{if } n<0.
  \end{array}
  \right.
\end{equation*}
For $p\in\Z$, we define the hitting time of $p$ by $(S_n)_n$ by:
\begin{equation*}
\tau(p):=\inf\{k\in\N,\ S_k=p\}.
\end{equation*}
We now recall some basic estimates that will be useful throughout the paper.

\begin{lemma}\label{LemmaProbaAtteinte}
(see e.g. Zeitouni \cite{Z1} formula (2.1.4) p. 196)
If $p<q<r$, then
\begin{equation}\label{eqProbaAtteinte}
\po^q[\tau(r)<\tau(p)]=\Big(\sum_{k=p}^{q-1}e^{V(k)}\Big)\Big(\sum_{k=p}^{r-1}e^{V(k)}\Big)^{-1}.
\end{equation}
\end{lemma}

\smallskip

\begin{lemma}(see e.g. Zeitouni \cite{Z1} p. 250)
If $g<h<i$,
\begin{equation}\label{InegEsperanceZeitouni}
    \eo^{h} [\tau(g)    \wedge\tau(i)]
\leq
    \sum_{k={h}}^{i-1}\sum_{\ell=g}^{k}\frac{\exp[V(k)-V(\ell)]}{\o_{\ell}}
\leq
    \e_0^{-1}(i-g)^2\exp\left[\max_{g\leq\ell\leq k\leq i-1}(V(k)-V(\ell))\right].
\end{equation}
\end{lemma}

\smallskip

\noindent{\bf Proof:} This formula \eqref{InegEsperanceZeitouni} is proved by Zeitouni \cite{Z1} p. 250,
in the particular case $h=0$. Indeed, Zeitouni uses the notation
$\omega_x^+:=\omega_x$, $\rho_x:=(1-\omega_x)/\omega_x$, $x\in\Z$ (see \cite{Z1} p. 194 and p. 195),
$\overline{T}_{b,n}:=\tau(a_\delta^n)\wedge \tau(b^n)$
for some $a_\delta^n<0<b_\delta^n=b^n$,
and proves in the fourth formula of \cite{Z1} p. 250 that
\begin{equation}\label{eqZeitouniP250}
    \eo^{0} \big[\tau(a_\delta^n)\wedge \tau(b^n)\big]
=
    \eo^{0} \big[\overline{T}_{b,n}\big]
\leq
    \sum_{i=1}^{b^n}\, \sum_{j=0}^{i-1-a_\delta^n}\frac{\prod_{k=1}^{j}\rho_{i-k}}{\omega_{i-j-1}}
=
    \sum_{k=0}^{b^n-1}\,\sum_{\ell=a_\delta^n}^{k}\frac{\exp[V(k)-V(\ell)]}{\omega_\ell}.
\end{equation}
Since the proof of this formula does not use any property of $a_\delta^n$ and $b^n$ except $a_\delta^n<0<b_n^\delta=b^n$,
it is true for any integers $a_\delta^n<0<b^n$. The general case
\eqref{InegEsperanceZeitouni} follows from \eqref{eqZeitouniP250} by translation, since
$
    \eo^h[\tau(g) \wedge  \tau(i)]
=
    E_{\widehat{\omega}}^0[\tau(g-h) \wedge  \tau(i-h)]
$
for $g<h<i$, with $\widehat{\omega}_x:=\omega_{x+h}$ for every $x\in\Z$.
\hfill$\Box$

\smallskip

Moreover, the following estimate can be found in Andreoletti (\cite{Andreoletti_PhD} p. 22) and is in the spirit of R{\'e}v{\'e}sz (\cite{Revesz} p278-279).

\smallskip

\begin{lemma}\label{LemmaRevesz1}
If $p<z\leq q<r$ or  $p<q<z<r$,
\begin{equation}\label{eqLemmaRevesz1}
%\label{eqLemmaRevesz2}
    \eo^q[L(z,\tau(p)\wedge\tau(r))]
=
\frac{\po^q[\tau(z)<\tau(p)\wedge \tau(r)]}
{\o_z\po^{z+1}[\tau(z)>\tau(r)]+(1-\o_z)\po^{z-1}[\tau(z)>\tau(p)]}.
\end{equation}
\end{lemma}

\noindent
For the sake of completeness, we recall the proof:
\medskip

\noindent{\bf Proof of Lemma \ref{LemmaRevesz1}:}
By the strong Markov property,
\begin{eqnarray*}
\eo^q[L(z,\tau(p)\wedge\tau(r))]
&=&
\eo^q[L(z,\tau(p)\wedge\tau(r))\mathds{1}_{\{\tau(z)<\tau(p)\wedge\tau(r)\}}]\\
& = & \eo^z[L(z,\tau(p)\wedge\tau(r))]\po^q[\tau(z)<\tau(p)\wedge\tau(r)].
\end{eqnarray*}
Since $L(z,\tau(p)\wedge\tau(r))$ is under $\po^z$ a geometric random variable of parameter
$\o_z\po^{z+1}[\tau(z)>\tau(r)]+(1-\o_z)\po^{z-1}[\tau(z)>\tau(p)]$, we get \eqref{eqLemmaRevesz1}.
\hfill$\Box$

\subsection{$x$-extrema}
We now recall some definitions introduced by Neveu and Pitman \cite{NP}. %also used by Cheliotis \cite{Cheliotis}.
If $w$ is a continuous function $\R\to\R$, $x>0$, and $y_0\in\R$, it is said that $w$ admits an {\it $x$-minimum at $y_0$} if there exists real numbers $\alpha$ and $\beta$ such that $\alpha<y_0<\beta$, $w(y_0)=\inf\{w(y),\ y\in[\alpha,\beta]\}$, $w(\alpha)\geq w(y_0)+x$ and $w(\beta)\geq w(y_0)+x$. It is said that $w$ admits an {\it $x$-maximum} at $y_0$ if $-w$ admits an $x$-minimum at $y_0$. In these two cases we say that $w$ admits an {\it $x$-extremum} at $y_0$.

We denote by $\mathcal{W}$ the set of functions $w$ from $\R$ to $\R$ such that the three following conditions are satisfied:
{\bf (a)} $w$ is continuous on $\R$;
{\bf (b)} for every $x>0$, the set of $x$-extrema of $w$ can  be written
$\{x_k(w,x),\ k\in\Z\}$, where $(x_k(w,x))_{k\in\Z}$ is strictly increasing, unbounded from above and below, with $x_0(w,x)\leq 0<x_1(w,x)$;
{\bf (c)} for all $x>0$ and $k\in\Z$, $x_{k+1}(w,x)$ is an $x$-maximum if and only if $x_{k}(w,x)$ is an $x$-minimum.
We now consider a two-sided standard Brownian motion $W$.
We know from Cheliotis (\cite{Cheliotis}, Lemma 8) that $\eta(W\in\mathcal{W})=1$.

%Cheliotis \cite{Cheliotis} proves that almost surely, for every $x>0$, the set of $x$-extrema of $W$ can  be written
%$\{x_k(W,x),\ k\in\Z\}$, where $(x_k(W,x))_{k\in\Z}$ is strictly increasing, unbounded from above and below, and $x_0(W,x)\leq 0<x_1(W,x)$.
%Moreover, almost surely for all $x>0$ and $k\in\Z$, $x_{k+1}(W,x)$ is an $x$-maximum if and only if $x_{k}(W,x)$ is an $x$-minimum.

\noindent
For each $x>0$, $b_W(x)$, also denoted by $b(x)$ when no confusion is possible, is defined on $\{W\in\mathcal{W}\}$ as
\begin{equation*}
b_W(x):=\left\{\begin{array}{ll}
x_0(W,x)\ \ \ \ & \text{if } x_0(W,x) \text{ is an $x$-minimum},\\
x_1(W,x)\ \ \ \ & \text{otherwise}.\\
\end{array}
\right.
\end{equation*}
One interesting feature about $b_W$ is that the diffusion in the random potential $W$, defined by Schumacher \cite{Schumacher}, is localized in a small neighborhood of $b_W(\log t)$ at time $t$ with probability nearly one (see Brox \cite{Brox}, Tanaka \cite{{tanaka}} and  Hu \cite{HuLocal}).
Such a diffusion can be viewed as a continuous time analogue of Sinai's walk (see e.g. Shi \cite{S2}),
and a similar localization phenomenon arises for Sinai's walk (see Sinai \cite{S3}, Golosov \cite{Golosov} and more recently Andreoletti \cite{Papier_4_A}).

For $x>0$ and $k\in\Z$, the restriction of $W-W(x_k(W,x))$ to $[x_k(W,x),x_{k+1}(W,x)]$ is denoted by $T_k(x)$ and is called an {\it $x$-slope}. It is the translation of the trajectory of $W$ between two consecutive $x$-extrema.  If $x_k(W,x)$ is an $x$-minimum (resp. $x$-maximum), it is a nonnegative (resp. non-positive) function, and its maximum (resp. minimum) is attained at $x_{k+1}(W,x)$.
For each $x$-slope $T_k(x)$, we denote by $H(T_k(x))$ its {\it height} and by $e(T_k(x))$ its {\it excess height}, that is
$H(T_k(x)):=|W[x_{k+1}(W,x)]-W[x_k(W,x)]|\geq x$ and
$e(T_k(x)):=H(T_k(x))-x\geq 0$. We also define $e(T_k(0))=H(T_k(0))=0$, $k\in\Z$.
%$e(T_k(x)):=|W[x_{k+1}(W,x)]-W[x_k(W,x)]|-x$.

The point of view of $x$-extrema has been used in some recent studies of processes in random environment, see e.g.
Bovier and Faggionato \cite{Bovier_Faggionato} for Sinai's walk,
Cheliotis \cite{Cheliotis_Favorite} for (recurrent) diffusions in a Brownian potential,
and Andreoletti and Devulder \cite{Andreoletti_Devulder_1} for (transient) diffusions in a drifted Brownian potential.

%%%%%%%%%%%%%%%%%%%%%%%%%%%%%%%%%%%%%%%%%%%%%%%%%%%%%%%%%%%
%                                                         %
%              PROOF OF THE UPPER BOUND                   %
%                                                         %
%%%%%%%%%%%%%%%%%%%%%%%%%%%%%%%%%%%%%%%%%%%%%%%%%%%%%%%%%%%

\mysection{Proof of the upper bound}\label{SectionUpper}

\subsection{Sketch of the proof, and organization of this proof}\label{SubSectSketchProofPart3}
We approximate the potential $V$ in \eqref{eqDefB1}  by $\sigma W$, where $W$ is a suitable two-sided Brownian motion.

In many cases for Sinai's walk, the environment largely controls the behavior of the random walk.
This is due to the fact that the random walk tends to go to places with a low potential, and spend a large amount of time around these places.
So, heuristically speaking, the idea is to prove that for most environments,
the deepest location (in terms of potential) visited until time $n$
is $<0$ for at least one time $n\leq N$, and that
%with large quenched probability
the RWRE $(S_k)_k$ spends a large amount of time around this deepest location
before going back to the positive locations at some time $m\leq N$, making the sum
$\sum_{k=1}^m f(S_k)$ negative with large annealed probability.

One good candidate for this deepest location visited until time $n$ seems to be $b_{\sigma W}(\log n)$,
 that is, $b_{\sigma W}(x)$ for some $x$ much bigger that $1$ and much smaller than $\log N$ such that $b_W(x)<0$.
 %which exists with probability nearly $1$.
However, the existence of such an $x$ with $b_{\sigma W}(x)<0$ is not enough to ensure that
with a large quenched probability the random walk $(S_k)_k$ will go quickly to this (negative) place and spend a great
amount of time around it
before going back to $0$.
%before spending too much time in positive locations.
This is why we introduce, in Definition \ref{DefStrongChangeofSign} below, the notion of {\it $a$-strong change of sign} for $b_W$,
in order to push the walk to go quickly to negative locations and spend a large amount of time there.

We first study the potential in Subsections \ref{SubSectChangeofSign} and \ref{SubSectBadEnvironments}.
We prove in Lemma \ref{LemmeProbaBad} that with a very large probability, the environment is what we call a {\it bad environment}:
it satisfies some technical conditions, but also,
there are many changes of sign $X_k$ of $b_W$ in   $[(\log N)^\e, (\log N)^{1-\e}]$ (see \eqref{eqDefB2}),
%or equivalently $b_{\sigma W}$, for $\e>0$,
and among them, at least one is a "strong" change of sign $\mathbf{h}_N:=X_{k_N}$ of $b_W$
(see \eqref{eqDefB3} and Lemma \ref{LemmaProbaStrongChangeOfSign}),
as defined in Definition \ref{DefStrongChangeofSign} below, with 
%$b_W(\mathbf{h}_N)>0$ and 
$b_W(X_{k_{N+1}})\leq0$.
A schema representing the potential $V$ of a typical "bad environment" is given in Figure \ref{figureBad} page \pageref{figureBad}.

Then in Subsection \ref{SubSectRWinBadEnvironment}, we consider a random walk $(S_k)_k$ in such a bad environment $\o$.
Due to the conditions defining our strong change of sign $\mathbf{h}_N$, we prove that with a large quenched probability,
the random walk $(S_k)_k$ goes quickly to $x_{-1}:=x_{-1}(W,\mathbf{h}_N)\leq -1$
before going to some $v_2\leq x_{2}(W,\mathbf{h}_N)=:x_{2}$ (see Figure \ref{figureBad}).
Moreover, it stays a long time in $\Z_-^*$ before going back to $0$.
It stays such a long time in $\Z_-^*$, on which $f<0$, that
$\sum_{k=1}^n f(S_k)$ becomes $\leq u$ for some $1\leq n<N$, with large quenched probability
uniformly on bad environments (see Lemma \ref{LemmeComportementdansBad}), and so with a large annealed probability.
This leads to the upper bound of Theorem \ref{TheoremePersistenceAdditive}.

%%%%%%%%%%%%%%%%%%%%%%%%%%%%%%%%%%%%%%%%%%%%%%%%%%%%%%%%%%%%%%

%%%%%%%%%%%%%%%%%%%%%%%%%%%%%%%%%%%%%%%%%%%%%%%%%%%%%%%%%%%%%%

\subsection{Strong change of sign}\label{SubSectChangeofSign}
%Let $c\geq 1$ and let $X_1$ be the smaller
%$x\geq c$ for which the function $b_W(.)$ changes its sign.
%For $k\in\N^*$,
%$X_{k+1}$ is defined by induction as the smaller $x>X_k$ for which $b$ changes its sign.
Let $c>0$. Similarly as in Cheliotis (\cite{Cheliotis} Corollary 2), we denote by $(X_k)_{k\geq 1}$ the strictly
increasing sequence of points for which $b_W(.)$ changes its sign in $[c,+\infty)$.
The proof of the following fact is deferred to Section \ref{SectionPreuveduLemme21}:

\begin{fact} \label{FactExcessHeight0}
Almost surely,
\begin{eqnarray*}
X_1 & = & \inf\{x\geq c,\ e(T_0(x))=0\},\\
X_{k+1} & = & \inf\{x> X_k,\ e(T_0(x))=0\},\qquad k\in\N^*.
\end{eqnarray*}
Moreover, the sign of $b_W(.)$ is constant on every interval $[c,X_1]$, $(X_k, X_{k+1}]$, $k\in\N^*$.
\end{fact}

%Hence for each $k\in\N^*$, $X_k$ is a stopping time for the filtration $(\F_x)_{x\geq 0}$.
As a consequence, a.s. for every $x>0$,
$b_W$ changes its sign at $x$ if and only if $e(T_0(x))=0$.
We can now define {\it strong changes of sign} of $b_W$ as follows:

\begin{deff}\label{DefStrongChangeofSign}
Consider $x>0$.
For $a>0$, we say that $x$ is an {\it $a$-strong change of sign} of $b_W$ if and only if
\begin{equation*}
e(T_0(x))=0,\qquad  e(T_{-1}(x))\geq a x,\quad \text{and}\quad e(T_1(x))\geq a x.
\end{equation*}
%We denote by $b(x^-)$ the limit $\lim_{y\to x, y<x} b(y)$ for $x>0$.
\end{deff}

%We notice that for such an $x$, $b(x)=\lim_{y\to x, y<x} b(y)$.
In the following lemma, we evaluate the probability that there is no
$a$-strong change of sign $x$ such that
%$b(x^-)>0$
$b_W(x)>0$ in $[c,X_{2k+1})$.

\bigskip
\begin{lemma}\label{LemmaProbaStrongChangeOfSign}
For $a>0$, $c\geq 1$ and $k\in\N^*$, we define $A(k,a,c)$ also denoted by $A_{k,a,c}$ as follows:
\begin{equation*}
A_{k,a,c}:=\{\forall i\in\{1,\dots, 2k\},\ b_W(X_i)>0\Rightarrow
(e(T_{-1}(X_i))<a X_i \text{ or } e(T_1(X_i))<a X_i)\}.
\end{equation*}
We have,
\begin{equation}\label{eqLemmaProbaAka}
    \eta(A_{k,a,c})
\leq
    \eta(A_{1,a,c})
    \big(1-e^{-2a}\big)^{k-1}.
%\eta(A_{k,a,c})=\eta(A_{1,a,c})(1-e^{-2a})^{k-1}.
\end{equation}
\end{lemma}

The proof of this lemma is deferred to Section \ref{SectionPreuveduLemme21}.

\subsection{Bad environments}\label{SubSectBadEnvironments}
Let $(\o_i)_{i\in\Z}$ be a collection of independent and identically distributed random variables satisfying
\eqref{elliptic}, \eqref{eqrecA} and \eqref{variance}.
%We introduce the potential $V$, which depends only on the environment $\o$ as follows:
We now fix $\e\in(0,1/2)$. Let $K\geq 1$. In order to transfer to our random potential $V$, with some approximations,
some results such as the ones of Cheliotis \cite{Cheliotis}, which are available for Brownian motion, but unavailable for $V$ to the extent of our knowledge,
we use the following coupling.
According to the Koml\'os--Major--Tusn\'ady strong approximation
theorem (see Koml\'os et al. \cite{KMT}), there exist (strictly) positive constants $C_1$, $C_2$ and
$C_3$, independent of $K\in\N^*$, such that, possibly in an enlarged probability
space, there exists a two-sided standard Brownian
motion $(W(t),\ t\in\R)$, such that
\begin{equation}\label{eqDefB1}
    \mathcal{B}_1(K)
:=
    \left\{\sup_{-K\leq i\leq K} \big|V(i)-\sigma W(i)\big|\leq C_1\log
K\right\}
\end{equation}
satisfies $\eta([\mathcal{B}_1(K)]^c)\leq C_2K^{-C_3}$.

Throughout the proof, we set $a:=\cfrac{1}{2}\exp\Big(\cfrac{\sqrt{5}-3}{2\e}\Big)$.
Moreover, for $u\in\R$, $\lfloor u \rfloor$ denotes the integer part of $u$.
We define for $N>2$ the events
\begin{eqnarray}
    \mathcal{B}_2(N)
& :=&
    \{\text{the number of sign changes of }b_W
\nonumber\\
& &
    \quad \quad \text{ in } [(\log N)^\e,(\log N)^{1-\e}] \text{ is at least }2\e \log_2 N +2\},
\label{eqDefB2}\\
    \mathcal{B}_3(N)
& :=&
    [A(\lfloor \e\log_2 N\rfloor, a, (\log N)^\e)]^c,
    %\\ \mathcal{B}_4(N) & := & \{x_{-1}(W,(\log N)^\e)\leq -1\},   ({\bf a\ supprimer, inutile ? })
\label{eqDefB3}
\end{eqnarray}
where $\log_2 x:=\log\log x$ for $x>1$.

\medskip\noindent
We now introduce, for every continuous process $(Z(t),\ t\geq 0)$,
\begin{eqnarray}
\underline{Z}(t)   & := & \inf\{Z(u),\quad 0\leq u\leq t\},\qquad t\geq 0,\label{eqdefXbarre}\\
d_Z(r) & := &  \inf\{t\geq 0,\quad Z(t)-\underline{Z}(t)\geq r\},\qquad  r\geq 0.\nonumber
%X^\sharp(t) & := & \sup\{X(u)-\underline{X}(u),\quad 0\leq u\leq t\},\quad t\geq 0, \nonumber\\
%d_X(r) & := &  \inf\{t\geq 0,\quad X^\sharp(t)\geq r\},\quad  r\geq 0.\nonumber
\end{eqnarray}
Then we set $W^+(t):=W(t)$ and $W^-(t):=W(-t)$ for $t\geq 0$, and consider  for $N>1$:
\begin{equation*}
\mathcal{B}_4(N):=\left\{d_{\sigma W^+}(5\log N)\leq (\log N)^4\right\}, \quad
\mathcal{B}_5(N):=\left\{d_{\sigma W^-}(5\log N)\leq (\log N)^4\right\}.
\end{equation*}
%and $\mathcal{B}_5(N):=\mathcal{B}_5^+(N)\cap \mathcal{B}_5^-(N)$.
For technical reasons, we also introduce
\begin{equation*}
\mathcal{B}_6(N):=\left\{\forall k\in\Z\cap[-\log^4 N-1, \log^4 N),\ \forall t\in[k,k+1],\quad  |W(t)-W(k)|\leq \log_2 N\right\}. \quad
\end{equation*}
This enables us to define the set $\mathcal{B}(N)$ of
{\it bad environments} as follows:
\begin{equation*}
    \mathcal{B}(N)
:=
    %\{W\in\mathcal{W}\}
    %\cap
    \mathcal{B}_1\Big[\big\lfloor (\log N)^{\frac{3-\sqrt{5}}{2C_3}+4}\big\rfloor\Big]
    \cap\bigcap_{i=2}^6 \mathcal{B}_i(N).
\end{equation*}
We now estimate the probability of bad environments with the following lemma:

\begin{lemma} \label{LemmeProbaBad}
If $\e>0$ is small enough, we have for large $N$,
\begin{equation}\label{e1p14}
\eta(\mathcal{B}(N)^c)\leq \frac{3}{(\log N)^{\frac{3-\sqrt{5}}{2}-\zeta(\e)}},
\end{equation}
where $\zeta$ is a function $(0,1/3)\to\R$ such that $\zeta(t)\to_{t\to 0} 0$ and $\zeta(t)>0$ for $t>0$ small enough, which is
defined just after \eqref{eqProbaB1c}.
\end{lemma}

\medskip
\noindent{\bf Proof:}
Denote by $k_W(e^t)$ the number of sign changes of $b_W$ in $[1,e^t]$ for $t>0$.
Cheliotis (\cite{Cheliotis} Corollary 5) proves that the laws of $k_W(e^t)/t$, $t>0$ satisfy a large
deviation principle with speed $t$ and good rate function $I$, defined by
$
    I(x)
:=
    x\log\big(2x\big(x+\sqrt{x^2+5/4}\big)\big)+3/2-\big(x+\sqrt{x^2+5/4}\big)
$
for $x> 0$, $I(x):=+\infty$ for $x<0$, and
$I(0):=(3-\sqrt{5})/2$. Hence by scaling, for $N$ large enough,
%(revoir les $\e$; et le +2 ?)
\begin{eqnarray}
\eta(\mathcal{B}_2(N)^c)
    %& = & \eta(\text{the number of sign changes of }b \text{ in } [1,(\log N)^{1-2\e}]\text{ is less than }2\e \log_2 N+1 %)\nonumber\\
& \leq  &
     \eta\big(k_W\big(e^{(1-2\e)\log_2 N}\big)\leq 3\e\log_2 N\big)\nonumber\\
& \leq & 
    \exp\left\{-\left[I\left(3\e/(1-2\e)\right)-\e\right](1-2\e)\log_2 N\right\} \nonumber\\
    % & \leq & \exp\left(-\left(I(0)-\zeta(\e)\right)\log_2 N\right)\nonumber\\
& = & 
    (\log N)^{\zeta(\e)-\frac{3-\sqrt{5}}{2}},
\label{eqProbaB1c}
\end{eqnarray}
where $\zeta(t):=I(0)-[I(3t/(1-2t))-t](1-2t)$ for $t\in(0,1/3)$.
Notice that $\zeta(t)>0$ for small $t>0$ since $0<I(u)<I(0)$ for small $u>0$. Moreover,
 $\zeta(t)\to0$ as $t\to 0$, $t>0$, since $I$ is right-continuous at $0$.
Lemma \ref{LemmaProbaStrongChangeOfSign} gives  since $1-e^{-t}\leq t$ for $t\in\R$,
for $N$ large enough so that $\lfloor \e\log_2 N\rfloor-1>0$,
\begin{equation*}
    \eta\big[\mathcal{B}_3(N)^c\big]
=
    \eta\big[A\big(\lfloor \e\log_2 N\rfloor, a, (\log N)^\e\big)\big]
\leq
    \big(1-e^{-2a}\big)^{\lfloor \e\log_2 N\rfloor-1}
\leq
    \big(2a\big)^{\lfloor \e\log_2 N\rfloor-1}.
\end{equation*}
So,  since $2a=\exp\big([\sqrt{5}-3]/(2\e)\big)\in(0,1)$,
\begin{equation*}
    \eta\big[\mathcal{B}_3(N)^c\big]
\leq
    \big(2a\big)^{\e\log_2 N-2}
=
    \big[\exp\big((3-\sqrt{5})/\e\big)\big]  (\log N)^{\frac{\sqrt{5}-3}{2}}.
\end{equation*}
Consequently, for every (fixed) $\e>0$ small enough so that $\zeta(\e)>0$, we have for $N$ large enough,
$\exp\big((3-\sqrt{5})/\e\big)\leq (\log N)^{\zeta(\e)}$ and then
\begin{equation}\label{eqProbaE9}
    \eta\big[\mathcal{B}_3(N)^c\big]
\leq
    (\log N)^{\zeta(\e)-\frac{3-\sqrt{5}}{2}}.
\end{equation}
Notice that for $r\geq 0$ and $T>0$,
\begin{equation*}
%\label{eqProbaDW}
    \eta(d_{W^+}(r)> T)
%=
%    \eta\Big(\sup_{0\leq t\leq T}W^{\sharp}_+(t)\leq r\Big)
\leq
    %\eta( W^{\sharp}_+(T)\leq r)
    \eta\left(W^+(T)-\underline{W^+}(T)\leq r\right)
=
    \eta(|W(T)|\leq r)
\leq
    2r/\sqrt{T},
\end{equation*}
since
%$(W^{\sharp}_+(t),\ 0\leq t\leq T)=_{law}(|W(t)|,\ 0\leq t\leq T)$
$W^+(T)-\underline{W^+}(T)=_{law}|W(T)|$
%$W^{\sharp}_+(T)=_{law}|W(T)|$
(see L\'évy's theorem e.g. in Revuz and Yor \cite{RevuzYor} th VI.2.3).
This gives
\begin{equation}\label{eqProbaE5}
\eta[\mathcal{B}_4(N)^c]=\eta[\mathcal{B}_5(N)^c]\leq 10/(\sigma \log N).
%\eta(\mathcal{B}_5(N)^c)\leq 20/(\sigma\log N).
\end{equation}
%$$\P(d_{W^+}(r)\geq T)=\P(\sup_{0\leq t\leq T}W^{\sharp}_+(t)<r)=
%\P(\sup_{0\leq t\leq T}|B(t)|<r)
%\leq \P(|B(T)|<r)$$
Moreover for large $N$, we get since $\sup_{0\leq t\leq 1}W(t)=_{\text{law}} |W(1)|$
and $\eta[W(1)\geq x]\leq e^{-x^2/2}$ for $x\geq 1$,
%(utile ? ie KMT vrai pour reels aussi ?) + OK, $S_t=_l L_t$ et $S_t=_l |B_t|$ donc OK.
\begin{equation}
    \eta(\mathcal{B}_6(N)^c)
\leq
    3(\log^4 N) \eta\Big(\sup_{0\leq t\leq 1}|W(t)|> \log_2 N\Big)
        %\nonumber\\
        % & \leq & 4 (\log^4 N) \eta(\sup_{0\leq t\leq 1}W(t)>  \log_2 N)
        %\nonumber\\
        % & = & 4 (\log^4 N) \eta(|W(1)|> \log_2 N)\nonumber\\
\leq
    12(\log^4 N) \exp(- (\log_2 N)^2/2)
\leq
    (\log N)^{-2}.\label{eqProbaE6}
\end{equation}

Combining this with  \eqref{eqProbaB1c}, \eqref{eqProbaE9}, \eqref{eqProbaE5}
%, $\eta(W\in\mathcal{W})=1$
and
$\eta(\mathcal{B}_1(K)^c)\leq \frac{C_2}{K^{C_3}}$
proves the lemma. \hfill$\Box$

\subsection{Random walk in a bad environment}\label{SubSectRWinBadEnvironment}

In the following lemma, we show that in a bad environment, the quenched probability that $\sum_{k=0}^n f(S_k)$
is greater than $u\leq 0$ for all $n$ between $1$ and $N$ is small:

\begin{lemma}\label{LemmeComportementdansBad}
Let $f$ be as in Theorem \ref{TheoremePersistenceAdditive}, and $u\leq 0$. For large $N$,
\begin{equation}\label{eqLemmeComportementdansBad}
\forall \o\in \mathcal{B}(N),\qquad \po\left(\forall n\in[1,N],\ \  \sum_{k=0}^n f(S_k)>u\right)\leq 4 (\log N)^{-2}.
\end{equation}
\end{lemma}

\medskip\noindent{\bf Proof of Lemma \ref{LemmeComportementdansBad}:}
We  assume that $\o\in \mathcal{B}(N)$, and we prove that in such a bad environment, there exists a time $t\in[1,N]$ such that $\sum_{k=1}^t f(S_k)\leq u$, with a large enough quenched probability.

First, define $C_4:=\sigma+\frac{3-\sqrt{5}}{2C_3}C_1+4C_1$.
Since $\o\in \mathcal{B}_6(N)\cap \mathcal{B}_1\big[\big\lfloor (\log N)^{\frac{3-\sqrt{5}}{2C_3}+4}\big\rfloor \big]$, we have
\begin{equation}\label{eqC4}
    \forall u\in[-\log^4 N, \log^4 N],
\qquad
    |V(\lfloor u\rfloor)-\sigma W(u)|\leq C_4\log_2 N.
\end{equation}
Notice that since $\o\in\mathcal{B}_3(N)$, there exists $k_N\in\{1,\dots,2\lfloor \e\log_2 N\rfloor \}$ such that
$\mathbf{h}_N:=X_{k_N}$ is an $a$-strong change of sign of $b_W$ and $b_W(\mathbf{h}_N)>0$, where the $(X_k)_k$ are the ones in
Fact \ref{FactExcessHeight0} with $c=(\log N)^\e$.
Moreover, since $\o\in \mathcal{B}_2(N)\cap \mathcal{B}_3(N)$,
$$
(\log N)^\e\leq \mathbf{h}_N=X_{k_N}< X_{k_N+1}< X_{k_N+2}\leq X_{2\lfloor \e\log_2 N\rfloor +2}\leq (\log N)^{1-\e}.
$$
%There exists $x=X_{k_N}\in[(\log N)^\e,(\log N)^{1-\e}]$, such that $x$ is an
%$a$-strong change of sign for $W$, $b(x)>0$, and $X_{k_N+1}\leq (\log N)^{1-\e}$.
%$e(T_1(x/\sigma))\geq a x/\sigma$ and $e(T_{-1}(x/\sigma))\geq a x/\sigma$ (for $W$; for $\sigma W$ we multiply by %$\sigma$...).

\begin{figure}[h]
\begin{center}
\begin{tikzpicture}
\path[use as bounding box] (3.2,17.7) rectangle (18,26.0);
%\draw[use as bounding box] (3,17) rectangle (18,27);
\includegraphics{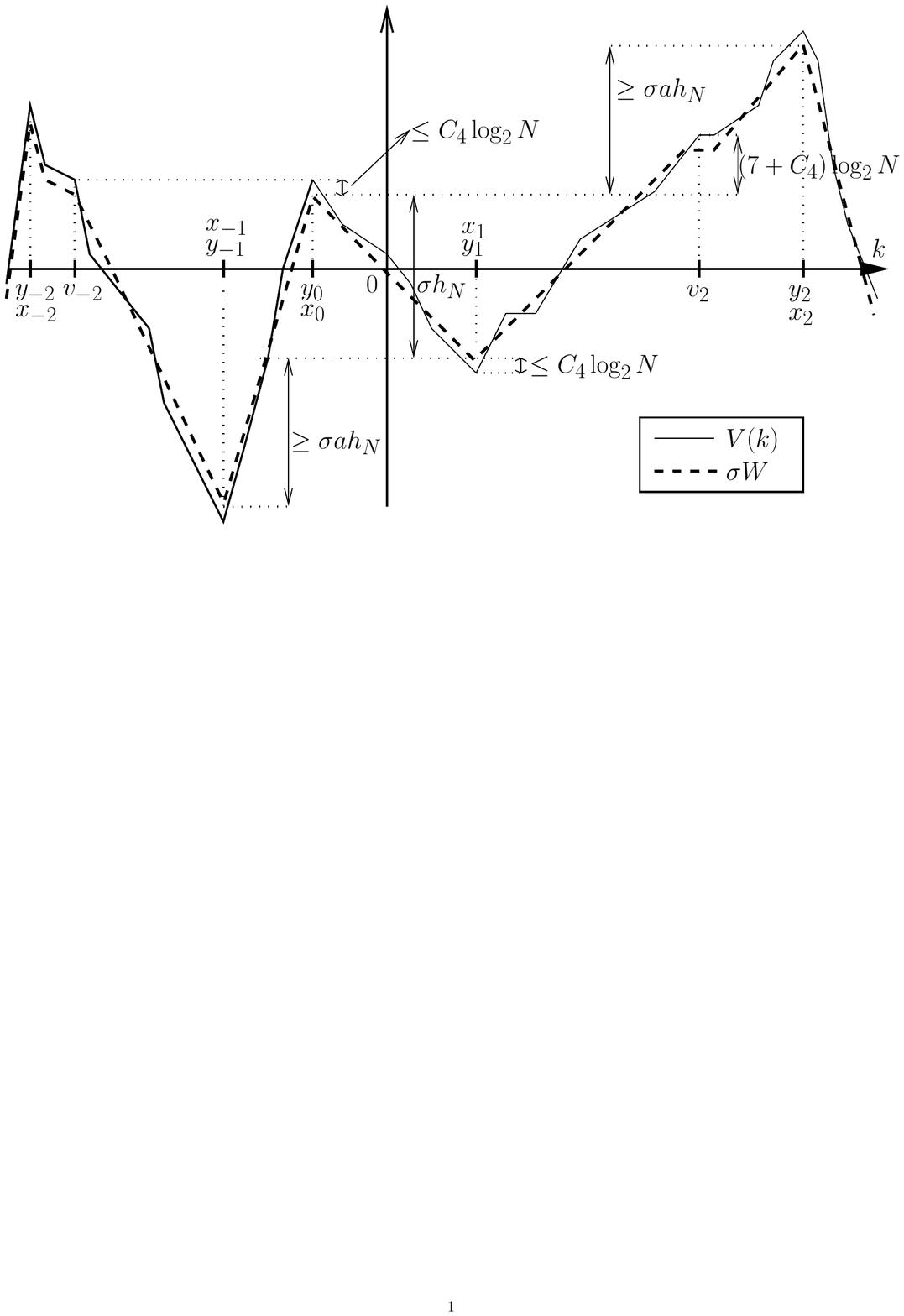}
\end{tikzpicture}
    \caption{Schema of the potential $V$ for a ``bad'' environment $\o\in \mathcal{B}(N)$ in the case $x_{-2}<v_{-2}$}
    \label{figureBad}
\end{center}
\end{figure}

To simplify the notation,  we set $x_i:=x_i(W,\mathbf{h}_N)$ and $y_i:=\lfloor x_i\rfloor$ for $i\in\{-2,\dots,2\}$.
We also define (see Figure \ref{figureBad})
%$C_4:=\sigma+\frac{3-\sqrt{5}}{2C_3}C_1+4 C_1$ and
\begin{eqnarray*}
v_{-2} & := & \max\{k\in\Z,\ k\leq y_{-1},\ V(k)\geq V(y_0)\},\\
v_2 & := & \min\{k\in\Z,\ k>y_1,\ V(k)\geq \sigma W(x_0)+(7+C_4)\log_2 N\}.
\end{eqnarray*}

\noindent

%\begin{figure}[h]
%\begin{center}
%    \scalebox{0.53}{\input{b18.tex}}
%    \caption{Schema of the potential $V$ for a ``bad'' environment $\o\in \mathcal{B}(N)$ in the case $x_{-2}<v_{-2}$}
%    \label{figureBad}
%\end{center}
%\end{figure}

\noindent
Since $b_W(\mathbf{h}_N)>0$, $x_1$ is an $\mathbf{h}_N$-minimum for $W$, and consequently $x_0$ and $x_2$ are $\mathbf{h}_N$-maxima for $W$ and $x_{-1}$ is an $\mathbf{h}_N$-minimum for $W$. Moreover, $e(T_0(\mathbf{h}_N))=0$, $e(T_1(\mathbf{h}_N))\geq a \mathbf{h}_N$ and $e(T_{-1}(\mathbf{h}_N))\geq a \mathbf{h}_N$ since $\mathbf{h}_N$ is an $a$-strong change of sign of $b_W$. Due to these properties, we get
\begin{eqnarray}
W(x_0) & = & \sup\{W(t),\ t\in[x_{-1},x_1]\}\geq 0,\label{eq55}\\
W(x_1) & = & \inf\{W(t),\ t\in[x_{0},x_2]\}\leq 0,\label{eqSupW}\\
W(x_0)-W(x_1) & = & \mathbf{h}_N,\label{ep7}\\
W(x_2)-W(x_1) & \geq & (1+a)\mathbf{h}_N,\label{eq57}\\
W(x_0)-W(x_{-1}) & \geq & (1+a)\mathbf{h}_N,\label{eq56}\\
%|V(\lfloor t\rfloor)-\sigma W(t) | & \leq & C_4\log_2 N, \qquad t\in[-(\log N)^4, (\log N)^4].
W(x_{-1}) & = & \inf\{W(t),\ t\in[x_{-2},x_0]\}< W(x_1).\label{eqp7xm2}
\end{eqnarray}
%Notice that for large $N$, $v_2\leq d_{\sigma W^+}(\sigma X_{k_N}+(7+2C_4)\log_2 N)\leq d_{\sigma W^+}(5\log N)\leq (\log N)^4 $.

The following lemma will allow us to apply \eqref{eqC4} to some $x_i$, $y_i$ and $v_i$.

\begin{lemma}\label{LemmaEncadrementxi} For $N$ large enough,
\begin{eqnarray}
\forall\o\in\mathcal{B}(N), & &   -(\log N)^4\leq v_{-2}\leq x_{-1}< x_0\leq 0< x_1<v_2< x_2\leq (\log N)^4,\nonumber\\
\forall\o\in\mathcal{B}(N), & &
v_{-2}+3\leq y_{-1}\leq y_0-3\leq -3. \label{InegsLemma28}
\end{eqnarray}
\end{lemma}

\noindent{\bf Proof:}
First, it is clear by definition that $x_{-2}<x_{-1}<x_0\leq 0<x_1<x_2$.

Moreover, $x_1$ is an $X_{k_N}$-minimum, whereas $x_1(W,X_{k_N+1})$ is an $(X_{k_N+1})$-maximum
because $X_{k_N}=\mathbf{h}_N$ and $X_{k_N+1}$ are consecutive changes of sign for $b_W$.
So $x_1\neq x_1(W,X_{k_N+1})$. Since $x_1(W,X_{k_N+1})$ is also an $X_{k_N}$-maximum, and $x_2$ is the smallest positive
$X_{k_N}$-maximum, we get $x_2\leq x_1(W,X_{k_N+1})$.
Now, if $0\leq t \leq x_1(W,X_{k_N+1})$, $W^+(t)-\underline{W^+}(t)$ is less than or equal to
$$W^+[x_1(W,X_{k_N+1})]-\underline{W^+}[x_1(W,X_{k_N+1})]\leq W[x_1(W,X_{k_N+1})]-W[x_0(W,X_{k_N+1})]=X_{k_N+1}.$$
Since $X_{k_N+1}\leq (\log N)^{1-\e}\leq (5/\sigma)\log N$ for $N$ such that $5\log N>\sigma (\log N)^{1-\e}$,
and $\o\in\mathcal{B}_4(N)$, this yields
$$0< x_2\leq x_1(W,X_{k_N+1})\leq d_{\sigma W^+}(5\log N)\leq (\log N)^4.$$

Since $v_2>y_1=\lfloor x_1\rfloor$, we have $x_1<v_2$.
Moreover, we can now apply \eqref{eqC4} to $x_2$ together with \eqref{ep7} and \eqref{eq57},
which gives
$V(y_2)\geq \sigma W(x_2)-C_4\log_2 N\geq \sigma W(x_0)+\sigma a \mathbf{h}_N-C_4\log_2 N$, which is greater than
$\sigma W(x_0)+(7+C_4)\log_2 N+2\log[(1-\e_0)/\e_0]$ uniformly on $\mathcal{B}(N)$ for $N$ large enough. This gives $v_2<y_2\leq x_2$.

Moreover, $v_{-2}\leq y_{-1}\leq x_{-1}$.
Now, similarly as before, $x_0(W,X_{k_N+2})<x_0(W,X_{k_N+1})<x_0(W,X_{k_N})$, and since all of them are $X_{k_N}$-extrema, this yields  $x_0(W,X_{k_N+2})\leq x_{-2}$.
%Now, we have $W^-(x_{-2})-\underline{W^-}(x_{-2})\leq W[x_0(W,X_{k_N+2})]-W[x_1(W,X_{k_N+2})]
%=X_{k_N+2}\leq (\log N)^{1-\e}$, which gives as previously $x_{-2}\geq d_{\sigma W^-}(5\log N)\geq -(\log N)^4$.
Now, we have
$
    W^-(-x_0(W,X_{k_N+2}))-\underline{W^-}(-x_0(W,X_{k_N+2}))
\leq
%W[x_0(W,X_{k_N+2})]-W[x_1(W,X_{k_N+2})]
H(T_0(X_{k_N+2}))
=
    X_{k_N+2}
\leq
    (\log N)^{1-\e}
$,
which gives as previously
$x_{-2}\geq x_0(W,X_{k_N+2}) \geq -d_{\sigma W^-}(5\log N)\geq -(\log N)^4$.

%{\bf I do not know if we need the inequality for $v_{-2}$}.
We already know that $x_0(W, X_{k_N+2})\leq x_{-1}< 0 < x_2\leq x_1(W, X_{k_N+1})<x_1(W, X_{k_N+2})$, which leads to
$W[x_0(W, X_{k_N+2})]\geq W(x_2)\geq W(x_0)+a \mathbf{h}_N$
since $x_0(W, X_{k_N+2})$ is an $(X_{k_N+2})$-maximum.
%with $X_{k_N+2}\geq (1+a) \mathbf{h}_N$.
Applying \eqref{eqC4} to $x_0(W, X_{k_N+2})\geq -\log^4 N$ and to $x_0\geq x_0(W, X_{k_N+2})\geq -\log^4 N$, this gives
$V(\lfloor x_0(W, X_{k_N+2})\rfloor)\geq \sigma W(x_0(W, X_{k_N+2}))-C_4\log_2 N\geq
\sigma W(x_0)+C_4\log_2 N\geq V(y_0)$ for $N$ such that $\sigma a \mathbf{h}_N \geq 2C_4\log_2 N$,
which yields  $v_{-2}\geq \lfloor x_0(W, X_{k_N+2})\rfloor \geq -(\log N)^4$.

Finally, notice that by \eqref{eqC4} and \eqref{eq56},
$$V(y_0)-V(y_{-1})\geq \sigma W(x_0)-\sigma W(x_{-1})-2C_4\log_2 N \geq \sigma (1+a)\mathbf{h}_N-2C_4\log_2 N,$$ which is, for large $N$ uniformly on $\mathcal{B(N)}$,  strictly larger
than $-3\log \e_0\geq 3\sup_{k\in\Z}|V(k)-V(k-1)|$ since $\mathbf{h}_N\geq (\log N)^\e$. This and $x_{-1}\leq x_0\leq 0$ give
the second inequality in \eqref{InegsLemma28}. The first one is obtained similarly.
\hfill$\Box$

%Consequently for large $N$, $\sigma W(x_0)+3C_4\log_2 N\leq \sigma W(x_2)$. We also notice that $W(x_2)-W(x_0(W,X_{k_N+1}))\leq X_{k_N+1}$.
%Furthermore, since $0\geq x_0\geq -d_{\sigma W^-}(5\log N)\geq -\log^4 N$, \eqref{eqC4} yields
%$$V(v_{-2})\leq \sigma W(x_0)-\log\e_0+ C_4\log_2 N\leq \sigma W(x_0(W,X_{k_N+1}))+\sigma X_{k_N+1}-C_4\log_2 N.
%$$
%Hence $W(v_{-2})-W(x_0(W,X_{k_O+1}))\leq X_{k_N+1}\leq (\log N)^{1-\e}$ which gives
%$|v_{-2}|\leq d_{\sigma W^-}(\sigma X_{k_N+1})\leq d_{\sigma W^-}(5\log N)\leq (\log N)^4$.
%%Notice that $|x_0(\sigma W^{\pm},\mathbf{h}_N)|\leq d_{\sigma W^{\pm}}(\mathbf{h}_N)$.
%Since $W(x_2)-W(x_1)\geq (1+a)\mathbf{h}_N$ and $\log_2 N=o(\mathbf{h}_N)$ as $n\to+\infty$,
%we also get $v_2\leq x_2\leq x_1(W, X_{k_N+1})\leq d_{\sigma W^+}(\sigma X_{k_N+1})\leq (\log N)^4$ for large $N$.
%Consequently, $|x_i|\leq (\log N)^4$
%for $i\in\{-1,0,1,2\}$ on $\mathcal{B}_5$.
%These properties are satisfied for every $\o\in\mathcal{B}(N)$ if $N$ is large enough.

%Hence, $|x_i|\leq (\log N)^4$
%for $i\in\{-1,0,1\}$ on $\mathcal{B}_5$, and $x_2\leq d_{\sigma W^+}(\sigma X_{k_N+1})\leq (\log N)^4$.
%Since $W(x_2)-W(x_1)\geq (1+a)\mathbf{h}_N$ and $\log_2 N=o(a \mathbf{h}_N)$ as $n\to+\infty$,
%we also get $v_2\leq x_2$ for large $N$. These properties are satisfied for every $\o\in\mathcal{B}$ is $N$ is large enough.

%and the same inequality holds for $|y_{i}|$ and $|y_{-2}|$. (utilité de $x_2$ ?)

Let
$$
    E_1
:=
    \{\tau(y_{-1})<\tau(v_2)\},
\qquad
    E_2
:=
    \big\{L((0,v_2],\tau(y_{-1})\wedge \tau(v_2))\leq (\log N)^{18+2C_4}e^{\sigma \mathbf{h}_N}\big\}.
$$
We prove the following lemma:
\begin{lemma}\label{LemmaProbasE1E2}
For large $N$,
\begin{equation*}
\forall \o\in\mathcal{B}(N),\qquad
\po(E_1^c)\leq (\log N)^{-2}, \qquad
\po(E_2^c)\leq (\log N)^{-2}.
\end{equation*}
\end{lemma}

\medskip
\noindent{\bf Proof:}
%Since $\{x_k(W,t_2), k\in\Z\}\subset  \{x_k(W,t_1), k\in\Z\}$ for $t_1<t_2$, $\mathbf{h}_N>(\log N)^\e$ and $\o\in \mathcal{B}_3$, we have
% $y_{-1}\leq x_{-1}(W,\mathbf{h}_N)\leq x_{-1}(W,(\log N)^\e)\leq -1$.
First, due to the previous lemma, $-(\log N)^4 \leq y_{-1}\leq -3$ uniformly on $\mathcal{B}(N)$ for large $N$, and equations \eqref{elliptic}, \eqref{eqProbaAtteinte}, \eqref{eqC4} and \eqref{eq55} yield
\begin{equation*}
    \po(E_1^c)
%=
%    \Big(\sum_{k=y_{-1}}^{-1}{e^{V(k)}}\Big)\Big(\sum_{k=y_{-1}}^{v_2-1}{e^{V(k)}}\Big)^{-1}
\leq
    |y_{-1}|\max_{y_{-1}\leq k \leq -1}e^{V(k)-V(v_2-1)}
\leq
    \e_0^{-1}(\log N)^{4+C_4} \exp[\sigma W(x_0)-V(v_2)]
\leq
    (\log N)^{-2},
\end{equation*}
for every $\o\in\mathcal{B}(N)$ for large $N$, which proves the first part of the lemma.

%If $v_2\geq 2$, let $z\in(0, v_2)$. Since $x_0\leq v_2\leq x_2$, equation \eqref{eqSupW} gives Since $x_0\leq v_2\leq x_2$, which yields
Thanks to $x_0\leq 0< v_2< x_2$ and to \eqref{eqSupW}, we have
$W(z)\geq W(x_1)$ for all $z\in(0, v_2)$.
Moreover, $V(k)\leq \sigma W(x_0)+(7+C_4)\log_2 N$ for every $k\in[y_{-1},v_2-1]$ by the definition of $v_2$, \eqref{eq55}, and \eqref{eqC4}.
This, Lemma \ref{LemmaEncadrementxi}, \eqref{eqProbaAtteinte},  \eqref{ep7} and  \eqref{eqC4} again give for $z\in(0, v_2)$,
\begin{equation*}
    \po^{z-1}[\tau(z)>\tau(y_{-1})]
=
    \frac{e^{V(z-1)}}{\sum_{k=y_{-1}}^{z-1}e^{V(k)}}
\geq
    \frac{\e_0e^{\sigma W(x_1)-C_4\log_2 N}}{2(\log N)^4 e^{\sigma W(x_0)+(7+C_4)\log_2 N}}
\geq
    \frac{\e_0 e^{-\sigma \mathbf{h}_N}}{2(\log N)^{11+2C_4}}.
\end{equation*}
%Applying \eqref{eqLemmaRevesz2} and observing that
Applying \eqref{eqLemmaRevesz1} and observing that
$v_2\leq (\log N)^4$, $\po[\tau(z)<\tau(y_{-1})]\leq 1$ and $y_{-1}\leq -1$,
we obtain for every $\omega\in\mathcal{B}(N)$ for large $N$,
\begin{align*}
\eo[L((0,v_2],\tau(y_{-1})\wedge \tau(v_2))]
\ & \leq \  \sum_{z=1}^{v_2-1}\frac{\po[\tau(z)<\tau(y_{-1})]}
{\o_z\po^{z+1}[\tau(z)>\tau(v_2)]+(1-\o_z)\po^{z-1}[\tau(z)>\tau(y_{-1})]}+1\\
\ & \leq \  2\e_0^{-2}(\log N)^{15+2C_4}e^{\sigma \mathbf{h}_N}+1.
\end{align*}
Using Markov's inequality, we get $\po(E_2^c)\leq (\log N)^{-2}$ for large $N$. \hfill$\Box$

\medskip

Now, let
%$\tau'(v_{-2})$ (and $\tau'(y_0)$) be the first hitting time
%of $v_{-2}$ (respectively $y_0$) after $\tau(y_{-1})$.
$T:=\inf\{k>\tau(y_{-1}),\ S_k\in\{v_{-2}, y_0-1\}\}$
%$T:=(\tau(v_{-2})\wedge \tau(y_0-1))\circ\Theta_{\tau(y_{-1})}+\tau(y_{-1})$
be the first exit time from the interval $(v_{-2}, y_0-1)$ by the random walk $S$ after $\tau(y_{-1})$.
We introduce $n_1:=\big\lfloor \frac{\e_0^2\exp(\sigma(1+a)\mathbf{h}_N)}{2(\log N)^{2+2C_4}}\big\rfloor$ and the events
$$
E_3:=\{T\geq \tau(y_{-1})+n_1\},
\qquad
E_4:=\{\tau(y_{-1})+n_1 < N\}.
$$

%$E_3':=\{L(y_{-1}, T)\geq n_1\}$ a \'eliminer.

\begin{lemma} \label{LemmaProbasE3E4}
For $N$ large enough,
\begin{equation*}
\forall \o\in\mathcal{B}(N),
\qquad
    \po(E_3^c)\leq (\log N)^{-2},
\qquad
    \po(E_4^c\cap E_1)\leq (\log N)^{-2}.
\end{equation*}
\end{lemma}

\medskip
\noindent
{\bf Proof:}
%We know (see R{\'e}v{\'e}sz \cite{Revesz} p278-279 ? V\'erifier)
%If $N$ is large enough, $\log((1-\e_0)/\e_0)<(\log N)^\e-C_4\log_2 N$,   so $v_{-2}<y_{-1}<y_0-1$ by \eqref{InegsLemma28}.
Recall that  $v_{-2}<y_{-1}<y_0-1$ on $\mathcal{B}(N)$ for $N$ large enough by \eqref{InegsLemma28}, and that $\tau(y_{-1})<\infty$ $\P$-a.s. since $(S_n)_n$ is recurrent.
We first consider
$L(y_{-1}, T)$ and notice that it is under $\po$ a geometric random variable of parameter
\begin{eqnarray*}
p_1
    & := &
    \o_{y_{-1}}\po^{y_{-1}+1}[\tau(y_{-1})>\tau(y_0-1)]+
    (1-\o_{y_{-1}})\po^{y_{-1}-1}[\tau(y_{-1})>\tau(v_{-2})]\\
& = &
    \o_{y_{-1}} e^{V(y_{-1})}\Big(\sum_{k=y_{-1}}^{y_0-2}e^{V(k)}\Big)^{-1}
    +
    (1-\o_{y_{-1}})e^{V(y_{-1}-1)}\Big(\sum_{k=v_{-2}}^{y_{-1}-1}e^{V(k)}\Big)^{-1}\\
& \leq &
    \e_0^{-2} \exp[V(y_{-1})-V(y_0)]\\
& \leq &
    \e_0^{-2}  e^{-\sigma(1+a)\mathbf{h}_N}(\log N)^{2 C_4}=:p_2,
\end{eqnarray*}
thanks to \eqref{eqProbaAtteinte} and the definition of $v_{-2}$, and where the last inequality comes from \eqref{eqC4} and \eqref{eq56}.
This ensures that for large $N$, uniformly on $\mathcal{B}(N)$ since $\mathbf{h}_N\geq (\log N)^\e$,
\begin{equation*}
\log\po[L(y_{-1}, T)\geq n_1]= (n_1-1)\log(1-p_1)\geq -2 n_1 p_1\geq -2 n_1 p_2\geq-(\log N)^{-2}.
\end{equation*}
Since $1-e^{-t}\leq t$ for $t\in\R$, this yields $\po[L(y_{-1}, T)< n_1]\leq (\log N)^{-2}$. Finally, we have
$T\geq \tau(y_{-1})+L(y_{-1}, T)$, which gives
$\po(E_3^c)\leq \po[L(y_{-1}, T)< n_1]
\leq (\log N)^{-2}$.

We now turn to $E_4$.
Notice that uniformly on $\mathcal{B}(N)$ for large $N$,
%since $-(\log N)^4\leq x_{-1}\leq v_2\leq x_2\leq(\log N)^4$ by Lemma \ref{LemmaEncadrementxi} and
thanks to Lemma \ref{LemmaEncadrementxi}, \eqref{eqC4}, \eqref{eq55}, \eqref{eqSupW}, \eqref{eqp7xm2}
and the definition of $v_2$, we have
\begin{equation}
\forall k\in[y_{-1},v_2-1],\qquad
\sigma W(x_{-1})-C_4\log_2 N\leq V(k)\leq \sigma W(x_0)+(7+C_4)\log_2 N.
\label{eqEncadrementVk}
\end{equation}
Since $H(T_0(X_{k_N}))=X_{k_N}<X_{k_{N+1}}$, $x_0$ and $x_1$ are not $(X_{k_{N+1}})$-extrema.
Hence, $[x_{-1},x_2]\subset [x_0(W,X_{k_N+1}),x_1(W,X_{k_N+1})]$, and then $W(x_2)-W(x_{-1})\leq X_{k_N+1}$.
Moreover, $\log_2 N=o(\mathbf{h}_N)$ uniformly on $\mathcal{B}(N)$ and $W(x_0)\leq W(x_2)-a \mathbf{h}_N$
by \eqref{ep7} and \eqref{eq57}, so \eqref{eqEncadrementVk} gives for large $N$,
\begin{eqnarray*}
    \max\{V(k)-V(\ell),\ y_{-1}\leq\ell\leq k\leq v_2-1\}
& \leq &
    \sigma(W(x_0)-W(x_{-1}))+(7+2C_4) \log_2 N
\\
& \leq &
    \sigma (W(x_2)-W(x_{-1}))
\leq
    \sigma X_{k_N+1}
\leq
    \sigma(\log N)^{1-\e}.
\end{eqnarray*}
%$\max_{y_{-1}\leq\ell\leq k\leq v_2}(V(k)-V(\ell))\leq \sigma (W(x_2)-W(x_{-1}))+2C_4\log_2 N
%\leq\sigma X_{k_N+1}+2C_4\log_2 N\leq 2\sigma(\log N)^{1-\e}$,
This together with \eqref{InegEsperanceZeitouni} and $|v_2-y_{-1}|\leq 2(\log N)^4$ yield
$
    \eo(\tau(y_{-1})\mathds{1}_{E_1})
\leq
    \eo[\tau(y_{-1})\wedge \tau(v_2)]
<
    \sqrt{N}
$
uniformly on $\mathcal{B}(N)$ for large $N$.
Since $\eo(n_1\mathds{1}_{E_1})<\sqrt{N}$ because $\mathbf{h}_N\leq (\log N)^{1-\e}$ on $\mathcal{B}(N)$, this yields
$\po(E_4^c\cap E_1)\leq (\log N)^{-2}$ for every $\o\in\mathcal{B}(N)$ for large $N$ by Markov's inequality.
\hfill$\Box$

%Moreover, $\sup_{t\in[v_{-2}, y_{-1}]}  W(t)\leq \sigma^{-1}(V(y_0)+C_4\log_2 N)<W(x_{-1}$
%Similarly, $\eo(T-\tau(y_{-1}))<\sqrt{N}$. This yields for $N$ large enough, by Tchebychev's inequality,
%$
%\po(E_4^c\cap E_1)\leq (\log N)^{-2}.
%$.

We now consider $f$ satisfying the hypotheses of Theorem \ref{TheoremePersistenceAdditive}.
For every $\o\in\mathcal{B}(N)$, we have on $E_1\cap E_2$ and then on $E_5:=\cap_{i=1}^4 E_i$,  since  $f(x)\leq 0$ for every $x\leq 0$,
\begin{eqnarray}\label{e1p20}
\sum_{k=0}^{\tau(y_{-1})-1} f(S_k)
=
\sum_{x=-\infty}^{v_2-1} f(x) L(x,\tau(y_{-1})\wedge \tau(v_2)-1)
& \leq &
\left[\max_{k\in(0, v_2]} f(k)\right]  L((0,v_2],\tau(y_{-1})\wedge \tau(v_2))\nonumber\\
& \leq & \left[\max_{k\in(0, v_2]} f(k)\right] (\log N)^{18+2C_4}e^{\sigma \mathbf{h}_N}.~
 \end{eqnarray}
For every $\Delta\subset \Z$ and $0\leq s\leq t$, we define
$L(\Delta, s\rightsquigarrow t):=\sum_{k=s}^t \mathds{1}_{\{S_k\in \Delta\}}$, which is the number of visits of
$(S_n)_{n\in\N}$ to the set $\Delta$ between times $s$ and $t$.
%Notice that  $S_k\leq -1$ and then $f(S_k)\leq -1$ for every $\tau(y_{-1})\leq k \leq T$.
%Moreover, $y_{-1}\leq -1$ then $f(y_{-1})\leq -1$, which gives
%\begin{eqnarray}\label{e2p20}
%\sum_{k=\tau(y_{-1})}^T f(S_k) & \leq &
%-L(y_{-1}, \tau(y_{-1})\rightsquigarrow T)\leq -n_1
%=-\e_0\frac{e^{\sigma(\mathbf{h}_N+a \mathbf{h}_N)}}{2(\log N)^{2+2C_4}}.
%\end{eqnarray}

For every $\omega\in\mathcal{B}(N)$ and each integer $k\in[\tau(y_{-1}), \tau(y_{-1})+n_1]$, we have $\tau(y_{-1})\leq k \leq T$ on $E_3$, so $S_k\leq y_0-1\leq -1$,
thus $f(S_k)\leq -1$. As a consequence on $E_5$ for large $N$,
\begin{equation}\label{e2p20}
\sum_{k=\tau(y_{-1})}^{\tau(y_{-1})+n_1} f(S_k)
% \leq -L(y_{-1}, \tau(y_{-1})\rightsquigarrow T)
\leq -n_1-1
\leq -\e_0^2\frac{\exp[\sigma \mathbf{h}_N+\sigma a(\log N)^\e]}{2(\log N)^{2+2C_4}},
\end{equation}
since $ \mathbf{h}_N\geq(\log N)^\e$.
Combining \eqref{e1p20}, \eqref{e2p20}, and
$\max_{k\in(0, v_2]} f(k)\leq \max_{k\in(0, (\log N)^4]} f(k)\leq\exp((\log N)^{\e/2})$ for large $N$,
we get $\sum_{k=0}^{\tau(y_{-1})+n_1} f(S_k) \leq u$ on $E_5$ for every $\o\in \mathcal{B}(N)$ for large $N$.
Moreover, $1\leq \tau(y_{-1})+n_1\leq N$ on $E_5$,
hence for large $N$, for every $\o\in \mathcal{B}(N)$, we have $E_5\subset \{\exists n\in[1,N],\ \sum_{k=0}^n f(S_k)\leq u\}$.
Consequently, the left hand side of \eqref{eqLemmeComportementdansBad} is less than
$\po(E_5^c)\leq 4 (\log N)^{-2}$ for every $\o\in\mathcal{B}(N)$ for large $N$
by Lemmas \ref{LemmaProbasE1E2} and \ref{LemmaProbasE3E4},
which proves Lemma \ref{LemmeComportementdansBad}.
\hfill$\Box$

\medskip
Finally, integrating \eqref{eqLemmeComportementdansBad} on the set of bad environments $\mathcal{B}(N)$ gives by Lemma \ref{LemmeProbaBad}:
\begin{align*}
    \P\left(\forall n\in[1,N],\ \sum_{k=0}^n f(S_k)>u\right)
\ & \leq \  \int_{\mathcal{B}(N)}
    \po\left(\forall n\in[1,N],\ \sum_{k=0}^n f(S_k)>u\right) \eta(\text{d}\o)+\eta(\mathcal{B}(N)^c)\\
\ & \leq \
    4(\log N)^{-2}+ \frac{3}{(\log N)^{\frac{3-\sqrt{5}}{2}-\zeta(\e)}}
\leq
\frac{4}{(\log N)^{\frac{3-\sqrt{5}}{2}-\zeta(\e)}}
% \frac{1}{(\log N)^{\frac{3-\sqrt{5}}{2}+o(1)}}
\end{align*}
for large $N$. Now, let $\e\to 0$, so $\zeta(\e)\to 0$. This gives the upper bound in Theorem \ref{TheoremePersistenceAdditive}.
\hfill$\Box$

%$\mathcal{B}_2$

%Now, let $E_2:=\{\tau_{y_{-1}}\wedge \tau_{v_2}\leq \}$. We have,
%\begin{equation*}
%\eo(\tau_{y_{-1}}\wedge \tau_{v_2})\leq \frac{(1-\e_0)(2\log^4 N)^2}{\e_0}e^{V(y_2)-V(y_0)}
%%\leq \frac{4(1-\e_0)\log^8 N}{\e_0} e^{2(\log N)^{1-\e}}(\log N)^{C_4}
%\end{equation*}
%which is $\leq \frac{N}{3\log^2 N}$. Hence by Tchebychev, $\po(E_2)\leq (\log N)^{-2}$.

%%%%%%%%%%%%%%%%%%%%%%%%%%%%%%%%%%%%%%%%%%%%%%%%%%%%%%%%%%%
%                                                         %
%              GOOD    ENVIRONMENTS                      %
%                                                         %
%%%%%%%%%%%%%%%%%%%%%%%%%%%%%%%%%%%%%%%%%%%%%%%%%%%%%%%%%%%

%\mysection{Good environments}\label{good}

\mysection{Proof of the lower bound}\label{good}

\subsection{Sketch of the proof, and organization of this proof}\label{SubSectSketchProofPart4}
We give in this subsection some non-rigorous heuristics, for which we invite the reader to look at Figure \ref{figure1};
 everything will be proved in details in the next subsections.

Let $N\geq 2$.
We build in Subsection \ref{SubSectGoodEnvironments} a set $\mathcal{G}(N)$ of "good environments".
We would like that uniformly  on these good environments $\o\in\mathcal{G}(N)$,
$\sum_{k=0}^n f(S_k)>0$ for all $1\leq n\leq N$ with large quenched probability (see Lemma \ref{LemmaRWinGoodEnv}).
To this aim, we first require that the potential $V$ of such good environments decreases quickly between $0$
and $\e\log_2 N$ and then remains low up to some random $\theta_0$,
which is the smallest $k>0$ such that $V(k)\leq -5h(N)$ ($h(N)$ being defined in \eqref{eqDefHdeN} below).
We then make a coupling between the potential outside this interval $[0,\theta_0]$, called $\VV$ and defined in \eqref{eqDefVChapeau},
and a two-sided Brownian motion $W$ (see \eqref{eqCouplingVchapeau} below).
We then require that $b_{\sigma W}(x)>0$ for all $1\leq x \leq 5\log N$, and add some technical conditions.
Such environments are called {\it good environments} $\o\in\mathcal{G}(N)$.
A schema of the  potential of a good environment is given in Figure \ref{figure1}.

We then show in Subsection \ref{SubSectProbaGoodEnv} that loosely speaking, the probability of the set of good environments is
$
    \eta[\mathcal{G}(N)]
\geq
    1/(\log N)^{\frac{3-\sqrt{5}}{2}+o(1)}.
$

Finally, we study in Subsection \ref{SectionLower} a random walk $(S_k)_k$ in a good environment $\o\in\mathcal{G}(N)$.
We introduce the location
$
    \theta_i\approx\inf\{k\geq \theta_0,\ V(k)-\inf_{0\leq \ell\leq k} V(\ell)\geq i h(N)\}
$,
$i\geq 1$, which is approximatively the first location
where there is an increase of at least $i h(N)$ for the potential $V$ restricted to $[\theta_0, +\infty)$
(see Figure \ref{figure1}, and \eqref{eqDefthetai} below).
We first show in Lemma \ref{LemmaLowerBoundProbaFpouri1} that, because the potential $V$ decreases quickly in
$[0, \e \log_2 N]$ and remains low up to $\theta_1$ with $V(\theta_1)$ much lower than $0$,
with a large quenched probability
the random walk $(S_k)_k$ goes to $\theta_1$ before going to $-1$, and then
$
    \sum_{k=0}^n f(S_k)
\geq
    f(S_1)
=
    f(1)
>
    0
$
for all $1\leq n\leq \tau(\theta_1)$.
Moreover we prove that
$
    \sum_{k=0}^{\tau(\theta_1)}f(S_k)
%\geq
%    \tau(\theta_1)
\geq
    L(m_1, \tau(-1)\wedge\tau(\theta_1))
\geq
    e^{h(N)}/[2(\log N)^{\nu}]
$
for some $\nu>0$
with large quenched probability,
that is, the sum of $f(S_k)$ has accumulated some large positive quantity at time $\tau(\theta_1)$.

We then prove by induction in Lemma \ref{LemmaLowerBoundProbaFpouriGeneral} (see also \eqref{eqDefFi})
that for every $i\geq 1$ such that $ i h(N)\leq 4 \log N$,
 with large quenched probability uniformly on all good environments $\o\in\mathcal{G}(N)$,
$\sum_{k=0}^n f(S_k)>0$ for all $1\leq n\leq \tau(\theta_i)$,
and the sum of $f(S_k)$ has accumulated some large positive quantity at time $\tau(\theta_i)$,
that is, $\sum_{k=0}^{\tau(\theta_i)} f(S_k) \geq e^{i h(N)}/[2(\log N)^{\nu}]$.

Assume that this is true for such an $i$, and fix a good environment $\omega\in\mathcal{G}(N)$.
Loosely speaking,
%$x_1(\sigma W, i h(N))$ is an $i h(N)$-minimum for $\sigma W$,
since $b_{\sigma W}[i h(N)]>0$ and $b_{\sigma W}[(i+1) h(N)]>0$,
the deepest location (in terms of potential)
that $(S_k)_k$ can visit with large quenched probability between times $\tau(\theta_i)$ and $\tau(\theta_{i+1})$ is
$m_{i+1}\approx \theta_0+x_1(\sigma W, (i+1) h(N))>\theta_0>0$ (see \eqref{eqDefmi} and Figure \ref{figure1}).
Moreover, our hypotheses for $V(x)$, $0\leq x \leq \theta_0$ have "lowered" the potential $V$ in $[\theta_0, +\infty)$
compared to the potential $V$ in $\Z_-^*$.
In particular, the potential $V(x)$ for locations $x<0$ that the random walk $(S_k)_k$ may visit between $\tau(\theta_i)$
and $\tau(\theta_{i+1})$, that is, $x\in[x_0(\sigma W, (i+1) h(N)),-1]$, satisfy by definition of $x_1(\sigma W,.)$,
$$
    V(x)
\approx
    \sigma W(x)
\geq
    \sigma W[x_1(\sigma W, (i+1) h(N))]
\approx
    V(m_{i+1})-V(\theta_0)
\approx
    V(m_{i+1})+5h(N).
$$
%Due to our hypotheses that
%%$V$ decreases in $[0,\theta_0]$ and since
%$V(\theta_0)\leq - 5 h(N)$ is quite low,
%$V(m_{i+1})\approx V(\theta_0)+\sigma W[x_1(\sigma W, (i+1) h(N))]$
%is much lower than $\sigma W[x_1(\sigma W, (i+1) h(N))]$,
%which is itself lower than $\sigma W\approx V$ in the negative locations
%that $(S_k)_k$ can visit between times $\tau(\theta_i)$ and $\tau(\theta_{i+1})$, which are included in $[x_0(\sigma W, (i+1) h(N)),-1]$.
Hence, $V(m_{i+1})$ is much lower than the potential $V(x)$ in the negative locations $x$ the random walk $(S_k)_k$ may visit
between times $\tau(\theta_i)$ and $\tau(\theta_{i+1})$, so the random walk can go to these negative locations, where $f<0$,
but the total amount of time it spends there is small, with large quenched probability
(this is proved in details in the second step of the proof of Lemma \ref{LemmaLowerBoundProbaFpouriGeneral}).

Consequently, $\big|\sum_{k=\tau(\theta_i)+1}^{\tau(\theta_{i+1})} f(S_k)\un_{f(S_k)<0}\big|$
is very small compared to
the quite large (positive) sum
$
    \sum_{k=0}^{\tau(\theta_i)} f(S_k)
\geq
    e^{i h(N)}/[2(\log N)^{\nu}]
$
already accumulated by induction hypothesis.
%Since we have already accumulated a quite large positive sum
%$\sum_{k=0}^{\tau(\theta_i)} f(S_k) \geq e^{i h(N)}/(\log N)^{o(1)}$ by induction hypothesis,
%$f(x)\geq 0$ for $x\geq 0$,
This allows us to prove that
$\sum_{k=0}^n f(S_k)>0$ for all $\tau(\theta_i)< n\leq \tau(\theta_{i+1})$ (recall that $f(x)\geq 0$ for $x\geq 0$).
Finally we prove (in the third step) that $(S_k)_k$ spends a large amount of time in the deepest location $m_{i+1}$
between times $\tau(\theta_i)$ and $\tau(\theta_{i+1})$.
This leads to  $\sum_{k=0}^{\tau(\theta_{i+1})} f(S_k) \geq e^{(i+1) h(N)}/[2(\log N)^{\nu}]$
with large quenched probability, which ends the induction.
Since we can choose $i$ so large that $\tau(\theta_i)\geq N$ with large probability,
this leads to the lower bound of Theorem \ref{TheoremePersistenceAdditive}.

\subsection{Definition of the set $\mathcal{G}(N)$ of good environments}\label{SubSectGoodEnvironments}

We consider a collection $(\o_i)_{i\in\Z}$ of independent and identically distributed random variables, satisfying
\eqref{elliptic}, \eqref{eqrecA} and \eqref{variance}.

We notice that due to \eqref{eqrecA} and \eqref{variance}, there exist $\gamma>0$ and $\delta>0$ such that $\eta(-2\delta\leq \log\frac{1-\o_0}{\o_0}\leq -\delta)=:e^{-\gamma}>0$.
We fix $\e>0$ such that $\e\d/4< 4$. Let $N\in\N$ such that  $N\geq 3$.
In the spirit of Devulder \cite{Devulder},
we first define
\begin{equation*}
\mathcal{G}_1(N):=\left\{\forall k\in\{1,\dots,\lfloor \e \log_2 N\rfloor\},\quad -2\delta \leq \log\frac{1-\o_k}{\o_k}\leq -\delta\right\},
\end{equation*}
and we introduce
\begin{eqnarray}
    h(N)
& := &
    (\log N)^{\e\delta/32},
    %\geq (22+5C_4)\log_2 N+\log \e_0,
\label{eqDefHdeN}
\\
    \t_0
& := &
    \inf\{k\geq \lfloor \e\log_2 N\rfloor, \quad V(k)\leq -5 h(N)\},
\label{eqdeft0}\\
        %\tt_0 & := & \sup\{k\leq -\e\log_2 N, \quad V(k)\geq 5 h(N)\}.\label{eqdeftt0}
    \mathcal{G}_2(N)
& := &
    \big\{\forall k\in\{\lfloor \e\log_2 N\rfloor,\dots,\t_0\},\quad V(k)\leq -(\d\e/2)\log_2 N\big\},
\nonumber\\
    \mathcal{G}_3(N)
& := &
    \big\{\t_0\leq \lfloor\e\log_2 N\rfloor+(\log N)^{\e\delta/4}\big\}.
\nonumber
\end{eqnarray}
%\begin{eqnarray*}
%\mathcal{G}_2^+   & := & \{\forall k\in\{\e\log_2 N,\dots,\t_0\},\quad V(k)\leq -\frac{\d\e}{2}\log_2 N\},\\
%\mathcal{G}_2^- & := & \{\forall k\in\{\tt_0,\dots,-\e\log_2 N\},\quad V(k)\geq \frac{\d\e}{2}\log_2 N\},\\
%\mathcal{G}_3^+   & := & \{\t_0\leq \e\log_2 N+(\log N)^{\e\delta/4}\},\\
%\mathcal{G}_3^- & := & \{\tt_0\geq -\e\log_2 N-(\log N)^{\e\delta/4}\}.
%\end{eqnarray*}
%We also set $\mathcal{G}_2:=\mathcal{G}_2^+\cap \mathcal{G}_2^-$ and $\mathcal{G}_3:=\mathcal{G}_3^+\cap \mathcal{G}_3^-$.
We also set
\begin{equation}\label{eqDefVChapeau}
    \VV(i)
:=
\left\{
\begin{array}{ll}
    V(i+\t_0)-V(\t_0) & \text{ if } i\geq 0,\\
    V(i) & \text{ if } i<0.\\
\end{array}
\right.
\end{equation}
%$$\VV(i):=V(i+\t_0)-V(\t_0), i\geq 0$$, and
%$\VV(i):=V(i)$ for $i<0$.
By the strong Markov property, $\VV$ has the same law as $V$ and is independent of
$(V(i),\ 0\leq i \leq \t_0)$.
%$(V(i), \t_0^-\leq i \leq \t_0)$.
Let $K\geq 1$. As before, according to the Koml\'os--Major--Tusn\'ady strong approximation
theorem (see Koml\'os et al. \cite{KMT}), possibly in an enlarged probability
space, there exists
%(strictly) positive constants $C_1$, $C_2$ and $C_3$, independent of $K$, and
a standard two-sided Brownian
motion $(W(t),\ t\in\R)$ such that the set
\begin{equation}\label{eqCouplingVchapeau}
    \mathcal{G}_4(K)
:=
    \left\{\sup_{-K\leq i\leq K} \Big|\VV(i)-\sigma W(i)\Big|\leq C_1\log K\right\}
\end{equation}
satisfies $\eta(\mathcal{G}_4(K)^c)\leq C_2/K^{C_3}$. Moreover, we can choose $(W(t), \ t\in\R)$
so that it is independent of $(V(i),\ 0\leq i \leq \t_0)$ since $\VV$ is independent of $(V(i),\ 0\leq i \leq \t_0)$.
In the following, we take $K=(\log N)^{\frac{3-\sqrt{5}}{2C_3}+4}$.
We introduce
\begin{eqnarray*}
    \mathcal{G}_5(N)
& := & \left\{\max\{d_{\sigma W^+}(5\log N), d_{\sigma W^-}(5\log N),d_{-\sigma W^-}(5\log N)\}
    \leq (\log N)^4\right\},\\
\mathcal{G}_7(N)
& := &
    \{\forall x\in[1/\sigma, 5(\log N)/\sigma], \quad b_W(x)>0 \},
\end{eqnarray*}
and define $\mathcal{G}_6(N)$ by the same formula as $\mathcal{B}_6(N)$.
%$\mathcal{G}_5(N)$ and $\mathcal{G}_6(N)$ by the same formulas as $\mathcal{B}_5(N)$ and $\mathcal{B}_6(N)$
%respectively.
%and $\mathcal{B}_5(N):=\mathcal{B}_5^+(N)\cap \mathcal{B}_5^-(N)$.
%We also introduce
%\begin{equation*}
%\end{equation*}

%Remettre ici le dessin Good environment06.pdf

\noindent
\begin{figure}[h]
\begin{center}
\begin{tikzpicture}
\path[use as bounding box] (3,17) rectangle (18,27);
%\draw[use as bounding box] (3,17) rectangle (18,27);
\includegraphics{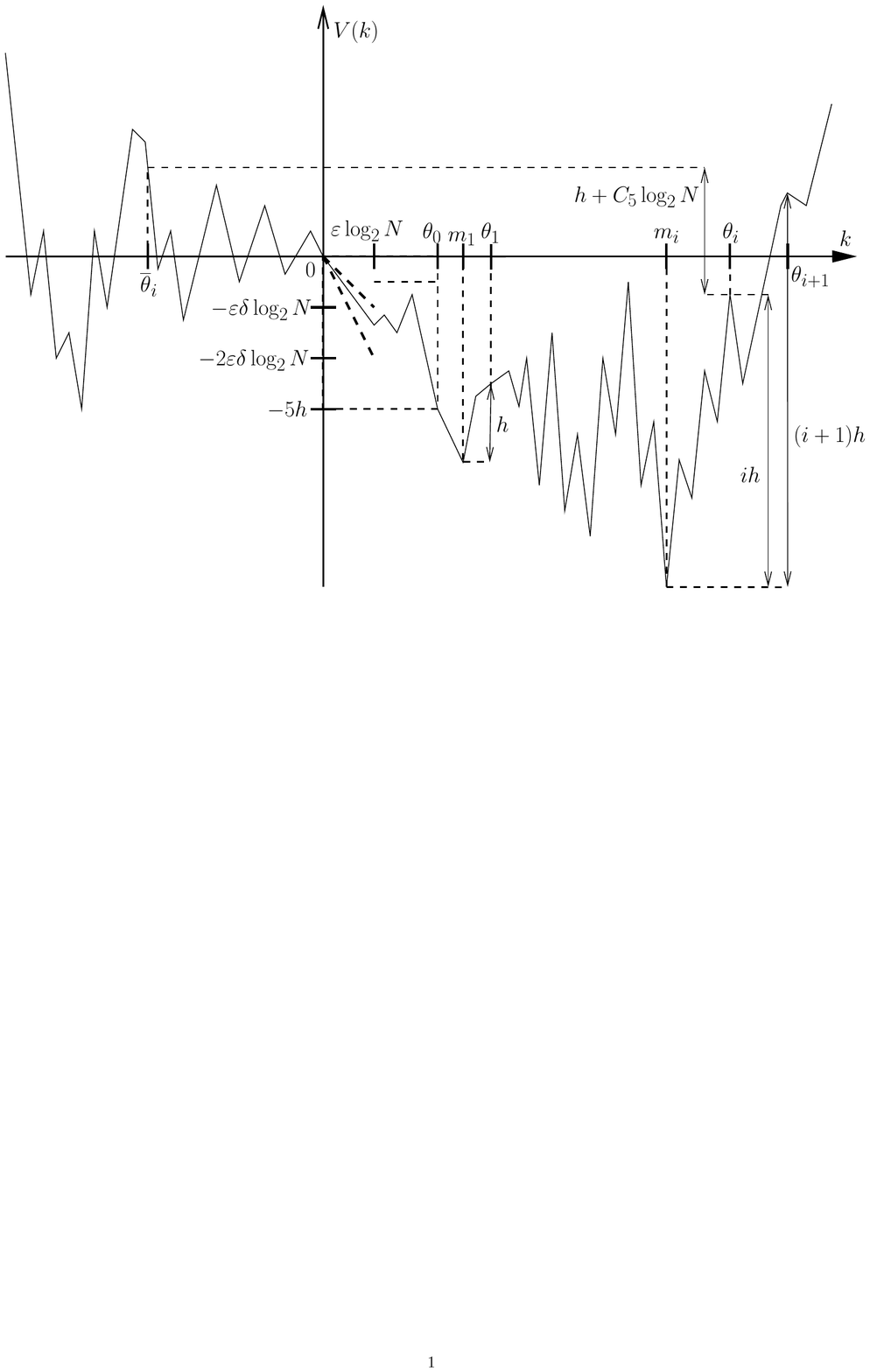}
\end{tikzpicture}
\caption{Schema of the potential $V$ for a ``good'' environment $\o\in \mathcal{G}(N)$ in the case $m_i=m_{i+1}$, where $h$ denotes $h(N)$.}
\label{figure1}
\end{center}
\end{figure}

%\begin{figure}[h]
%\begin{center}
%\scalebox{0.435}{\input{g17.tex}}
%\caption{Schema of the potential $V$ for a ``good'' environment $\o\in \mathcal{G}(N)$ in the case $m_i=m_{i+1}$, where $h$ denotes $h(N)$.}
%\label{figure1}
%\end{center}
%\end{figure}

\noindent
We can now define the set $\mathcal{G}(N)$ of {\it good environments} as follows (see Figure \ref{figure1}):
\begin{equation*}
\mathcal{G}(N):=\mathcal{G}_4\Big[(\log N)^{\frac{3-\sqrt{5}}{2C_3}+4}\Big]\cap
\bigcap_{i=1,\dots,7,\, i\neq 4 }\mathcal{G}_i(N).
\end{equation*}
When no confusion is possible we write $\mathcal{G}$ instead of $\mathcal{G}(N)$ and $\mathcal{G}_i$ instead of $\mathcal{G}_i(N)$, $i\neq 4$.

\subsection{Probability of the set $\mathcal{G}(N)$ of good environments}\label{SubSectProbaGoodEnv}

\begin{lemma} \label{LemmaProbaGood}
We have for large $N$,
\begin{equation}\label{eqprobaBonsEnvironnements}
    \eta(\mathcal{G}(N))
\geq
    \frac{c_1\e\log_2 N}{(\log N)^{\frac{3-\sqrt{5}}{2}+\e(\gamma+\delta/32)}}.
\end{equation}
\end{lemma}
\bigskip

\noindent{\bf Proof:}
First, observe that
\begin{equation*}
\eta(\mathcal{G}_1)\geq (e^{-\gamma})^{\e\log_2 N}=(\log N)^{-\e\gamma}.
\end{equation*}
We now prove that
\begin{equation}\label{eqProbaE2E3}
\eta(\mathcal{G}_2\cap \mathcal{G}_3 \mid \mathcal{G}_1)\geq \frac{\delta\e \log_2 N}{40h(N)}
\end{equation}
for large $N$. To this aim, we define $A:=\log((1-\e_0)/\e_0)$, so $|V(k+1)-V(k)|\leq A$ a.s.
for every $k\in\Z$ thanks to \eqref{elliptic}.
For $a\in\R$ and $b\in\R$ such that $a<0<b$, let $T_{a,b}:=\inf\{k\geq 0, \ V(k)\notin(a,b)\}<\infty$ a.s.
We recall that thanks to the optimal stopping theorem, $\eta[V(T_{a,b})<0]\geq b/(b-a+A)$ (see e.g. Zindy \cite{Zindy} Lemma 2.1 and apply it to $-V$).
In particular, we get on $\mathcal{G}_1$ uniformly for $N$ large enough,
$$
\eta[\mathcal{G}_2\mid V(\lfloor \e\log_2 N\rfloor)]\geq \d\e\log_2 N/(20h(N)),
$$
which yields  $\eta(\mathcal{G}_2| \mathcal{G}_1)\geq \d\e\log_2 N/(20h(N))$.
%First, $(V(k), \ k\geq 0)$ is a martingale for the natural filtration. For $a\in\R$ and $b\in\R$, with $a<0<b$,
%$T_V:=\inf\{k\geq 0, \ V(k)\notin(a,b)\}<\infty$ p.s. and $T_V$ is a stopping time.
%Hence, $\E(V(T_V))=0$ by the stopping theorem for bounded martingales. Define $A:=\log(\e_0/(1-\e_0))$ and observe that $V(k+1)-V(k)\in[-A,A]$ for %every $k\in\N$.
%This gives $0=\E(V(T_V))\geq b\P(V(T_V)>0)+(a-A)\P(V(T_V)<0)$, then $\P(V(T_V)<0)\geq b/(b-a+A)$.
%This gives $0=\E(V(T_V))\leq (a+A)\P(V(T_V)>0)+(-b)\P(V(T_V)<0)$, and  then $\P(V(T_V)>0)\geq b/(a+b+A)$. Similarly, we
%have $\P(V(T_V)>0)\leq (b+A)/(a+b+A)$. In particular, we get for $N$ large enough,
%$$\eta(\mathcal{G}_2)\geq \d\e\log_2 N/(20h(N)).$$
Moreover,  we have on $\mathcal{G}_1$ by the Markov property
\begin{eqnarray*}
    \eta(\mathcal{G}_2\cap \mathcal{G}_3^c|V(\lfloor \e\log_2 N\rfloor))
        %& \leq & \eta(V(k)\in[-5h(N),2\delta\e\log_2 N],\ \forall k\in[0,(\log N)^{\delta\e/4}])\\
& \leq &
    \eta\big(V\big(\big\lfloor(\log N)^{\delta\e/4}\big\rfloor\big)\in[-5h(N),2\delta\e\log_2 N]\big)
\\
& \leq &
    \eta\left(\frac{\big|V\big(\big\lfloor(\log N)^{\delta\e/4}\big\rfloor\big)\big|}{\sigma\sqrt{\lfloor(\log N)^{\delta\e/4}\rfloor}}
            \leq
            \frac{5h(N)}{\sigma\sqrt{\lfloor(\log N)^{\delta\e/4}\rfloor}}
        \right)
\end{eqnarray*}
for $N$ large enough. By Berry-Esseen, we get with $Y=_{\text{law}}\mathcal{N}(0,1)$,
\begin{eqnarray*}
\eta(\mathcal{G}_2\cap \mathcal{G}_3^c|V(\lfloor \e\log_2 N\rfloor))
& \leq & \eta\left(|Y|\leq
\frac{5h(N)}{\sigma\sqrt{\lfloor(\log N)^{\delta\e/4}\rfloor}}\right)
+\frac{c_2}{\sqrt{\lfloor(\log N)^{\delta\e/4}\rfloor}}\\
& \leq & \frac{11h(N)}{\sigma\sqrt{2\pi}(\log N)^{\delta\e/8}}+\frac{2c_2}{(\log N)^{\delta\e/8}}=
o(\eta(\mathcal{G}_2\mid \mathcal{G}_1))
\end{eqnarray*}
as $N\to +\infty$. Consequently $\eta(\mathcal{G}_2\cap \mathcal{G}_3^c\mid \mathcal{G}_1)=o(\eta(\mathcal{G}_2\mid \mathcal{G}_1))$, which gives \eqref{eqProbaE2E3} for large $N$.
%$$\P(\mathcal{G}_2\cap \mathcal{G}_3)=\P(\mathcal{G}_2)-\P(\mathcal{G}_2\cap \mathcal{G}_3^c)\geq \P(\mathcal{G}_2)-\P(\mathcal{G}_3^c).$$

\noindent
Since $W$ and $\VV$ are independent of $(V(i),\ 0\leq i \leq \t_0)$, we get
\begin{equation}\label{eqBonsEnvironnementsMarkov}
\eta(\mathcal{G}(N))=\eta(\mathcal{G}_1\cap \mathcal{G}_2\cap \mathcal{G}_3)\eta(\mathcal{G}_8)
\geq \frac{c_3 \e\log_2 N}{(\log N)^{\e \gamma}h(N)}\eta(\mathcal{G}_8)
\end{equation}
where $\mathcal{G}_8(N):=\mathcal{G}_4\big[(\log N)^{\frac{3-\sqrt{5}}{2C_3}+4}\big]\cap\bigcap_{i=5\dots 7}\mathcal{G}_i(N)$.
We now need the following result:

\begin{theo} \label{PropCheliotis}
(Cheliotis \cite{Cheliotis}, Corollary 1)
\begin{equation*}
\eta(\{(t\mapsto b_W(t)) \text{ keeps the same sign in }[1,x]\})/x^{(-3+\sqrt{5})/{2}}\longrightarrow_{x\to+\infty} 1/2 +7\sqrt{5}/30=:c_4.
\end{equation*}
\end{theo}

\noindent
Hence, $\eta(\mathcal{G}_7)\sim_{N\to+\infty} c_4/[2(5\log N)^{(3-\sqrt{5})/2}]$,
due to the scaling property of $b_W$, that is for fixed $r>0$,
$$
\big(b_W(r x), \ x> 0\big) =_{\text{law}}  \big(r^2 b_W( x), \ x> 0\big).
$$
Moreover,
$\eta[\mathcal{G}_5^c]\leq 30/(\sigma \log N)$
by \eqref{eqProbaE5},
$\eta[(\mathcal{G}_4(K))^c]\leq C_2/K^{C_3}$,
and
$\eta[\mathcal{G}_6(N)^c] \leq (\log N)^{-2}$
by \eqref{eqProbaE6}, so
\begin{eqnarray*}
    \eta(\mathcal{G}_8^c)
%\leq
%    \eta(\mathcal{G}_4[(\log N)^{\frac{3-\sqrt{5}}{2C_3}+4}]^c)
%    +\sum_{i=5\dots 7}\eta(\mathcal{G}_i(N)^c)
\leq
    1-  c_5/(\log N)^{\frac{3-\sqrt{5}}{2}}
\end{eqnarray*}
for $N$ large enough for some $c_5>0$, since $\frac{3-\sqrt{5}}{2}<1$. Hence,
$\eta(\mathcal{G}_8)\geq c_5/(\log N)^{\frac{3-\sqrt{5}}{2}}$ for large $N$.
This, combined with \eqref{eqBonsEnvironnementsMarkov}, gives \eqref{eqprobaBonsEnvironnements}.
\hfill$\Box$

%%%%%%%%%%%%%%%%%%%%%%%%%%%%%%%%%%%%%%%%%%%%%%%%%%%%%%%%%%%
%                                                         %
%              PROOF OF THE LOWER BOUND                   %
%                                                         %
%%%%%%%%%%%%%%%%%%%%%%%%%%%%%%%%%%%%%%%%%%%%%%%%%%%%%%%%%%%

\subsection{Random walk in a good environment}\label{SectionLower}

In this subsection, we prove the following lemma, and then the lower bound of
Theorem \ref{TheoremePersistenceAdditive}. Notice that we just have to consider the case $u=0$.
In all the rest of this section, the function $f$ satisfies the hypotheses of Theorem \ref{TheoremePersistenceAdditive}.

\begin{lemma} \label{LemmaRWinGoodEnv}
There exists a constant $c_6>0$ such that for $N$ large enough,
\begin{equation}\label{EqLemmaRWinGoodEnv}
\forall\o\in\mathcal{G}(N),
\qquad \po\left(\sum_{k=0}^n f(S_k)>0\quad \forall 1\leq n\leq N\right)
\geq c_6.
\end{equation}
\end{lemma}

\bigskip
%\noindent {\bf Proof:}
Before proving this lemma, we introduce some more notation.
We consider $N\geq 3$ and a good environment $\o\in \mathcal{G}(N)$.
We introduce for $i\in\N^*$ (see Figure \ref{figure1}),
\begin{eqnarray}
    t_i
& := &
    \inf\{t>0, \ \sigma W(t)-\sigma \underline{W}(t)\geq i h(N)\}=d_{\sigma W^+}(i h(N)),
\nonumber\\
    \t_i
& := &
    \lfloor t_i\rfloor+\t_0,
\label{eqDefthetai}\\
        %\t_i & := & \inf\{k>\t_0,\quad V(k)-\underline{V}(k)\geq  i h(N)/2\},\\
    m_i
& := &
    \inf\Big\{k\in\N,\quad  V(k)=\inf_{0\leq \ell\leq \t_{i}}V(\ell)\Big\},
\label{eqDefmi}
        %t_i^- & := & 0\wedge (\sup\{t<x_1(\sigma W, i h(N)),\ \sigma W(t)-\sigma W(x_1(\sigma W, i h(N)))\geq i h(N)\}),\\
        %t_i^- & := & x_0(\sigma W, i h(N)),\\
        %\tt_i & := & \sup\{k\in\Z,\ k\leq m_i,\quad V(k)-\inf_{k\leq \ell\leq m_i}V(\ell)\geq i h(N)/2+2h(N)\},
        %\tt_i & := & £\lfloor t_i^- \rfloor+\tt_0,
\end{eqnarray}
where $\t_0$ is defined in \eqref{eqdeft0}.
%where $\t_0$ and $\tt_0$ are defined in \eqref{eqdeft0} and \eqref{eqdeftt0}.
In particular, $\sigma W(t_i)=\sigma \underline{W}(t_i)+i h(N)$ by continuity of $W$.
Moreover,  $\omega\in \mathcal{G}_7$, so  $x_0(\sigma W, i h(N))=x_0(W, i h(N)/\sigma)$ is an $i h(N)$-maximum
for $\sigma W$ and $x_1(\sigma W, i h(N))$ an $i h(N)$-minimum
for $\sigma W$ for every integer $i\geq 1$ such that $1\leq i h(N)\leq 5\log N$.
Consequently for such $i$,
$t_i\geq x_1(\sigma W, i h(N))$, otherwise there would be an $i h(N)$-maximum for $\sigma W$
in $(0, x_1(\sigma W, i h(N)))$, which is not possible.
Moreover, $\sigma W[x_2(\sigma W,  i h(N))]-\sigma W[x_1(\sigma W,  i h(N))]\geq i h(N)$, which gives
$t_i\leq x_2(\sigma W,  i h(N))$.
Hence,
\begin{equation}\label{inegalitesTi}
x_0(\sigma W,  i h(N))
\leq 0< x_1(\sigma W, i h(N))< t_i\leq x_2(\sigma W, i h(N)),
\end{equation}
%then $\underline{W}(t_i)=W[x_1(\sigma W, i h(N))]$, and
then
\begin{eqnarray}
    \inf\{W(t),\ x_0(\sigma W,  i h(N))\leq t \leq t_i\} & = & W[x_1(\sigma W, i h(N))]
    \label{eq41},
\\
    \sup\{W(t),\ x_0(\sigma W,  i h(N))\leq t \leq t_i\} & = & W[x_0(\sigma W,  i h(N))],
    \label{eq42}
\end{eqnarray}
since
%({\bf devenu inutile ? Supprimer ?})
$\sigma W[x_0(\sigma W,  i h(N))]\geq \sigma W[x_1(\sigma W,  i h(N))]+ih(N)=\sigma W(t_i).$
%and similarly $\sigma W[x_2(\sigma W,  i h(N))]\geq \sigma W(t_i)$.
We set similarly as in \eqref{eqdefXbarre},
\begin{equation*}
\underline{V}(n)    :=  \inf\{V(k),\ 0\leq k\leq n\},\quad n\in \N.\\
\end{equation*}
We recall that $C_4=\sigma+\frac{3-\sqrt{5}}{2C_3}C_1+4C_1$ and notice that
similarly as in \eqref{eqC4},
%\eqref{eqC4} still holds with $\VV$ instead of $V$.
\begin{equation}\label{eqC4Hat}
    \forall u\in\big[-(\log N)^4 , (\log N)^4\big],
\qquad
    \big|\VV(\lfloor u\rfloor)-\sigma W(u)\big|
\leq
    C_4\log_2 N.
\end{equation}
We also introduce  $i_{max}(N):=\max\{i\in\N,\  i h(N)\leq 4\log N\}$.
Since $\e\d/4< 4$ and $\mathcal{G}(N)\subset  \mathcal{G}_3(N)\cap \mathcal{G}_5(N)$, we get uniformly on $\mathcal{G}(N)$ for large $N$,
\begin{equation}\label{inegthetaii}
    \forall 1\leq i\leq i_{max}(N),
\qquad
    0\leq m_i\leq \t_i\leq \lfloor d_{\sigma W^+}(5\log N)\rfloor+\t_0\leq 2(\log N)^4;
\qquad
    0\leq t_i\leq (\log N)^4.
\end{equation}
%for $N$ large enough (say $N\geq N_1$). (faire de meme pour les i n\'egatifs ?)
% voir 08/11/08.
We now define for $1\leq i\leq i_{max}(N)$, with $\nu:=8+2C_4$,
\begin{equation}\label{eqDefFi}
    F_i(N)
:=
    \left\{\sum_{k=0}^n f(S_k)>0\quad  \forall 1\leq n\leq \tau(\t_i)\right\}
    \cap
    \Bigg\{ \sum_{k=0}^{\tau(\t_i)} f(S_k)\geq\frac{\exp(i h(N))}{2(\log N)^{\nu}}\Bigg\}.
\end{equation}
Our aim in the following is to prove, by induction on $i$, a lower bound for $\po(F_i(N))$ for $1\leq i\leq i_{max}(N)$.
We also prove that $\tau(\theta_i)\geq N$ for $i=i_{max}(N)$ with high probability.
We start with $i=1$.

\begin{lemma}\label{LemmaLowerBoundProbaFpouri1}
There exists a constant $c_7>0$ such that for $N$ large enough,
\begin{equation}\label{eqrec1}
\forall\o\in \mathcal{G}(N),\qquad \po(F_1(N))\geq c_7-4(\log N)^{-6}.
%\po(F_i(N))\geq c_7-3i(\log N)^{-\a''-2},
\end{equation}
\end{lemma}

\noindent{\bf Proof:}
Recall that $\e_0\leq e^{V(-1)}\leq \e_0^{-1}$. Moreover,  we have
for $\o\in\mathcal{G}(N)$,
$V(k)\leq -\delta k$ for $0\leq k \leq \lfloor \e\log_2 N\rfloor$,
whereas
$
    V(k)
\leq
    -(\delta  \e/2) \log_2 N
$
for
$\lfloor \e\log_2 N\rfloor< k \leq \theta_0$,
and  for
$\theta_0< k \leq \theta_1$,
$$
V(k)= V(\theta_0)+\VV(k-\theta_0)\leq -5h(N)+\sigma W(k-\theta_0)+C_4\log_2 N\leq -4h(N)+C_4\log_2 N,
$$
thanks to \eqref{eqC4Hat} since $t_1\leq (\log N)^4$ by \eqref{inegthetaii}
for $N$ large enough so that $i_{max}(N)\geq 1$.
Let
$c_7:=\e_0\big(\e_0^{-1}+2\big(1-e^{-\delta}\big)^{-1}\big)^{-1}$.
We have
$
    \po[\tau(\t_1)<\tau(-1)]
=
    e^{V(-1)}\big(\sum_{k=-1}^{\t_1-1}e^{V(k)}\big)^{-1},
$
which is, due to the previous remarks, greater than or equal to
\begin{equation}\label{eqDefN0a}
\e_0\left[\e_0^{-1}+\sum_{k=0}^{\lfloor \e \log_2 N\rfloor}e^{-\d k}
+(\t_0-\lfloor \e\log_2 N\rfloor)(\log N)^{-(\d\e/2)}
+(\t_1-\t_0)e^{-4h(N)}(\log N)^{C_4}
\right]^{-1}\geq c_7,
\end{equation}
for every $\o\in\mathcal{G}(N)$ for large $N$ since $\theta_0-\lfloor \e\log_2 N\rfloor\leq (\log N)^{\e\d/4}$ on $\mathcal{G}_3(N)$
and due to \eqref{inegthetaii}.

Moreover on $\mathcal{G}_1(N)$,
$\theta_1\geq m_1\geq \theta_0 \geq \lfloor \e \log_2 N\rfloor$, which is greater than $1$ for large $N$, so $f(m_1)\geq 1$.
Observe that  on $\{\tau(\t_1)<\tau(-1)\}$, due to \eqref{eqFormuleOccupation} and since $f(m_1)\geq 1$ and $f\geq 0$ on $\N$,
\begin{equation}\label{eqInegSommeetTempsLocalenm1}
\sum_{k=0}^{\tau(\t_1)} f(S_k)\geq L(m_1,\tau(\t_1)\wedge \tau(-1)),
\qquad
\sum_{k=0}^{n} f(S_k)\geq f(1)>0, \quad 1\leq n\leq \tau(\t_1).
\end{equation}
%This and \eqref{eqDefN0a} lead to
%\begin{equation*}
%\po\left(\sum_{k=0}^n f(S_k)>0\quad  \forall 1\leq n\leq \tau(\t_1)\right)\geq
%\po[\tau(\t_1)<\tau(-1)]\geq c_7>0.
%\end{equation*}
In order to give a lower bound of $L(m_1,\tau(\t_1)\wedge \tau(-1))$, notice that thanks to \eqref{eqC4Hat}
and since $t_1\leq (\log N)^4$
and $\sigma W(t_1)=\sigma \underline{W}(t_1)+h(N)$, we have for $\o\in\mathcal{G}(N)$,
$$
    \VV(m_1-\theta_0)
%\leq
%    \sigma W(m_1-\theta_0)+C_4\log_2 N
\leq
    \sigma\underline{W}(t_1)+C_4\log_2 N
\leq
    \VV(\lfloor t_1\rfloor)-h(N)+2C_4\log_2 N.
$$
Consequently,
uniformly on $\mathcal{G}(N)$  for large $N$,
we have $m_1+1<\theta_1$ and
\begin{eqnarray*}
    \po^{m_1+1}[\tau(\t_1)<\tau(m_1)]
& = &
    e^{V(m_1)}\bigg(\sum_{k=m_1}^{\t_1-1}e^{V(k)}\bigg)^{-1}\leq
    e^{V(m_1)-V(\t_1-1)}\leq \e_0^{-1}e^{-h(N)}(\log N)^{2C_4},
\\
    \po^{m_1-1}[\tau(-1)<\tau(m_1)]
& = &
    e^{V(m_1-1)}\bigg(\sum_{k=-1}^{m_1-1}e^{V(k)}\bigg)^{-1}
    \leq e^{V(m_1-1)-V(0)}\leq \e_0^{-1}e^{-h(N)}(\log N)^{2C_4}
    \end{eqnarray*}
since $V(m_1)\leq V(\t_0)\leq -5h(N)\leq -h(N)+2C_4\log_2 N$.
We know that $L(m_1, \tau(-1)\wedge\tau(\theta_1))$ is under $\po^{m_1}$ a geometric random variable of parameter
$\po^{m_1}[\tau(-1)\wedge\tau(\theta_1)<\tau^*(m_1)]$,
where $\tau^*(m_1):=\inf\{k\in\N^*,\ S_k=m_1\}$ is the first return time to $m_1$.
Hence,
$$
%\po^{m_1}(\tau(0) \wedge\tau(\t_1)\geq k)
%\geq
    \po^{m_1}\big[L(m_1,\tau(-1)\wedge\tau(\theta_1))>k\big]
\geq
    \big(\po^{m_1}[\tau(-1) \wedge\tau(\t_1)> \tau^*(m_1)]\big)^k
\geq
    \left(1-\frac{\e_0^{-1}(\log N)^{2C_4}}{ e^{h(N)}}\right)^k.
$$
Taking $k=k_1:=\big\lfloor \frac{\exp(h(N))}{2(\log N)^{\nu}}\big\rfloor $, we obtain uniformly on $\mathcal{G}(N)$ for large $N$,
\begin{equation*}
%    \log \po^{m_1}(\tau(0)\wedge \tau(\theta_1)\geq k_1)
    \log \po^{m_1}\big[L(m_1,\tau(-1)\wedge\tau(\theta_1))> k_1\big]
% \geq
%    k_1 \log(1-\e_0^{-1} e^{-h(N)}(\log N)^{2C_4})\\
 \geq
    -2 k_1 \e_0^{-1} e^{-h(N)}(\log N)^{2C_4}\\
 \geq
    -(\log N)^{-6}.
\end{equation*}
%$$\log \po^{m_1}(\tau(0) \wedge\tau(\t_1)\geq k_1)\geq-(\log N)^{-6}.$$
Hence,
\begin{equation}\label{eqProbaLog-6}
\po^{m_1}\big[L(m_1,\tau(-1)\wedge\tau(\theta_1))\leq k_1\big]
%=
%    1-\po^{m_1}[\tau(0) \wedge\tau(\t_1)\geq k_1]
\leq
    1-\exp(-(\log N)^{-6})
\leq
    (\log N)^{-6}.
\end{equation}
Since $f(k)\geq 1$ for $k\geq 1$ and $f(0)=0$, we have, using twice \eqref{eqInegSommeetTempsLocalenm1},
\begin{eqnarray*}
    \po[\tau(\t_1)<\tau(-1)]
& = &
    \po\left[\sum_{k=0}^n f(S_k)>0 \ \forall 1\leq n\leq \tau(\theta_1), \ \tau(\t_1)<\tau(-1)\right]\\
& \leq &
    \po[F_1(N)]
    +\po[\tau(\t_1)<\tau(-1),L(m_1,\tau(\t_{1})\wedge\tau(-1))\leq k_1].
\end{eqnarray*}
We get in particular for large $N$ by the strong Markov property, \eqref{eqDefN0a} and \eqref{eqProbaLog-6},
\begin{eqnarray*}
\forall \omega\in\mathcal{G}(N), \qquad
    \po[F_1(N)]
& \geq &
    \po[\tau(\t_1)<\tau(-1)]-\po^{m_1}\big[L(m_1,\tau(\t_{1})\wedge\tau(-1))\leq  k_1\big]\\
& \geq &
    c_7-(\log N)^{-6}.
\end{eqnarray*}
This gives \eqref{eqrec1} for  $N$ large enough.
\hfill$\Box$
\medskip

We now set $C_5:=11+2C_4$.
By Lemma \ref{LemmaLowerBoundProbaFpouri1}, there exists $N_\e\in\N$ such that for every $N\geq N_\e$,
inequality \eqref{eqrec} holds for $i=1$, \eqref{inegthetaii} holds for every $\o\in\mathcal{G}(N)$,
$\lfloor \e \log_2 N\rfloor\geq 1$,
and the following conditions are satisfied:
\begin{eqnarray}\label{DefNepsilon1}
\forall N\geq N_\e, & & \quad \log N\geq h(N)\geq (C_5+17+8C_4)\log_2 N\geq
%\log\frac{1-\e_0}{\e_0}+
4+3\e_0^{-2},
\\
\label{DefNepsilon2}
    \forall N\geq N_\e,
& & \quad
    \min_{[-(\log N)^4,0]}f\geq -\exp\big((\log^4 N)^{\e\delta/2^7}\big)
=
    -e^{h(N)}.
\end{eqnarray}
We prove by induction on $i$ the following lemma:

\begin{lemma}\label{LemmaLowerBoundProbaFpouriGeneral}
For all $N\geq N_\e$ and for every $1\leq i\leq i_{max}(N)$,
\begin{equation}\label{eqrec}
\forall \o\in\mathcal{G}(N),\qquad
\po[F_i(N)]\geq c_7-4i(\log N)^{-6}.
%\po(F_i(N))\geq c_7-3i(\log N)^{-\a''-2},
\end{equation}
% where ($\a''=4$),
Moreover for all $N\geq N_\e$,
\begin{equation}\label{eqTauThetaImax}
\forall \o\in\mathcal{G}(N),\qquad
    \po\big[\tau\big(\theta_{i_{max}(N)}\big)\geq N\big]
\geq
    1-2(\log N)^{-6}.
\end{equation}
\end{lemma}

\bigskip
\noindent{\bf Proof:}
We fix $N\geq N_\e$. We already know that \eqref{eqrec} is true for $i=1$.
Now, assume \eqref{eqrec} is true for an integer $i$ such that $1\leq i\leq i_{max}(N)-1$, and let us prove it is true for $i+1$. We fix $\o\in\mathcal{G}(N)$.

We notice that $\theta_i<\theta_{i+1}$. Indeed, if $\underline{W}(t_i)=\underline{W}(t_{i+1})$, we have
$\sigma W(t_{i+1})=\sigma W(t_i)+h(N)$, which gives, since $N\geq N_\e$,
$\VV(\lfloor t_{i+1}\rfloor)
%\geq \sigma W(t_{i+1})-C_4\log_2 N
%= \sigma W(t_{i})+h(N)-C_4\log_2 N
\geq \VV(\lfloor t_{i}\rfloor)+h(N)-2C_4\log_2 N
>\VV(\lfloor t_{i}\rfloor)$
by \eqref{eqC4Hat} and \eqref{inegthetaii},
so $\theta_{i+1}\neq\theta_i$.
If $\underline{W}(t_i)\neq\underline{W}(t_{i+1})$, there exists $u\in[t_i, t_{i+1}]$ such that
$|\sigma W(u)-\sigma W(t_i)|>i h(N)$, and $\theta_i=\theta_{i+1}$ would imply
$|u-t_i|\leq |t_{i+1}-t_i|\leq 1$
and then contradict $\omega\in\mathcal{G}_6(N)$ for $N\geq N_\e$. So, $\theta_i<\theta_{i+1}$.

\noindent
{\bf First step :}
Define (see Figure \ref{figure1})
\begin{eqnarray}\label{eqDefThetaiBarre}
    \tn_i
& := &
    \max\{k\in\Z,\ k<\t_i,\ V(k)\geq V(\t_i)+h(N)+C_5\log_2 N\},
\\
    E_{6,i}
& := &
    \left\{
    \inf\{k\geq \tau(\theta_i),\ S_k=\theta_{i+1}\}
    <
    \inf\{k\geq \tau(\theta_i),\ S_k=\tt_i\}
    \right\}
    =\{\tau(\theta_{i+1})<\tau(\theta_i,\tt_i)\},
\nonumber
\end{eqnarray}
where
$$
\forall (a,b)\in\Z^2,\qquad \tau(a,b):=\inf\{k\geq \tau(a),\ S_k=b\}.
$$
%+1
%and $E_{6,i}:=\{\tau(\tn_{i})\circ\tau(\theta_i)>\tau(\t_{i+1})\circ\tau(\theta_i)\}$.
We prove that
%\begin{equation}\label{eqStep1}
$
\po(E_{6,i}^c)\leq (\log N)^{-6}.
$
%\end{equation}
First, notice that since
%$\inf_{[0,t]}W \leq \inf_{[0,t_{i}]}W$
$\underline{W}(t)\leq \underline{W}(t_i)$
for $t_i\leq t\leq t_{i+1}$,
applying twice \eqref{eqC4Hat} gives
\begin{equation}
\max_{[\t_i, \t_{i+1}]}V  \leq  V(\t_i)+h(N)+2C_4\log_2 N.\label{eqMaxVEntreLesThetas}
\end{equation}
Hence, applying the Markov property at time $\tau(\theta_i)$, we get since $\theta_{i+1}\leq 2(\log N)^4$ by \eqref{inegthetaii},
\begin{equation}\label{eqStep1}
\po(E_{6,i}^c)
%=
%\po^{\theta_i}[\tau(\tn_{i})<\tau(\t_{i+1})]
=
\frac{\sum_{k=\t_i}^{\t_{i+1}-1}e^{V(k)}}{\sum_{k=\tn_{i}}^{\t_{i+1}-1}e^{V(k)}}
\leq
\frac{2(\log N)^{4+2C_4} e^{V(\t_i)+h(N)}}{e^{V(\t_{i})+h(N)}(\log N)^{C_5}}
%\leq
%\frac{2(\log N)^4 \max_{k\in[\t_i,\t_{i+1}] }e^{V(k)}}{e^{V(\tn_{i})}}
\leq
(\log N)^{-6}.
%\leq (\log N)^{4+2C_4-{C_5}}.
\end{equation}
%which gives \eqref{eqStep1}.

%\begin{equation*}
%\po^{\t_i}(\tau(\tt_{i})<\tau(\t_{i+1}))
%=\frac{\sum_{k=\t_i}^{\t_{i+1}-1}e^{V(k)}}{\sum_{k=\tt_{i}}^{\t_{i+1}-1}e^{V(k)}}
%\leq  \frac{(\log N)^4 e^{V(\t_0)+\sigma W(t_{i})+h(N)+C_4\log_2 N}}{e^{V(\tt_{i})}}
%\leq \frac{(\log N)^{4+2C_4}}{h(N)^9}.
%\end{equation*}

\bigskip
\noindent
{\bf Second step :}
We recall that for every $\Delta\subset \Z$ and $0\leq s\leq t$,
$L(\Delta, s\rightsquigarrow t)=\sum_{k=s}^t \mathds{1}_{\{S_k\in \Delta\}}$
is the number of visits of
$(S_n)_{n\in\N}$ to the set $\Delta$ between times $s$ and $t$,
as defined after \eqref{e1p20}.
In this step,
we consider
\begin{equation*}\label{eqStep2}
    E_{7,i}
:=
    \left\{L\big((\tn_i,0), \tau(\theta_i)\rightsquigarrow\tau(\theta_i,\tn_{i})\wedge\tau(\t_{i+1})\big)< \exp[(i-3)h(N)]\right\},
\end{equation*}
and we show that
\begin{equation}\label{eqSecondStep}
    \po(E_{7,i}^c)
\leq
    (\log N)^{-6}.
\end{equation}

We consider separately two cases.
%, depending on our (fixed) $N$, $i$ and $\o$.

\noindent \underline{First case:} Assume that $\tn_i\geq -1$. Then, $(\tn_i,0)\cap \Z=\emptyset$,
hence
$$
    L\big((\tn_i,0), \tau(\theta_i)\rightsquigarrow\tau(\theta_i,\tn_{i})\wedge\tau(\t_{i+1})\big)
=
    0
<
    \exp[(i-3)h(N)].
$$
Consequently in this case,
$
    \po(E_{7,i}^c)
=
    \po(\emptyset)
=
0
\leq
    (\log N)^{-6}
$,
which proves \eqref{eqSecondStep} and then the second step in this first case.
We notice in particular that for $i\in\{1,2,3\}$,
since $V(\theta_0)\leq -5h(N)$ by \eqref{eqdeft0} and $\sigma W(t_i)=\sigma \underline{W}(t_i)+i h(N)$,
using \eqref{eqC4Hat} applied to $t_i$ (because $t_i\leq (\log N)^4$ by  \eqref{inegthetaii}),
$$
    V(\theta_i)
=
    V(\theta_0+\lfloor t_i\rfloor)
=
    V(\theta_0)+\VV(\lfloor t_i\rfloor)
\leq
    -5h(N)+\sigma \underline{W}(t_i)+i h(N)+C_4\log_2 N.
$$
Since $\sigma \underline{W}(t_i)\leq 0$ and $(C_5+C_4)\log_2 N\leq h(N)$ by \eqref{DefNepsilon1}, this gives for $i\in\{1,2,3\}$,
$$
    V(\t_i)+h(N)+C_5\log_2 N
\leq
    (i-4)h(N)+(C_5+C_4)\log_2 N
\leq
    -h(N)+(C_5+C_4)\log_2 N
\leq
    0,
$$
and so $\tn_i\geq 0$.
So if $i\leq 3$, we are automatically in the first case. Heuristically, this is due to the fact that we
have lowered the potential in $[\theta_0,+\infty)$
%forced the potential of the good environments to decrease to
by the quantity $|V(\theta_0)|\geq 5h(N)$, which is quite large, in our definitions \eqref{eqdeft0} and \eqref{eqDefVChapeau}
of $\theta_0$ and $\VV$.

%Since this is obvious when $\tn_i\geq -1$, we just have to consider the case
%$\tn_i<-1$.

\noindent \underline{Second case:} Assume that $\tn_i<-1$, which implies that $i\geq 4$ due to the previous remark.
First, notice that since $x_0(\sigma W, i h(N))$ is a $i h(N)$-maximum for $\sigma W$, we have by \eqref{eqC4Hat}
since $x_0(\sigma W, i h(N))\geq -d_{-\sigma W^-}(5\log N)\geq -(\log N)^4$ (where we used $i\leq i_{max}(N)$),
%({\bf a justifier})
\begin{eqnarray}
% V(\tt_i) & \geq & 5 h(N)+\sigma W[x_0(\sigma W, i h(N))]-C_4\log_2 N\\
%         & \geq & 10h(N)+V(\t_0)+\sigma W[x_1(\sigma W, i h(N))]+i h(N)-C_4\log_2 N\\
%         & \geq & V(\t_i)+10h(N)-2C_4\log_2 N\\
%         & \geq & V(\tn_i),
V(\lfloor x_0(\sigma W, i h(N))\rfloor)
 & \geq & \sigma W[x_0(\sigma W, i h(N))]-C_4\log_2 N\nonumber\\
         & \geq & \sigma W[x_1(\sigma W, i h(N))]+i h(N)-C_4\log_2 N\nonumber\\
         & \geq & \sigma W[x_1(\sigma W, i h(N))]+i h(N)-C_4\log_2 N+5h(N)+V(\t_0).~~~~ \label{exProvisoireVx0}
\end{eqnarray}
Moreover, $\sigma W(t_i)=\sigma \underline{W}(t_i)+i h(N)$, and $\underline{W}(t_i)=W[x_1(\sigma W, i h(N))]$ due to \eqref{eq41}.
This together with \eqref{eqC4Hat} and $t_i\leq (\log N)^4$ (see \eqref{inegthetaii}) gives
\begin{equation}\label{InegTthetaI}
V(\theta_i)-V(\theta_0)=\VV(\lfloor t_i\rfloor)
\leq \sigma W(t_i)+C_4\log_2 N
=\sigma W[x_1(\sigma W, i h(N))]+i h(N)+C_4\log_2 N.
\end{equation}
Hence, \eqref{exProvisoireVx0} and then $N\geq N_\e$ and \eqref{DefNepsilon1} lead to
\begin{equation*}
V(\lfloor x_0(\sigma W, i h(N))\rfloor)
          \geq  V(\t_i)+5h(N)-2C_4\log_2 N
          >  V(\tn_i).
\end{equation*}
Consequently,    $\lfloor x_0(\sigma W, i h(N))\rfloor < \tn_i<\t_i<\t_{i+1}$ by definition of $\tn_i$.
%Moreover,
%$V(\tt_{i})
%%\geq \sigma W(t_{i}^-) -C_4\log_2 N+V(\tt_0)
%\geq \sigma W(t_i) -C_4\log_2 N-5h(N)
%%\geq \sigma W(t_{i+1}) -C_4\log_2 N+4h(N)
%$.
Recalling that $\tn_i<-1$ in this second case, we can consider $z\in(\tn_i,0)\cap \Z$. We get by Lemma \ref{LemmaRevesz1},
\begin{eqnarray}
    \eo^{\t_i}\big[L\big(z, \tau(\tn_{i})\wedge\tau(\t_{i+1})\big)\big]
& = &
    \frac{\po^{\t_i}[\tau(z)<\tau(\t_{i+1})]}{\o_z e^{V(z)}\big(\sum_{k=z}^{\t_{i+1}-1}e^{V(k)}\big)^{-1}
    +(1-\o_z)\po^{z-1}[\tau(z)>\tau(\tn_i)]}\nonumber\\
& \leq &
    \e_0^{-1}e^{-V(z)}\sum_{k=z}^{\t_{i+1}-1}e^{V(k)}\nonumber\\
& \leq &
    3\e_0^{-1}(\log N)^4 \exp\Big(-V(z)+\max_{[z,\t_{i+1}]} V \Big),\label{eqInegEspL1}
\end{eqnarray}
since $\t_{i+1}\leq 2(\log N)^4$ by \eqref{inegthetaii} and $z>\tn_i\geq \lfloor x_0(\sigma W,i h(N))\rfloor\geq-(\log N)^4$.
We notice that by \eqref{InegTthetaI} and since $V(\theta_0)\leq -5h(N)$ by \eqref{eqdeft0},
\begin{eqnarray}\label{InegsVThetaIbis}
V(\t_i)=V(\theta_0)+\VV(\t_i-\theta_0)\leq -5 h(N)+\sigma W[x_1(\sigma W, i h (N))]+i h(N)+C_4\log_2 N.
%+\log(\e_0^{-1}).
\end{eqnarray}
%Notice that for every
%$x\in[x_0(\sigma W, i h (N)),x_1(\sigma W, i h (N))]$, we have $\sigma W(x_1(\sigma W, i h (N)))\leq \sigma W(x)\leq \sigma W(x_0(\sigma W, i h (N)))$, because $x_1(\sigma W, i h (N))$ is an $i h(N)$-minimum and $x_0(\sigma W, i h (N))$ and $x_1(\sigma W, i h (N))$ are two consecutive $i h(N)$-extrema.
%Hence,
Since $-\log^4 N\leq x_0(\sigma W, i h(N))\leq \tt_i< z<0< x_1(\sigma W, i h(N))\leq t_i\leq \log^4 N$
by \eqref{inegalitesTi} and \eqref{inegthetaii},
equations \eqref{eq41},
\eqref{eqC4Hat} and \eqref{InegsVThetaIbis} give
%for $\tn_i< z\leq 0$,
\begin{eqnarray}
%V(z)& =& V(\tt_0)+\VV(z-\tt_0)\geq 5h(N)+\sigma W(x_1(\sigma W, i h (N)))-C_4\log_2 N\nonumber\\
%& \geq & (10-i) h(N)+V(\t_i)-2C_4\log_2 N\label{eqInegEspL2}
    V(z)
=
    \VV(z)
\geq
%=\VV(z)
    \sigma W[x_1(\sigma W, i h (N))]-C_4\log_2 N
\geq
    (5-i) h(N)+V(\t_i)-2C_4\log_2 N.
\label{eqInegEspL2}
\end{eqnarray}
%If $\tt_0< z\leq 0$ then $V(z)\geq 0\geq V(\t_i)+(5-i)h(N)-C_4\log_2 N$.
Moreover by definition \eqref{eqDefThetaiBarre} of $\tn_i$ and \eqref{elliptic},
\begin{equation*}
\max_{[\tn_i,\t_i]}V \leq  V(\t_i)+h(N)+C_5\log_2 N-\log\e_0.
\end{equation*}
Combining this with \eqref{eqMaxVEntreLesThetas},
 \eqref{eqInegEspL1} and \eqref{eqInegEspL2} gives since $N\geq N_\e$,
\begin{eqnarray*}
    \eo^{\t_i}\big[L\big(z, \tau(\tn_{i})\wedge\tau(\t_{i+1})\big)\big]
& \leq &
    3\e_0^{-1}(\log N)^4 \exp\big(-V(z)+V(\t_i)+h(N)+C_5\log_2 N-\log \e_0\big)\\
& \leq &
    3\e_0^{-2}(\log N)^{C_5+2C_4+4}e^{(i-4)h(N)} \\
& \leq &
    (\log N)^{-10}e^{(i-3)h(N)}.
\end{eqnarray*}
Summing this over $z$ gives
$
    \eo^{\t_i}\big[L\big((\tn_i, 0), \tau(\tn_{i})\wedge\tau(\t_{i+1})\big)\big]
\leq
    (\log N)^{-6}e^{(i-3)h(N)}
$
since $\tn_i\geq -(\log N)^4$. We get $\po(E_{7,i}^c)=\po^{\t_i}(E_{7,i}^c)\leq (\log N)^{-6}$ by Markov's inequality and property.
This proves \eqref{eqSecondStep} in this second case, which ends the second step.

\bigskip
\noindent
{\bf Third step :} we define
\begin{equation*}
    E_{8,i}
:=
    \left\{L\big(m_{i+1}, \tau(\theta_i)\rightsquigarrow \tau(\theta_i,\tt_{i})\wedge\tau(\t_{i+1})\big)
                    > \frac{\exp[(i+1)h(N)]}{(\log N)^\nu}
    \right\}.
\end{equation*}
We prove that %(with $\nu=8+2C_4$) %(remplacer $\t_i^-$ par $\tn_i$ ?)
\begin{equation}\label{eqStep3}
\po\left(E_{8,i}^c\right)\leq 2(\log N)^{-6}.
\end{equation}
To this aim, we first show that
\begin{equation}\label{ep16}
    \po^{\t_i}\big[\tau(m_{i+1})>\tau\big(\tt_{i}\big)\wedge\tau(\t_{i+1})\big]
\leq
    (\log N)^{-6}.
\end{equation}
 This is true if $\t_i\leq m_{i+1}\leq \t_{i+1}$ by \eqref{eqStep1}.
%(attention a $C$...).
Else, $m_i=m_{i+1}< \theta_i$ and then
$
    \sigma \underline{W}(t_{i+1})
\geq
    \sigma \underline{W}(t_i)-2C_4\log_2 N
$
by \eqref{eqC4Hat} and \eqref{inegthetaii}, which leads to
\begin{equation}\label{InegVThetaiThetaSuivant}
    V(\t_{i+1})
\geq
    V(\t_i)+h(N)-4C_4\log_2 N.
\end{equation}
We get successively, again
by \eqref{eqC4Hat} and \eqref{inegthetaii}, for every $m_i\leq k\leq \theta_i$,
\begin{align}
\nonumber
    \sigma W(m_i-\theta_0)
& \leq
    \VV(m_i-\theta_0)+C_4\log_2 N
=
    \inf\nolimits_{[0, \lfloor t_i]\rfloor}\VV+C_4\log_2 N
\leq
    \sigma\underline{W}(t_i)+2C_4\log_2 N,
\\
    \VV(k-\theta_0)
& \leq
    \sigma[W(k-\theta_0)-W(m_i-\theta_0)]+\sigma W(m_i-\theta_0)+C_4\log_2 N
\nonumber\\
& \leq
    i h(N) +\sigma\underline{W}(t_i)+3C_4\log_2 N
=
    \sigma W(t_i)+3C_4\log_2 N,
\label{InegVTheta0}
\end{align}
where we used the definition of $t_i$ in the last inequality.
Using \eqref{InegVTheta0},
then \eqref{eqC4Hat} and \eqref{inegthetaii},
then the definitions \eqref{eqDefVChapeau} and  \eqref{eqDefthetai}  of $\VV$ and $\theta_i$,
and finally \eqref{InegVThetaiThetaSuivant}, we get
\begin{align}
    \max_{[m_i,\theta_i]}V
& \leq
    V(\theta_0)+\sigma W(t_i)+3C_4\log_2 N
\leq
    V(\theta_0)+\VV(\lfloor t_i\rfloor )+4C_4\log_2 N
\nonumber\\
&  =
    V(\theta_i)+4C_4\log_2 N
\label{InegPourmietThetaBarre}\\
& \leq
    V(\theta_{i+1})-h(N)+8C_4\log_2 N.
\label{eqMaxV}
\end{align}
In particular, \eqref{InegPourmietThetaBarre} combined with \eqref{DefNepsilon1}
and the definition \eqref{eqDefThetaiBarre} of  $\tt_{i}$ leads to  $\tt_{i}<m_i$.
Now, in this case $m_i=m_{i+1}<\theta_i<\theta_{i+1}$, we have
since $N\geq N_\e$,
$\theta_i\leq 2(\log N)^4$ by \eqref{inegthetaii},
and by \eqref{eqProbaAtteinte},
\begin{eqnarray*}
    \po^{\t_i}\big[\tau(m_{i+1})>\tau\big(\tt_{i}\big)\wedge\tau(\t_{i+1})\big]
& = &
    \po^{\t_i}[\tau(m_{i+1})>\tau(\t_{i+1})]\nonumber\\
& \leq &
    2(\log N)^4\e_0^{-1}\exp\Big[\max_{[m_i,\theta_i]}V-V(\theta_{i+1})\Big],
%& \leq &
%    2\e_0^{-1}  (\log N)^{8C_4+4} e^{-h(N)},
%& \leq &
%    (\log N)^{-6}.
\end{eqnarray*}
which together with \eqref{eqMaxV} gives \eqref{ep16} since $N\geq N_\e$.

Moreover, we prove that $\po^{m_{i+1}}(E_{9,i}^c)\leq (\log N)^{-6}$, where
$$
    E_{9,i}
:=
    \left\{
        L\big(m_{i+1}, \tau\big(\tt_{i}\big)\wedge\tau(\t_{i+1})\big)
        >
        \frac{\exp[(i+1) h(N)]}{(\log N)^{\nu}}
    \right\}.
$$
We know that $\tt_{i}<m_i\leq m_{i+1}<\theta_{i+1}$ thanks to \eqref{InegPourmietThetaBarre},
which is true in every case, as is \eqref{InegVTheta0}.
%We know that $\tt_{i}<m_{i+1}<\theta_{i+1}$ because $V(m_{i+1})\leq V(\theta_i)<\inf\{V(\tt_i), V(\theta_i)\}$
%and $\tt_i\leq m_i\leq m_{i+1}$ by definition of $t_i$, $\tt_i$ and \eqref{eqC4Hat}.
So,
$L\big(m_{i+1}, \tau\big(\tt_{i}\big)\wedge\tau(\t_{i+1})\big)$ is under $\po^{m_{i+1}}$ a geometric r.v. with parameter
\begin{eqnarray}
q_1 & := &
    \o_{m_{i+1}}\po^{m_{i+1}+1}[\tau(m_{i+1})>\tau(\t_{i+1})]
    +(1-\o_{m_{i+1}})\po^{m_{i+1}-1}\big[\tau(m_{i+1})>\tau\big(\tt_{i}\big)\big]
    \nonumber\\
& \leq &
    \o_{m_{i+1}}\e_0^{-1}e^{V(m_{i+1})-V(\theta_{i+1})}+(1-\o_{m_{i+1}})\e_0^{-1}e^{V(m_{i+1})-V(\tt_i)}
    \label{Inegq1},
\end{eqnarray}
by \eqref{eqProbaAtteinte}.
%We notice that
%$V(m_{i+1})\leq V(\theta_0)+\sigma \underline{W}(t_{i+1})+C_4\log_2 N
%\leq V(\theta_{i+1})+2C_4\log_2 N -(i+1)h(N)$. Moreover, since $\underline{W}(t_i)\geq \underline{W}(t_{i+1})$,
%we have $V(\theta_i)\geq V(\theta_0)+\sigma \underline{W}(t_{i+1})+i h(N)-C_4\log_2 N
%\geq V(m_{i+1})-2C_4\log_2 N+i h(N)$. This gives
%$V(\tt_i)\geq V(m_{i+1})+(i+1)h(N)+11\log_2 N$.
Moreover, we obtain successively the following inequalities:
\begin{align}
    V(m_{i+1})
&
    \leq V(\theta_0)+\sigma \underline{W}(t_{i+1})+C_4\log_2 N
=
    V(\theta_0)+\sigma W(t_{i+1})-(i+1)h(N)+C_4\log_2 N
\nonumber\\
&
    \leq V(\theta_{i+1})-(i+1)h(N)+2C_4\log_2 N,
\label{eqAlign1}\\
    V(\theta_i)
&
\geq
    V(\theta_0)+\sigma W(t_i)-C_4\log_2 N
=
    V(\theta_0)+\sigma \underline{W}(t_i)+i h(N)-C_4\log_2 N
\nonumber\\
&
\geq
    V(\theta_0)+\sigma \underline{W}(t_{i+1})+i h(N)-C_4\log_2 N
\geq
    V(m_{i+1})+i h(N)-2C_4\log_2 N,
\label{eqAlign2}\\
    V\big(\tt_i\big)
&
    \geq V(\theta_i)+h(N)+C_5\log_2 N
    \geq V(m_{i+1})+(i+1)h(N)+11\log_2 N,
\label{eqAlign3}
\end{align}
where we used
$
    V(m_{i+1})
\leq
    V(\theta_0+\lfloor z_{i+1}\rfloor)
=
    V(\theta_0)+\VV(\lfloor z_{i+1}\rfloor)
$ with $z_{i+1}\in[0,t_{i+1}]$ such that $W(z_{i+1})=\underline{W}(t_{i+1})$
and \eqref{eqC4Hat}
in the first inequality of \eqref{eqAlign1},
$\underline{W}(t_i)\geq \underline{W}(t_{i+1})$ in \eqref{eqAlign2},
and the definition \eqref{eqDefThetaiBarre} of $\tt_{i}$ in \eqref{eqAlign3}. It follows from
\eqref{Inegq1}, \eqref{eqAlign1} and \eqref{eqAlign3} that
\begin{eqnarray*}
q_1
& \leq & {\e_0}^{-1} \exp(-(i+1)h(N)+2C_4\log_2 N)=:q_2.
\end{eqnarray*}
%Now, define $n_2:=\frac{1}{2(\log N)^6 q_2}=\frac{\e_0^2\exp((i+1)h(N))}{2(\log N)^{6+2C_4}}$
%and $E_{21}:=\{L(m_{i+1},\tau(\tt_{i})\wedge\tau(\t_{i+1})) > n_2\}$. We have for $N\geq N_\e$,
%\begin{equation*}
%\log\po^{m_{i+1}}(E_{21})= n_2\log(1-q_1)\geq n_2\log(1-q_2)\geq -2 n_2 q_2=-(\log N)^{-6}.
%\end{equation*}
\noindent
Now, define $n_2:=\big\lfloor \frac{\exp[(i+1) h(N)]}{(\log N)^{\nu}}\big\rfloor $. We have for $N\geq N_\e$,
\begin{equation*}
    \log\po^{m_{i+1}}(E_{9,i})
=
    n_2\log(1-q_1)
\geq
    n_2\log(1-q_2)
\geq
    -2 n_2 q_2
\geq
    -(\log N)^{-6}.
\end{equation*}
Indeed, $q_2\in(0,1/2)$ hence $\log(1-q_2)\geq -2q_2$.
Since $1-e^{-t}\leq t$ for $t\in\R$, this yields
$\po^{m_{i+1}}(E_{9,i}^c)
%\leq\po^{m_{i+1}}(E_{21}^c)
\leq (\log N)^{-6}$.
Hence by the strong Markov property,
\begin{eqnarray*}
    \po(E_{8,i}^c)
& = &
    \po^{\theta_i}\big[E_{8,i}^c, \tau(m_{i+1})>\tau(\tt_{i})\wedge\tau(\t_{i+1})\big]
    +\po^{\theta_i}\big[E_{8,i}^c,\tau(m_{i+1})\leq\tau(\tt_{i})\wedge\tau(\t_{i+1})\big]
\\
& \leq &
    \po^{\t_i}\big[\tau(m_{i+1})>\tau(\tt_{i})\wedge\tau(\t_{i+1})\big]+\po^{m_{i+1}}(E_{9,i}^c)
\\
& \leq &
    2(\log N)^{-6},
\end{eqnarray*}
where we used \eqref{ep16} in the last inequality.
This gives \eqref{eqStep3}.
Moreover, notice that in the particular case $i=i_{max}(N)-1$, we get on $E_{8,i}$ since $N\geq N_\e$,
\begin{equation}\label{eqTempsFinal}
    \tau\big(\theta_{i_{max}(N)}\big)
\geq
    L\big(m_{i_{max}(N)}, \tau\big(\theta_{i_{max}(N)-1}\big)
            \rightsquigarrow
            \tau\big(\theta_{i_{max}(N)-1},\tt_{i_{max}(N)-1}\big)\wedge\tau\big(\t_{i_{max}(N)}\big)
        \big)
\geq
    N.
\end{equation}
This and \eqref{eqStep3} already prove \eqref{eqTauThetaImax},
since we did not yet use our induction hypothesis.

\bigskip
\noindent{\bf Fourth step: conclusion.}
First, let $\tau(\t_i)<n\leq \tau(\t_{i+1})$. We have in the case $\tn_i<-1$,
\begin{equation}\label{EqDecomposition}
\sum_{k=0}^n f(S_k)
=
\sum_{k=0}^{\tau(\t_i)-1} f(S_k)
+
\left(\sum_{z\leq \tt_i}+\sum_{\tt_i<z<0}+\sum_{z\geq 0}
\right)
f(z)L(z, \tau(\theta_i)\rightsquigarrow n).
\end{equation}
The second sum of the right hand side is $0$ on $E_{6,i}$, and the last sum is at least $f(\theta_i)$ because $f\geq 0$ on $\N$.
Since $f<0$ on $\Z_-^*$ and $\tt_i\geq -(\log N)^4$, we get on $E_{6,i}$,
%$F_i(N)\cap E_{6,i}\cap E_{7,i}$,
\begin{eqnarray*}
    \sum_{k=0}^n f(S_k)
& \geq &
    \sum_{k=0}^{\tau(\t_i)} f(S_k)
    +\Big(\min_{[-(\log N)^4,0]}f\Big)L\big(\big(\tn_i,0\big), \tau(\theta_i)\rightsquigarrow\tau\big(\theta_i,\tn_{i}\big)\wedge\tau(\t_{i+1})\big).
%    \nonumber\\
\end{eqnarray*}
Since  for $N\geq N_\e$, $0>\min_{[-(\log N)^4,0]}f\geq -e^{h(N)}$ by \eqref{DefNepsilon2} and $e^{h(N)}\geq \min\{(\log N)^\nu,4\}$,
we get on $F_i(N)\cap E_{6,i}\cap E_{7,i}$,
\begin{eqnarray}
    \sum_{k=0}^n f(S_k)
& \geq &
    \frac{\exp(i h(N))}{2(\log N)^\nu}-\exp[h(N)]\exp[(i-3)h(N)]>0.
    \label{eqP16}
\end{eqnarray}
The proof is similar if $\tn_i\geq -1$, since in this case on $E_{6,i}$,
for all $\tau(\theta_i)\leq k\leq n\leq \tau(\theta_{i+1})$,
$S_k\geq \tt_i+1\geq 0$ and then $f(S_k)\geq 0$, which leads to
$\sum_{k=0}^n f(S_k)\geq \sum_{k=0}^{\tau(\theta_i)} f(S_k)\geq \frac{\exp(i h(N))}{2(\log N)^\nu} >0$
on $F_i(N)\cap E_{6,i}\cap E_{7,i}$, which gives \eqref{eqP16} also in this case.

We now consider $\sum_{k=0}^{\tau(\t_{i+1})} f(S_k)$, which is on $E_{6,i}$ equal to
(assuming first that $\tn_i<-1$)
\begin{equation}\label{EqDecomposition2}
\sum_{k=0}^{\tau(\t_{i})-1} f(S_k)
    +\left(\sum_{z\leq \tt_i}+\sum_{\tt_i<z<0}+\sum_{z\in\N-\{m_{i+1}\}}+\sum_{z\in\{m_{i+1}\}}
    \right)f(z)L\big(z,\tau(\theta_i) \rightsquigarrow\tau\big(\theta_i,\tt_{i}\big)\wedge\tau(\t_{i+1})\big).
\end{equation}
%, which is equal to \eqref{EqDecomposition} with $n=\tau(\t_{i+1})$.
The potential $V$ is decreasing on $[0,\lfloor \e \log_2 N\rfloor]$ since $\o\in\mathcal{G}_1(N)$,
hence  $m_{i+1}\geq \lfloor \e\log_2 N\rfloor\geq 1$ since $N\geq N_\e$, and then $f(m_{i+1})\geq 1$. Consequently,
the last sum in the right hand side of \eqref{EqDecomposition2} is at least
$L\big(m_{i+1},  \tau(\theta_i) \rightsquigarrow\tau\big(\theta_i,\tt_{i}\big)\wedge\tau(\t_{i+1})\big)$.
Moreover, the first term is positive on $F_i(N)$, the second one is $0$ on $E_{6,i}$, and the forth one is nonnegative
since $f\geq 0$ on $\N$. So,
we have on $F_i(N)\cap E_{6,i}\cap E_{7,i}\cap E_{8,i}$ for $N\geq N_\e$, since $\tt_i\geq -(\log N)^4$,
\begin{eqnarray*}
\sum_{k=0}^{\tau(\t_{i+1})} f(S_k)
%& = &
%    \sum_{k=0}^{\tau(\t_{i})-1} f(S_k)
%    +\sum_{k=\tau(\t_{i})}^{\tau(\tn_i)\wedge  \tau(\t_{i+1})} f(S_k)\\
%& = & \sum_{k=0}^{\tau(\t_{i})-1} f(S_k)
%    +\left(\sum_{z\leq \tt_i}+\sum_{\tt_i<z<0}+\sum_{z\in\N-\{m_{i+1}\}}+\sum_{z\in\{m_{i+1}\}}
%    \right)f(z)L(z,\tau(\theta_i) \rightsquigarrow\tau(\tt_{i})\wedge\tau(\t_{i+1}))
%    \\
& \geq &
    %f(m_{i+1})
    L\big(m_{i+1},  \tau(\theta_i) \rightsquigarrow\tau\big(\theta_i,\tt_{i}\big)\wedge\tau(\t_{i+1})\big)
\\[-2mm]
&&
%\qquad
    +\Big(\min_{[-(\log N)^4,0]}f\Big)
    L\big(\big(\tn_i,0\big), \tau(\theta_i) \rightsquigarrow\tau\big(\theta_i,\tn_{i}\big)\wedge\tau(\t_{i+1})\big)
\end{eqnarray*}
This gives on $F_i(N)\cap E_{6,i}\cap E_{7,i}\cap E_{8,i}$ for $N\geq N_\e$,
\begin{equation}
    \sum_{k=0}^{\tau(\t_{i+1})} f(S_k)
\geq
    \frac{\exp[(i+1) h(N)]}{(\log N)^{\nu}}
    -\exp[(i-2)h(N)]
\geq
    \frac{\exp[(i+1) h(N)]}{2(\log N)^{\nu}}.
\label{eqPart4Bis}
\end{equation}
We get \eqref{eqPart4Bis} similarly if $\tn_i\geq -1$, since in this case on $E_{6,i}$,
$f(S_k)\geq 0$ for all $\tau(\theta_i)\leq k\leq \tau(\theta_{i+1})$ as explained after \eqref{eqP16},
and so
$
    \sum_{k=0}^{\tau(\t_{i+1})} f(S_k)
\geq
    \sum_{k=0}^{\tau(\t_{i})-1} f(S_k)
    +
    L\big(m_{i+1},  \tau(\theta_i) \rightsquigarrow\tau\big(\theta_i,\tt_{i}\big)\wedge\tau(\t_{i+1})\big)
$,
which also leads as previously to \eqref{eqPart4Bis} in this case.

Now, \eqref{eqP16} and \eqref{eqPart4Bis} yield $F_i(N)\cap E_{6,i}\cap E_{7,i}\cap E_{8,i}\subset F_{i+1}(N)$. Consequently,
our induction hypothesis $\po[F_i(N)]\geq c_7-4i(\log N)^{-6}$
and inequalities \eqref{eqStep1}, \eqref{eqSecondStep} and \eqref{eqStep3}
% and $\po(E_{7,i}^c)\leq (\log N)^{-6}$
give for every $\o\in\mathcal{G}(N)$,
\begin{equation}\label{eqProbaInduction}
\po[F_{i+1}(N)]\geq
\po[F_{i}(N)]-\po(E_{6,i}^c)-\po(E_{7,i}^c)-\po(E_{8,i}^c)\geq c_7-4(i+1)(\log N)^{-6}.
\end{equation}
This ends the induction for all $N\geq N_\e$.
Hence \eqref{eqrec} is true for every $1\leq i\leq i_{max}(N)$ for each $N\geq N_\e$,
which ends the proof of Lemma \ref{LemmaLowerBoundProbaFpouriGeneral}.
\hfill$\Box$

\bigskip
\noindent{\bf Proof of Lemma \ref{LemmaRWinGoodEnv}:}
Notice that due to \eqref{eqrec} and \eqref{eqTauThetaImax} of Lemma \ref{LemmaLowerBoundProbaFpouriGeneral},
%Notice that thanks to \eqref{eqTempsFinal},  and the inclusion before \eqref{eqProbaInduction} for $i=i_{max}(N)-1$, we also show similarly as in %\eqref{eqProbaInduction} that
$
    \po\big[F_{i_{max}(N)}(N)\cap\big\{\tau\big(\theta_{i_{max}(N)}\big)\geq N\big\}\big]
\geq
    c_7-\frac{4i_{max}(N)}{(\log N)^6}-\frac{2}{(\log N)^6}
\geq
    c_7-\frac{18}{(\log N)^5}
$
for all $N\geq N_\e$ and $\o\in\mathcal{G}(N)$.
Consequently, for $N$ large enough,
\begin{eqnarray*}
    \forall\o\in\mathcal{G}(N),\quad
& &
    \po\left(\sum_{k=0}^n f(S_k)>0\ \forall 1\leq n\leq N\right)  \\
& \geq &
    \po\big[F_{i_{max}(N)}(N)\cap\big\{\tau\big(\theta_{i_{max}(N)}\big)\geq N\big\}\big]
% \geq  c_7-\frac{4i_{max}+1}{(\log N)^6}
  \geq
    \frac{c_7}{2}=:c_6,
\end{eqnarray*}
which proves Lemma \ref{LemmaRWinGoodEnv}.\hfill$\Box$

\medskip

%Consequently, for $N$ large enough, independent of $\o\in E$,
%\begin{eqnarray*}
%\po(\sum_{k=0}^n f(S_k)>0\ \forall 0<n\leq N) & \geq &
%\po(\sum_{k=0}^{n} f(S_k)>0\ \forall 0<n\leq \tau(\theta_{i_{max}}))-\po(\tau(\t_{i_{max})}<N)\\
%& \geq & c_7-3 i_{max}(\log N)^{-6}\geq c_7/2.
%\end{eqnarray*}
%This is true for every $\o\in \mathcal{G}(N)$.

Now, integrating \eqref{EqLemmaRWinGoodEnv} on $\mathcal{G}(N)$ and applying Lemma \ref{LemmaProbaGood} gives
%Consequently, by Lemma \ref{LemmaProbaGood}, we get
\begin{eqnarray*}
\P\left(\sum_{k=0}^n f(S_k)>0\ \forall 1\leq n\leq N\right) & \geq & c_6\eta(\mathcal{G}(N))\geq \frac{c_6c_1 \e\log_2 N}{(\log N)^{\frac{3-\sqrt{5}}{2}+\e(\gamma+\delta/32)}}
\end{eqnarray*}
for $N$ large enough. Now, let $\e\to 0$; this proves the lower bound of
Theorem \ref{TheoremePersistenceAdditive} for $u=0$ and then for every $u\leq 0$.

\mysection{Proofs of Fact \ref{FactExcessHeight0} and Lemma \ref{LemmaProbaStrongChangeOfSign}}\label{SectionPreuveduLemme21}
\subsection{Proof of Fact \ref{FactExcessHeight0} }

We first study the left continuity of some functions. The following lemma is more or less obvious, however we provide a proof for the sake of completeness.

\begin{lemma}\label{LemmaLeftcontinuity}
On $\{W\in\mathcal{W}\}$, for all $k\in\Z$, the functions $x_k(W,.)$,
$e(T_k(.))$ and $H(T_k(.))$ are left-continuous on $(0,+\infty)$.
More precisely, for all realization of $W$ in $\mathcal{W}$, for every $n\in\N^*$ and $x>0$, there exists $K_{x,n}\in(0,x)$ such that all
the functions $x_k(W,.)$, $k\in\{-n,\dots, n\}$, are constant on $[K_{x,n},x]$.
\end{lemma}

\noindent{\bf Proof:}
We assume throughout the proof that
%we are on $\{W\in\mathcal{W}\}$
$W\in \mathcal{W}$.
Let $x>0$. We first notice that $\lim_{k\to\pm\infty}|x_k(W,x/2)|=+\infty$, so there is a finite number of $(x/2)$-extrema on every compact set,
and in particular on $[x_0(W,x),x_1(W,x)]$.
Now, we can denote the $(x/2)$-extrema in this interval by $x_0(W,x)=x_{K_0}(W,x/2)<\dots<x_{K_1}(W,x/2)=x_1(W,x)$ for some integers $K_0<K_1$.

%Let $i\in\{K_0+1,\dots,K_1-1\}$.
Assume that $K_1>K_0+1$, and let $i\in\{K_0+1,\dots,K_1-1\}$.
We define $H_i:=\sup\{y>0,\ x_i(W,x/2) \text{ is an }y\text{-extremum}\}$. Assume for example that $x_i(W,x/2)$ is an $(x/2)$-minimum and that $x_0(W,x)$ is an $x$-minimum. There exists an increasing sequence $(y_n)_n$, converging to $H_i$ as $n\to+\infty$, and such that for every $n\in\N$,  $x_i(W,x/2)$ is an $y_n$-extremum, and so an $y_n$-minimum.
So, $W$ being continuous, there exist $\alpha_n<x_i(W,x/2)<\beta_n$ such that
$$
W[x_i(W,x/2)]=\inf_{[\alpha_n,\beta_n]}W, \quad  W(\alpha_n)= W[x_i(W,x/2)]+y_n=W(\beta_n).
$$
Since $x_0(W,x)<x_i(W,x/2)<x_1(W,x)$, $x_i(W,x/2)$ is not an $x$-extremum, so $x\geq H_i\geq y_n$.
If $\alpha_n<x_0(W,x)$, then $W[x_i(W,x/2)]\leq W[x_0(W,x)]$ so $x_i(W,x/2)$ would be
an $x$-minimum, which is not the case, so $\alpha_n\in[x_0(W,x),x_1(W,x))$.
If $W(\beta_n)\leq W[x_1(W,x)]$ and $\beta_n>x_1(W,x)$,
we can replace $\beta_n$ by another $\beta_n\leq x_1(W,x)$.
If $W(\beta_n)> W[x_1(W,x)]$ and $\beta_n>x_1(W,x)$,
we would have $W(\alpha_n)=W(\beta_n)>W[x_1(W,x)]$, which is the supremum
of $W$ in $[x_0(W,x),x_1(W,x)]$, and this is not possible.
%If $W(\beta_n)> W[x_1(W,x)]$ and $\beta_n>x_1(W,x)$,
%we would have $\beta_n>x_2(W,x)$,
%then $W[x_i(W,x/2)]\leq W[x_2(W,x)]\leq W[x_1(W,x)]-x<W(\beta_n)-x\leq W(\beta_n)-y_n=W[x_i(W,x/2)]$,
%which is not possible.
Hence $(\alpha_n,\beta_n)$ belongs to the compact
$[x_0(W,x),x_1(W,x)]^2$, thus there exists a strictly increasing sequence $n_p$ and $(\alpha,\beta)\in\R^2$ such that
$(\alpha_{n_p},\beta_{n_p})_{p\to+\infty}(\alpha,\beta)$.
By continuity of $W$, $W[x_i(W,x/2)]=\inf_{[\alpha,\beta]} W$,
and $W(\alpha)= W[x_i(W,x/2)]+H_i=W(\beta)$.
Hence $x_i(W,x/2)$ is an $H_i$-minimum.
Since $x_i(W,x/2)$ is not an $x$-extremum,
this gives $H_i<x$. The other cases are treated similarly.

Now, let $H_x':=\max_{K_0<i<K_1}H_i$; we have $x/2\leq H'_x<x$. For $y\in(H_x',x)$, the only possible $y$-extrema in
$(x_0(W,x),x_1(W,x))$ are the $(x/2)$-extrema, that is the $x_i(W,x/2)$, $K_0<i<K_1$,
but they are not $y$-extrema since $y>H_i$. So,
there is no $y$-extrema in $(x_0(W,x),x_1(W,x))$, and then $x_0(W,y)=x_0(W,x)$ and $x_1(W,y)=x_1(W,x)$,
for every $y\in(H_x',x)$.
This is also true with $H_x'=x/2$ in the case $K_1=K_0+1$.
Hence in every case, for every $x>0$, there exists $H_x''<x$ such that the functions $x_0(W,.)$ and $x_1(W,.)$ are constant on $[H_x'',x]$,
and consequently, they are left-continuous.
More generally, we prove similarly that for all $n\in\N^*$, there exists $K_{x,n}\in(0,x)$ such that all
the functions $x_k(W,.)$, $k\in\{-n,\dots, n\}$ are constant on $[K_{x,n},x]$.
Hence all the functions $x_k(W,.)$, $H(T_k(.))$ and $e(T_k(.))$, $k\in\Z$ are left-continuous.
$\hfill \Box$

\noindent {\bf Proof of Fact \ref{FactExcessHeight0}:}
Let $c>0$.
Assume that we are on $\{W\in\mathcal{W}\}$, and let $x>0$. We saw in  Lemma \ref{LemmaLeftcontinuity} that
there exists an interval $[y,x]$ with $0<y<x$ such that $x_0(W,.)$ and $x_1(W,.)$ are constant on this interval, and so is $b(.)$, therefore
$b(.)$ does not change its sign on $[y,x]$.

Define $H_{p,q}:=\big|\sum_{k=p}^{q-1}(-1)^k H(T_k(c))\big|$ for $p<q$ and
$\mathcal H :=\{\forall p<q\leq r<s, \ H_{p,q}\neq H_{r,s}\}\cap \{W\in \mathcal{W}\}$.
Since the r.v. $H(T_k(c))$, $k\in\Z$ are independent (see \cite{NP} Proposition of Section 1)
and have a density (see \cite{Cheliotis} (8) p. 1768 and (11) p. 1770), it follows that
the r.v. $H_{p,q}-H_{r,s}$, $p<q\leq r<s$ also have densities, thus $\eta(\mathcal H)=1$.
Moreover, for every trajectory $W\in\mathcal W$, every $x\geq c$ and $m<n$, there exist
$p<q \leq r<s$ such that $H(T_m(x))=H_{p,q}$ and $H(T_n(x))=H_{r,s}$.
Consequently, on $\mathcal H$, for every $x\geq c$, all the $H(T_i(x))$, $i\in\Z$ are different.

%\noindent
Now, assume we are on $\mathcal H$. Let $x\geq c$. The $e(T_i(x))$, $i\in\{-3, \dots 3\}$ are all different,
so for $\e>0$ small enough, at most one of them is less than $\e$.
As was shown in the proof of Lemma 2 of Cheliotis (\cite{Cheliotis} p. 1772), for such $\e>0$, $b(x)$ and $b(x+\e)$ have
a different sign iff $e[T_0(x)]<\e$.
So, if $e(T_0(x))>0$ (resp. $e(T_0(x))=0$), there exists $\e>0$ such that the sign of $b(.)$ in $(x,x+\e]$ is the sign of $b(x)$ (resp. of $-b(x)$).

%So there exists a non empty interval $[x,x+\e]$ without any change of sign of $b$ for some $\e>0$
%iff $e(T_0(x))>0$.

\noindent
Hence on $\mathcal H$ there is a change of sign of $b$ at $x$ iff $e(T_0(x))=0$, which proves Fact \ref{FactExcessHeight0}.
\hfill$\Box$

\subsection{Proof of Lemma \ref{LemmaProbaStrongChangeOfSign}}
We consider a two-sided Brownian motion $W$ defined on a probability space $(\Omega, \mathcal{A}, \eta)$.
We know that $\eta(\mathcal{H}\cap\{W\in\mathcal{W}\})=1$. This enables us to replace, in the rest of the paper, $\Omega$ by $\Omega \cap \mathcal{H}\cap\{W\in\mathcal{W}\}$.

%We first introduce some notation.
We denote by $\F_x$ the completion of the $\sigma$-field
$\sigma\big(W(s)\un_{\{x_0(W,x)\leq s \leq x_1(W,x)\}},s\in\R\big)$
for $x>0$,
and by $\mathcal{F}_0$ and $\mathcal{F}_\infty$ the completions of $\sigma(\emptyset)$ and $\sigma(W(s),s\in\R)$ respectively.
%and we set $\mathcal{F}_0=\sigma(\emptyset)$ and denote by $\mathcal{F}_\infty$ the completion of $\sigma(W(s),s\in\R)$.
For $0<y\leq x$, $[x_0(W,y),x_1(W,y)]\subset [x_0(W,x),x_1(W,x)]$ and
$x_0(W,y)$ and $x_1(W,y)$ are $\mathcal{F}_x$-measurable
(which we prove in details in Lemma \ref{LemmaMesurabilitex0x1} in Subsection \ref{SubSectAppendix} Appendix), so
$\mathcal{F}_y\subset \mathcal{F}_x$. Hence $(\mathcal{F}_x)_{x\geq 0}$ is a filtration.
%However this filtration is not right-continuous and then does not satisfy the so-called usual conditions,
Notice that $W$ is not adapted to
$(\mathcal{F}_x)_{x\geq 0}$. Moreover, for $k\in\Z$,  $x\mapsto e[T_k(x)]$ is left-continuous by Lemma \ref{LemmaLeftcontinuity},
but it is not right-continuous,
and $(\mathcal{F}_x)_{x\geq 0}$ is not the natural filtration of one of these processes.
We now give an elementary proof of Lemma \ref{LemmaProbaStrongChangeOfSign}. We start with the following lemma.

\begin{lemma}\label{eqLemmaStoppingTimes}
For every $k\geq 1$, $X_k$ is a $(\mathcal{F}_x)_{x\geq 0}$-stopping time.
\end{lemma}

\noindent{\bf Proof:}
%Since it is not straightforward
Instead of trying to prove
whether the filtration $(\mathcal{F}_x)_x$ is right-continuous, we give an elementary proof.
Notice that $e[T_0(y)]=(\sup_\R-\inf_\R)(W\un_{[x_0(W,y),x_1(W,y)]})-y$ is $\mathcal{F}_y$-measurable for every $y>0$,
that means, the processes $(e[T_0(y)])_y$ and then $(H[T_0(y)])_y$ are adapted to the filtration $(\mathcal{F}_y)_y$.
%and that $y\mapsto e(T_0(y))$ is left-continuous.
Moreover, the function $e[T_0(.)]$ has a jump at $y\in[c,x]$ if and only if $x_0(W,y)$ or $x_1(W,y)$ is a $y$-extremum but is not a $z$-extremum for $z>y$,
and in this case the number of $z$-extrema in $[x_0(W,x),x_1(W,x)]$ decreases by at least $1$ between $z=y$ and every $z>y$.
%Moreover, the function $e[T_0(.)]$ has a jump at $y$ if and only if $x_0(W,y)$ or $x_1(W,y)$, is a $y$-extremum but is not a $z$-extremum for $z>y$.
So, the number of discontinuities of $e[T_0(.)]$ in $[c,x]$ is less than the number of $c$-extrema in $[x_0(W,x),x_1(W,x)]$, which is finite
on $\{W\in\mathcal{W}\}$.

\noindent Hence,  the process $e(T_0(.))$ is left-continuous with  a finite number of discontinuities in $[c,x]$, is nonnegative, and it is strictly decreasing between two consecutive discontinuities and then has right limits.
Moreover on $\{W\in\mathcal{W}\}$, $H(T_0(.))$ is nondecreasing and so only has positive jumps, and then $e(T_0(.))$ also has only positive jumps. As a consequence, $e(T_0(.))$, which is left-continuous with right limits,
is lower semi-continuous on $(0,+\infty)$.

\noindent
Recalling that $\{X_1\leq x\}=\{\exists y\in[c,x],\ e[T_0(y)]=0\}$ by the proof of Fact \ref{FactExcessHeight0} since $\Omega\subset \mathcal H$,
we claim that for $x\geq c$,
\begin{eqnarray}
    \{X_1\leq x\}
& = & \label{eqLigne1X1}
    \cap_{p\in\N^*}\{\exists y\in[c,x],\ e[T_0(y)]<1/p\}
    \\
%& = &
%    \cap_{p\in\N^*}\{\exists y\in([c,x]\cap\mathbb{Q})\cup\{c\},\ e(T_0(y))<1/p\}\\
& = & \label{eqLigne2X1}
    \cap_{p\in\N^*}\cup_{y\in([c,x]\cap\mathbb{Q})\cup\{c\}}\{e[T_0(y)]<1/p\}
%    \cap \mathcal{W}
    .
\end{eqnarray}
Indeed for the first line, inclusion $\subset$ is clear. For the inclusion $\supset$,
on the event in RHS of \eqref{eqLigne1X1}, where RHS stands for right hand side,
there is a sequence $y_n\in[c,x],\ n\in\N^*$ such that $e[T_0(y_n)]<1/n$ for $n\in\N^*$.
Since $[c,x]$ is compact, there exists a subsequence
$(y_{p_n})_n$,
%such that $e(T_0(y_{p_n}))<1/{p_n}$,
which converges to an $y\in[c,x]$.
Hence, $0\leq e[T_0(y)]\leq \liminf_{n\to+\infty}e[T_0(y_{p_n})] =0$
by lower semi-continuity,
which proves the inclusion.
For line \eqref{eqLigne2X1}, inclusion
$(\text{RHS of }\eqref{eqLigne1X1})\supset(\text{RHS of }\eqref{eqLigne2X1})$ is clear, whereas inclusion $\subset$ follows from the
left-continuity of $e(T_0(.))$.

\noindent
Hence $\{X_1\leq x\}\in\mathcal{F}_x$ for every $x\geq c$, and $\{X_1\leq x\}=\emptyset\in\mathcal{F}_x$ for $0\leq x<c$,
so $X_1$ is a $(\mathcal{F}_x)_{x\geq 0}$-stopping time. Let $k\geq 1$.
Since $\lim_{u\to X_k,u>X_k}e[T_0(u)]>0$ because $e[T_0(X_k)]=0$ and so there is a positive jump at $x$ for $e[T_0(.)]$,
we show similarly that for $x\geq c$,
\begin{eqnarray*}
    \{X_{k+1}\leq x\}
    %\cap\mathcal{W}
& = &
    %\mathcal{W}\cap
    \{X_k<x\}\cap \cap_{p\in\N^*}\{\exists y\in(X_k,x],\ e(T_0(y))<1/p\}\\
%& = &
%    \{X_k<x\}\cap \cap_{p\in\N^*}\cup_{y\in((X_k,x]\cap\mathbb{Q})}\{e(T_0(y))<1/p\}\\
& = &
    %\mathcal{W}\cap
    \{X_k<x\}\cap \cap_{p\in\N^*}\cup_{y\in((c,x]\cap\mathbb{Q})}
    [\{y> X_k\}\cap\{e(T_0(y))<1/p\}].
\end{eqnarray*}
Hence it follows by induction that
$X_k$ is a $(\mathcal{F}_x)_{x\geq 0}$-stopping time for every $k\geq 1$.
\hfill$\Box$

\noindent
We can then consider the $\sigma$-fields $\mathcal{F}_{X_k}$ for $k\geq 1$.

We now fix $k\geq 1$. First, we notice that
$A_{k+1,a,c}=A_{k+1,a,c}^+\cup A_{k+1,a,c}^-$, where
$A_{k+1,a,c}^+:=A_{k+1,a,c}\cap\{b(c)>0\}$ and $A_{k+1,a,c}^-:=A_{k+1,a,c}\cap\{b(c)\leq 0\}$. We start with $A_{k+1,a,c}^+$,
and notice that
\begin{equation}
%A_{k+1,a,c}^-
%& := &
%    A_{k,a,c}\cap\{b(X_{1})\leq 0\}
%    \cap\left[\{e(T_{-1}(X_{2k+2}))<a X_{2k+2}\}\cup\{e(T_{1}(X_{2k+2}))<a X_{2k+2}\}\right],\nonumber\\
A_{k+1,a,c}^+
 =
    A_{k,a,c}\cap\{b(X_{1})>0\}
    \cap\left[\{e(T_{-1}(X_{2k+1}))<a X_{2k+1}\}\cup\{e(T_{1}(X_{2k+1}))<a X_{2k+1}\}\right].
    \label{eqFormuleA+}
\end{equation}
Let $n_0\in\N^*$. We define a sequence $(R_n)_{n\geq n_0}$ by induction as follows:
%such that $R_n\leq X_{2k+1}$ for every $n\in\N$,
%$R_n$ is measurable with respect to $\mathcal{F}_{X_{2k+1}}$,
%and $R_n\to_{n\to+\infty} X_{2k+1}$ on $\{R_{n_0}\neq 0\}$.
%We set $R_0=X_{2k}$.
% We first define
\begin{eqnarray*}
    R_{n_0}
& := &
    2^{-n_0}(\lfloor 2^{n_0}X_{2k}\rfloor +1)\un_{\{X_{2k+1}>2^{-n_0}(\lfloor 2^{n_0}X_{2k}\rfloor+1)\}
    %\cap \mathcal{W}
    },\\
    R_{n}
& := &
    2^{-n}\lfloor 2^nH[T_0(R_{n-1})]\rfloor
    \un_{\{X_{2k+1}>2^{-n_0}(\lfloor 2^{n_0}X_{2k}\rfloor+1)\}
    %\cap \mathcal{W}
    },\qquad n>n_0.
\end{eqnarray*}
In particular, we have $c\leq X_{2k}< R_{n_0}< X_{2k+1}$ on
$B_{k+1,a,c}^{+,n_0}:=\{X_{2k+1}>2^{-n_0}(\lfloor 2^{n_0}X_{2k}\rfloor+1)\}=\{R_{n_0}\neq 0\}$.
Moreover $R_{n}\in (2^{-n}\N)$ for all $n\geq n_0$.
We have
$R_n\leq H[T_0(R_{n-1})]\leq R_n+2^{-n}$ on  $B_{k+1,a,c}^{+,n_0}$
and $R_{n}=0$ on  $(B_{k+1,a,c}^{+,n_0})^c$ for $n\geq n_0$.
We now prove the two following lemmas:

\begin{lemma}\label{LemmaLimiteRn}
The sequence $(R_n)_{n\geq n_0}$ is nondecreasing. It converges a.s. to a r.v. $R_\infty$, and
$$
R_\infty=X_{2k+1}\un_{B_{k+1,a,c}^{+,n_0}}.
%    \quad \text{ on }
%B_{k+1,a,c}^{+,n_0}.
$$
\end{lemma}

\noindent{\bf Proof:}
Since $H[T_0(x)]\geq x$ for every $x\geq 0$ and
$(2^n R_{n-1})\in \N$ for $n>n_0$,
we get on $B_{k+1,a,c}^{+,n_0}$,
$$
    R_{n-1}
=
    2^{-n}\lfloor 2^{n}R_{n-1}\rfloor
\leq
    2^{-n}\lfloor 2^{n}H[T_0(R_{n-1})]\rfloor
=
    R_n, \qquad n>n_0.
$$
So, $(R_n)_{n\geq n_0}$ is a nondecreasing sequence on $B_{k+1,a,c}^{+,n_0}$, and also on
$(B_{k+1,a,c}^{+,n_0})^c$ on which $R_n=0$ for every $n\geq n_0$. Hence, it tends a.s. to
$R_\infty:=\lim_{n\to+\infty} R_n\in[R_{n_0},+\infty]$.

Let $n\geq n_0+1$. If $R_{n-1}<x<R_{n}$, then $R_n\neq 0$ and we have
\begin{equation}\label{EqPositifEntrelesRn}
    e[T_0(x)]=H[T_0(x)]-x
\geq
    H[T_0(R_{n-1})]-x
%=
%2^{-n} 2^n H[T_0(R_{n-1})]-x
%\geq
%2^{-n} \lfloor 2^n H[T_0(R_{n-1})]\rfloor -x
%$$
%that is,
%$$
%e[T_0(x)]
\geq R_n-x>0.
\end{equation}
Assume that $R_{n_0}\neq 0$ and that there exists $n\geq n_0$ such that $e[T_0(R_{n})]=0$, and let $n_1$ denote the smallest such $n$. Then,
$H[T_0(R_{n_1})]=R_{n_1}+e[T_0(R_{n_1})]=R_{n_1}$, so
$$
    R_{n_1+1}
=
    2^{-{(n_1+1)}}\lfloor 2^{n_1+1}H[T_0(R_{n_1})]\rfloor
=
    2^{-{(n_1+1)}}\lfloor 2^{n_1+1}R_{n_1}\rfloor
=
    R_{n_1}
$$
since $R_{n_1}\in 2^{-n_1}\N$.
We prove similarly by induction that $R_n=R_{n_1}$ for every $n\geq n_1$,  so $R_\infty=R_{n_1}$ and then $e[T_0(R_\infty)]=0$.
Moreover, by \eqref{EqPositifEntrelesRn}, $e(T_0(.))>0$ on $(R_{n_0},R_\infty)$.
Furthermore we know that on $B_{k+1,a,c}^{+,n_0}$, $X_{2k}<R_{n_0}<X_{2k+1}$,
so $e(T_0(.))>0$ on $(X_{2k},R_{n_0}]$ by Fact \ref{FactExcessHeight0} and then on
$(X_{2k},R_\infty)$.
Hence $R_\infty=\inf\{x>X_{2k},\ e[T_0(x)]=0\}=X_{2k+1}$ in this case.
% on $B_{k+1,a,c}^{+,n_0}$.
%Since $e[T_0(X_{2k+1})]=0$ and $X_{2k}< X_{2k+1}$, this yields  $R_\infty\leq X_{2k+1}<\infty$.

\noindent Else, assume that $R_{n_0}\neq 0$ and $e[T_0(R_{n})]\neq 0$ for every $n\geq n_0$.
%Then for every $n> n_0$, $R_{n-1}<R_n$ and $e[T_0(.)]>0$ on $(R_{n-1},R_n)$. Hence, $e[T_0(.)]>0$ on $[R_{n_0},R_\infty)$.
Then $(R_n)_{n\geq n_0}$ is a nondecreasing sequence such that  $e[T_0(.)]>0$ on each interval $(R_{n-1},R_n)$, $n>n_0$
by \eqref{EqPositifEntrelesRn}, and then $e[T_0(.)]>0$ on $[R_{n_0},R_\infty)$.
As in the previous case, we get $e[T_0(.)]>0$ on $(X_{2k},R_\infty)$.
Since $e[T_0(X_{2k+1})]=0$ and $X_{2k}< X_{2k+1}$, this yields  $R_\infty\leq X_{2k+1}<\infty$.
%\noindent
%Hence, combining these two cases, we have $R_\infty\leq X_{2k+1}<\infty$
% on $B_{k+1,a,c}^{+,n_0}$.

\noindent
Moreover in this case, as explained before Lemma \ref{LemmaLimiteRn},
$0< e[T_0(R_{n-1})]=H[T_0(R_{n-1})]-R_{n-1}\leq R_n+2^{-n}-R_{n-1}\to_{n\to+\infty} 0$ a.s.,
because $R_n\to_{n\to+\infty}R_\infty$.
Since $e[T_0(.)]$ is a left-continuous function on $\mathcal{W}$  and $(R_n)_n$ is nondecreasing
and converging to $R_\infty<\infty$, this gives
$e[T_0(R_\infty)]=\lim_{n\to+\infty}e[T_0(R_{n-1})]= 0$.
As in the previous case, we conclude that $R_\infty=X_{2k+1}$.
Since $R_n=0$ $\forall n\geq n_0$ if $R_{n_0}=0$, that is, on $\big(B_{k+1,a,c}^{+,n_0}\big)^c$, this proves the lemma.
\hfill$\Box$

\begin{lemma}\label{LemmaMeasurabilityRn}
For all $n\geq n_0$,
\begin{equation}\label{eqInductionHypotheseRn}
\forall m\in\N^*,\qquad \{R_{n}=m 2^{-n}\}\in\mathcal{F}_{m2^{-n}}.
\end{equation}
\end{lemma}

\noindent{\bf Proof:}
We prove this lemma by induction. We start with $R_{n_0}$, and  observe that for $m\in\N^*$,
\begin{eqnarray*}
\{R_{n_0}=m 2^{-n_0}\}
& = &
\{X_{2k+1}>2^{-n_0}(\lfloor 2^{n_0}X_{2k}\rfloor+1)\}\cap\{\lfloor 2^{n_0}X_{2k}\rfloor=m-1\}\\
& = &
\{X_{2k+1}>m2^{-n_0}\}\cap\{(m-1)2^{-n_0}\leq X_{2k}<m2^{-n_0}\},
\end{eqnarray*}
which belongs to $\mathcal{F}_{m2^{-n_0}}$
since $X_{2k}$ and $X_{2k+1}$ are $(\mathcal{F}_x)_{x\geq 0}$-stopping times by Lemma \ref{eqLemmaStoppingTimes}.
%$\{R_{n_0}=m 2^{-n_0}\}
%\in\mathcal{F}_{m2^{-n_0}}$ for every $m\in\N^*$.
This gives \eqref{eqInductionHypotheseRn} for $n=n_0$.
Now, assume that \eqref{eqInductionHypotheseRn} is true for some $n\geq n_0$. Then for $m\in\N^*$,
\noindent
%Assume that for some , $\{R_{n}=m 2^{-n}\}\in\mathcal{F}_{m2^{-n}}$ for every $m\in\N^*$.
\begin{eqnarray*}
    \{R_{n+1}=m 2^{-(n+1)}\}
& = &
    \{\lfloor 2^{n+1}H[T_0(R_{n})]\rfloor=m\}
    \cap B_{k+1,a,c}^{+,n_0}\\
%& = &
%    \{\lfloor 2^{n+1}H[T_0(R_{n})]\rfloor=m\}\\
& = &
    \cup_{p\in\N^*}\{R_n=p2^{-n}, \lfloor 2^{n+1}H[T_0(R_{n})]\rfloor=m\}
    \\
& = &
    \cup_{p\in\N^*, p2^{-n}\leq m2^{-(n+1)}}
    [\{ R_n=p2^{-n}\}\cap\{\lfloor 2^{n+1}H[T_0(p2^{-n})]\rfloor=m\}].
\end{eqnarray*}
The second equality comes from  $\{R_n \neq 0\}=\{R_{n+1}\neq 0\}=B_{k+1,a,c}^{+,n_0}$,
which itself is a consequence of $R_n\geq R_{n_0}>X_{2k}\geq c>0$
on $B_{k+1,a,c}^{+,n_0}$.
The third one is a consequence of
$R_n\leq R_{n+1}$. 
%or $x\leq H[T_0(x)]$.
%Due to our induction hypothesis,  for every  $p\in\N^*$ such that
If $0<p2^{-n}\leq m2^{-(n+1)}$,
our induction hypothesis gives
$\{ R_n=p2^{-n}\} \in\mathcal{F}_{p2^{-n}}\subset\mathcal{F}_{m2^{-(n+1)}}$,
%Moreover, $e(T_0(y))$ is $\mathcal{F}_y$-measurable for every $y>0$, so for such $p$,
and $\{\lfloor 2^{n+1}H[T_0(p2^{-n})]\rfloor=m\}\in\mathcal{F}_{p2^{-n}}\subset\mathcal{F}_{m2^{-(n+1)}}$ since $(H[T_0(y)], \ y\geq 0)$
is adapted to $(\mathcal F_y)_{y\geq 0}$. Consequently,
$\big\{R_{n+1}=m 2^{-(n+1)}\big\}\in\mathcal{F}_{m2^{-(n+1)}}$ for every $m\in\N^*$, which ends the induction.
\hfill$\Box$

%Recalling that $X_{2k+1}$ is the first zero of $e[T_0(.)]$ after $X_{2k}$, and since $X_{2k}< R_\infty$, $e[T_0(R_\infty)]=0$
%yields  $X_{2k+1}\leq R_\infty$. Since we also have in any case $R_\infty\leq X_{2k+1}$
%on $B_{k+1,a,c}^{+,n_0}$, we get
%$$
%R_\infty=X_{2k+1}
%    \quad \text{ on }
%B_{k+1,a,c}^{+,n_0}.
%$$
%Now, we have,
%\begin{equation*}
%    A_{k+1,a,c}^+
%=
%    A_{k,a,c}\cap\{b(X_{1})>0\}
%    \cap\left[\{e(T_{-1}(X_{2k+1}))<a X_{2k+1}\}\cup\{e(T_{1}(X_{2k+1}))<a X_{2k+1}\}\right],
%\end{equation*}
In view of \eqref{eqFormuleA+}, we define for $n\geq n_0$,
%(with $e(T_i(0))=0$ by convention)
$$
    C_{k+1,a,c}^{+,n}
:=
       A_{k,a,c}\cap\{b(X_{1})>0\}
    \cap\left[\{e[T_{-1}(R_n)]<a R_n\}\cup\{e[T_{1}(R_n)]<a R_n\}\right].
$$
Assume that we are on $B_{k+1,a,c}^{+,n_0}\cap A_{k+1,a,c}^+$.
There exists $i\in\{-1,1\}$ such that $e(T_{i}(X_{2k+1}))<a X_{2k+1}$,
that is $H[T_{i}(X_{2k+1})]<(a+1) X_{2k+1}$.
On the one hand,
$R_n\to_{n\to+\infty} X_{2k+1}$, $R_n\leq X_{2k+1}$
by Lemma \ref{LemmaLimiteRn},
%We have $R_n\to_{n\to+\infty} X_{2k+1}$,
%and since $e[T_i(.)]$ are left-continuous functions, we also have
%$e[T_i(R_n)]\to_{n\to+\infty} e[T_i(X_{2k+1})]$.
then by Lemma \ref{LemmaLeftcontinuity}, for $n$ large enough,
$R_n\in[K_{X_{2k+1},2},X_{2k+1}]$, then
$x_j(W, R_n)=x_j(W, X_{2k+1})$ for $-1\leq j\leq 2$ and so $H[T_i(R_n)]=H[T_i(X_{2k+1})]$.
On the other hand, $(a+1)(X_{2k+1}-R_n)$ tends to $0$ as $n\to\infty$ by Lemma \ref{LemmaLimiteRn}
and then is strictly less than
$(a+1)X_{2k+1}-H[T_i(X_{2k+1})]>0$ for $n$ large enough. So for large $n$,
\begin{equation*}
    H[T_i(R_n)]-(a+1)R_n
%& = &
%    H[T_i(X_{2k+1})]-(a+1)R_n\\
 =  (a+1)(X_{2k+1}-R_n)\\
 -[(a+1)X_{2k+1}-H[T_i(X_{2k+1})]]
 <  0,
\end{equation*}
and so 
$
e[T_i(R_n)]<a R_n
$.
Then for large $n$, $\un_{B_{k+1,a,c}^{+,n_0}\cap C_{k+1,a,c}^{+,n}}=1$.
%(virer $B_{k+1,a,c}^{+,n_0}$ ?)

\noindent
Hence, in every case,
$\un_{B_{k+1,a,c}^{+,n_0}\cap A_{k+1,a,c}^+}\leq \liminf_{n\to+\infty}\un_{B_{k+1,a,c}^{+,n_0}\cap C_{k+1,a,c}^{+,n}}$.
Then by Fatou's lemma,
\begin{equation}\label{eqFatou}
    \eta\big(B_{k+1,a,c}^{+,n_0}\cap A_{k+1,a,c}^+\big)
\leq
    \int_\Omega\Big(\liminf_{n\to+\infty}\un_{B_{k+1,a,c}^{+,n_0}\cap C_{k+1,a,c}^{+,n}}\Big)\dd \eta
\leq
    \liminf_{n\to+\infty}\eta\big({B_{k+1,a,c}^{+,n_0}\cap C_{k+1,a,c}^{+,n}}\big).
\end{equation}
Let $n\geq n_0$. We now have to estimate, recalling that 
$R_n\geq R_{n_0}>X_{2k}\geq c>0$
on $B_{k+1,a,c}^{+,n_0}$,
\begin{equation}
    \eta\big({B_{k+1,a,c}^{+,n_0}\cap C_{k+1,a,c}^{+,n}}\big)
%& = &
%    \eta({B_{k+1,a,c}^{+,n_0}\cap C_{k+1,a,c}^{+,n}}\cap\cup_{m=c 2^n}^\infty\{R_n=m 2^{-n}\})
 =
    \sum_{m\in\N, \ m\geq c 2^n}\eta\big({B_{k+1,a,c}^{+,n_0}\cap C_{k+1,a,c}^{+,n}}\cap\{R_n=m 2^{-n}\}\big)
%& = &
%    \sum_{m=c 2^n}^\infty\eta({B_{k+1,a,c}^{+,n_0}\cap C_{k+1,a,c}^{+,n}}\cap\{{\{m 2^{-n}\leq  H[T_0(R_{n-1})]   < (m+1)2^{-n}\}}\})
    \label{eqProbasSumBC}.
\end{equation}
For $m\geq c 2^n$, we have, since  $m>0$ and then $\{R_n=m 2^{-n}\}\subset B_{k+1,a,c}^{+,n_0}$,
\begin{align}
&
    \eta\big({B_{k+1,a,c}^{+,n_0}\cap C_{k+1,a,c}^{+,n}}\cap\{R_n=m 2^{-n}\}\big)
    \nonumber\\
& =
    \eta\big({C_{k+1,a,c}^{+,n}}\cap\{R_n=m 2^{-n}\}\big)
    \nonumber\\
& =
    \eta([A_{k,a,c}\cap\{b(X_{1})>0\}
    \cap\cup_{i=\pm1}\{e[T_{i}(R_n)]<a R_n\}]
    \cap\{R_n=m 2^{-n}\})\nonumber\\
& =
    \eta(A_{k,a,c}\cap\{b(c)>0\}\cap\{X_{2k}<m 2^{-n}\}
    \cap\{R_n=m 2^{-n}\}
    \cap\cup_{i=\pm1}\{e[T_{i}(m 2^{-n})]<a m 2^{-n}\}),
        \label{eqProbaBCR}
\end{align}
where the last equality comes from $X_{2k}< R_{n_0}\leq R_n$ on $\{R_n>0\}=B_{k+1,a,c}^{+,n_0}$.

\noindent
For $\ell\geq 1$, we have on $\{X_\ell< x\}$, $[x_{-1}(W,X_\ell),x_2(W,X_\ell)]\subset [x_0(W,x),x_1(W,x)]$
since $x_0(W,X_\ell)$ and $x_1(W,X_\ell)$ are not $x$-extrema on $\mathcal H$ due to $H[T_0(X_\ell)]=X_\ell<x$. Hence,
the random variables $e[T_i(X_\ell)]$, $i\in\{-1,1\}$ are measurable with respect to
$\mathcal{F}_{X_{\ell}+}=\{A\in\mathcal{F}_\infty,\ \forall x\geq 0,\ A\cap\{X_{\ell}<x\}\in\mathcal{F}_x\}$
(this is proved in details in Lemma \ref{LemmaMusurabiliteFX+} in Subsection \ref{SubSectAppendix} Appendix).
As a consequence, $A_{k,a,c}\in \mathcal{F}_{X_{2k}+}$ for every $k\geq 1$.
%{\bf We assume for the moment that} $A_{k,a,c}\in \mathcal{F}_{X_{2k}+}$,
which gives in particular
$[A_{k,a,c}\cap\{X_{2k}<m 2^{-n}\}]\in\mathcal{F}_{m 2^{-n}}$ for every $m\in\N$.

\noindent
Moreover, let $m\in\N$ such that $c\leq m 2^{-n}$.
We have  $\{b(c)>0\}\in\mathcal{F}_c\subset \mathcal{F}_{m 2^{-n}}$.
%otherwise $\{b(c)>0\}\cap\{X_{2k}<m 2^{-n}\}=\emptyset\in\mathcal{F}_{m 2^{-n}}$.
Since $\{R_n=m 2^{-n}\}\in\mathcal{F}_{m2^{-n}}$ by Lemma \ref{LemmaMeasurabilityRn}, we get
$[A_{k,a,c}\cap\{X_{2k}<m 2^{-n}\}\cap\{b(c)>0\}\cap\{R_n=m 2^{-n}\}]\in\mathcal{F}_{m2^{-n}}$.
%But  $\cup_{i=\pm1}\{e(T_{i}(m 2^{-n}))<a m 2^{-n}\}$ is independent of the slope $T_0(m 2^{-n})$ and then of
%$\mathcal{F}_{m2^{-n}}$ by Neveu, so
But  $e[T_{1}(m 2^{-n})]$, $e[T_{-1}(m 2^{-n})]$ and $\mathcal{F}_{m2^{-n}}$ are independent by Neveu et al. (\cite{NP}, Proposition of Section 1), so
\begin{eqnarray}
    \text{RHS of }\eqref {eqProbaBCR}
& = &
    \eta[A_{k,a,c}\cap\{b(c)>0\}\cap\{X_{2k}<m 2^{-n}\}
    \cap\{R_n=m 2^{-n}\}]\nonumber\\
&&
    \times
    \eta(\cup_{i=\pm1}\{e[T_{i}(m 2^{-n})]<a m 2^{-n}\})
    \nonumber\\
%& = &
%    (1-e^{-2a})
%    \eta[A_{k,a,c}\cap\{b(c)>0\}\cap\{X_{2k}<m 2^{-n}\}
%    \cap\{R_n=m 2^{-n}\}]\nonumber\\
& = &
    \big(1-e^{-2a}\big)
%    \sum_{m=c 2^n}^\infty
    \eta[A_{k,a,c}\cap\{b(c)>0\}
    \cap\{R_n=m 2^{-n}\}]
    \label{eqProbaBCR2}
\end{eqnarray}
since $e[T_{i}(m 2^{-n})]/(m 2^{-n})$, $i\neq 0$, are independent exponential r.v. with mean
$1$
(also by Neveu et al. \cite{NP}, prop. 1)
%(by Cheliotis \cite{Cheliotis} formula (8))
and $X_{2k}< R_{n}$ on $\{R_n\neq 0\}$.
%\begin{equation*}
%    \eta(\cup_{i=\pm1}\{e(T_{i}(m 2^{-n}))<a m 2^{-n}\})
% =
%    1-\eta(\cap_{i=\pm1}\{e(T_{i}(m 2^{-n}))\geq a m 2^{-n}\})\\
% =
%    1- \eta[e(T_{1}(m 2^{-n}))/(m 2^{-n})\geq a ]^2\\
% =
%    1-e^{-2a}.
%\end{equation*}
So, \eqref{eqProbasSumBC}, \eqref{eqProbaBCR} and \eqref{eqProbaBCR2},
%$R_n\in (2^{-n} \N)$ and $R_n\geq X_{2k}\geq c$ on $\{R_n\neq 0\}=B_{k+1,a,c}^{+,n_0}$
give
\begin{eqnarray*}
%&&
    \eta\big({B_{k+1,a,c}^{+,n_0}\cap C_{k+1,a,c}^{+,n}}\big)
%    \nonumber\\
%& = &
%    \sum_{m=c 2^n}^\infty\eta({B_{k+1,a,c}^{+,n_0}\cap C_{k+1,a,c}^{+,n}}\cap\{R_n=m 2^{-n}\})\nonumber\\
%& = &
%    (1-e^{-2a})
%    \sum_{m=c 2^n}^\infty
%    \eta[A_{k,a,c}\cap\{b(c)>0\}\cap\{X_{2k}<m 2^{-n}\}
%    \cap\{R_n=m 2^{-n}\}]\\
& = &
    (1-e^{-2a})
    \sum_{m\in\N,\ m\geq c 2^n}
    \eta[A_{k,a,c}\cap\{b(c)>0\}
    \cap\{R_n=m 2^{-n}\}]\\
%& = &
%    (1-e^{-2a})
%%    \eta[A_{k,a,c}\cap\{b(X_{1})>0\}].
%    \eta[A_{k,a,c}^+\cap B_{k+1,a,c}^{+,n_0}].
& \leq &
    \big(1-e^{-2a}\big)
%    \eta[A_{k,a,c}\cap\{b(X_{1})>0\}].
    \eta\big[A_{k,a,c}^+\big].
\end{eqnarray*}
Consequently, \eqref{eqFatou} leads to
%\begin{equation*}
%    \eta(B_{k+1,a,c}^{+,n_0}\cap A_{k+1,a,c}^+)
%\leq
%%    \liminf_{n\to+\infty}\eta({B_{k+1,a,c}^{+,n_0}\cap C_{k+1,a,c}^{+,n}})
%%=
%    (1-e^{-2a})
%%    \eta[A_{k,a,c}\cap\{b(X_{1})>0\}]
%    \eta[A_{k,a,c}^+].
%\end{equation*}
%Hence,
\begin{equation*}
    \eta\big(A_{k+1,a,c}^+\big)
\leq
    \eta\big(A_{k+1,a,c}^+\cap B_{k+1,a,c}^{+,n_0}\big)+\eta\big[\big(B_{k+1,a,c}^{+,n_0}\big)^c\big]
\leq
    \big(1-e^{-2a}\big)
%    \eta[A_{k,a,c}\cap\{b(X_{1})>0\}]
    \eta\big[A_{k,a,c}^+\big]
    +
    \eta\big[\big(B_{k+1,a,c}^{+,n_0}\big)^c\big].
\end{equation*}
But $c\leq X_{2k}$ and $X_{2k+1}/X_{2k}>1$ a.s., so
$$
    \eta\big[\big(B_{k+1,a,c}^{+,n_0}\big)^c\big]
%=
%\eta[X_{2k+1}\leq 2^{-n_0}(\lfloor 2^{n_0}X_{2k}\rfloor+1)]
\leq
    \eta\big[X_{2k+1}\leq X_{2k}+2^{-n_0}\big]
\leq
    \eta\big[X_{2k+1}/X_{2k}\leq 1+2^{-n_0}/c\big]
\to_{n_0\to+\infty}
    0.
$$
% and has a density which is given by (Cheliotis \cite{Cheliotis} Corollary 2). (maybe we do not need this !).
As a consequence,
\begin{equation*}
    \eta\big(A_{k+1,a,c}^+\big)
\leq
    \big(1-e^{-2a}\big)
    \eta\big(A_{k,a,c}^+\big).
\end{equation*}
We get similarly
$
    \eta\big(A_{k+1,a,c}^-\big)
\leq
    \big(1-e^{-2a}\big)
    \eta\big(A_{k,a,c}^-\big)
$.
These two inequalities yield
$
    \eta(A_{k+1,a,c})
\leq
    \big(1-e^{-2a}\big)
    \eta(A_{k,a,c})
$. Using this last inequality, we obtain \eqref{eqLemmaProbaAka} by induction on $k$,
%$$
%    \eta(A_{k,a,c})
%\leq
%    (1-e^{-2a})^{k-1}
%    \eta(A_{1,a,c})
%$$
which proves Lemma \ref{LemmaProbaStrongChangeOfSign}.
\hfill$\Box$

\subsection{Appendix : measurability}\label{SubSectAppendix}
We fix $x>0$. We define
\begin{equation}\label{eqDefZ}
    Z(s)
=
    W(s)\un_{\{x_0(W,x)\leq s \leq x_1(W,x)\}},
\end{equation}
so that $\mathcal{F}_x$ is the completion of $\sigma(Z(s),\ s\in\R)$. For the sake of completeness, we prove in this appendix
the measurability of some random variables. We start with the following lemma, which is used before Lemma \ref{eqLemmaStoppingTimes}
to prove that $(\mathcal{F}_x)_{x\geq 0}$ is a filtration.

\begin{lemma}\label{LemmaMesurabilitex0x1}
If $0<y\leq x$, then
$x_0(W,y)$ and $x_1(W,y)$ are $\mathcal{F}_x$-measurable.
\end{lemma}

\noindent{\bf Proof:}
%First, $x_1(W,x)=\sup\{s\geq 0,\ Z(s)\neq 0\}$, so $x_1(W,x)$  is  $\mathcal{F}_x$-measurable.
%Similarly, $x_0(W,x)$  is  $\mathcal{F}_x$-measurable.
Let $0< y<x$, and
\begin{eqnarray*}
z_0=z_0(y)  :=  \inf\{t\in\R,\ Z(t)\neq 0\}=x_0(W,x),\qquad
z_\infty  :=  \sup\{t\in\R,\ Z(t)\neq 0\}=x_1(W,x).
\end{eqnarray*}
This already proves that $x_0(W,x)$ and $x_1(W,x)$ are  $\mathcal{F}_x$-measurable.
We define recursively for $k\in\N$, (with $\inf\emptyset=+\infty$ and $\sup\emptyset=-\infty$)
\begin{align*}
    u_{2k+1}(y)
 := &
    \inf\{t>z_{2k}(y),\  Z(t)-\inf\{Z(u), z_{2k}(y)\leq u\leq t\}\geq y\}\un_{\{Z(z_0)>Z(z_\infty)\}}\\
&
    +\inf\{t>z_{2k}(y),\  \sup\{Z(u), z_{2k}(y)\leq u\leq t\}-Z(t)\geq y\}
    \un_{\{Z(z_0)<Z(z_\infty)\}}
    ,\\
    z_{2k+1}(y)
 := &
    [\inf\{t>z_{2k}(y),\ Z(t)= \inf\{Z(u), z_{2k}(y)\leq u \leq u_{2k+1}(y)\}\}\wedge z_\infty] \un_{\{Z(z_0)>Z(z_\infty)\}}\\
&
    +
    [\inf\{t>z_{2k}(y),\ Z(t)= \sup\{Z(u), z_{2k}(y)\leq u \leq u_{2k+1}(y)\}\}\wedge z_\infty]\un_{\{Z(z_0)<Z(z_\infty)\}}
    ,\\
    u_{2k+2}(y)
 := &
    \inf\{t>z_{2k+1}(y),\  \sup\{Z(u), z_{2k+1}(y)\leq u\leq t\}-Z(t)\geq y\}\un_{\{Z(z_0)>Z(z_\infty)\}}\\
&
    +\inf\{t>z_{2k+1}(y),\  Z(t)-\inf\{Z(u), z_{2k+1}(y)\leq u\leq t\}\geq y\}\un_{\{Z(z_0)<Z(z_\infty)\}}
    ,
\\
    z_{2k+2}(y)
 := &
    [\inf\{t>z_{2k+1}(y),\ Z(t)= \sup\{Z(u), z_{2k+1}(y)\leq u \leq u_{2k+2}(y)\}\}\wedge z_\infty]\un_{\{Z(z_0)>Z(z_\infty)\}}
    \\
&
    +[\inf\{t>z_{2k+1}(y),\ Z(t)= \inf\{Z(u), z_{2k+1}\leq u \leq u_{2k+2}(y)\}\}\wedge z_\infty]\un_{\{Z(z_0)<Z(z_\infty)\}}
    .
\end{align*}
Consequently, all these r.v. $z_i(y)$, $i\geq 0$ are $\mathcal{F}_x$-measurable
and so are the r.v. $Z(z_k(y))$, $k\in\N$.
Moreover it follows from the definition of $y$ and $y$-extrema that the
$y$-extrema in $[x_0(W,x),x_1(W,x)]$ are exactly the $z_k(y)$, $k\in\N$ (with repetitions at $z_\infty$).
In particular, $x_0(W,y)=\sum_{k\in\N} z_k(y) \un_{\{z_k(y)\leq 0<z_{k+1}(y)\}}$
and $x_1(W,y)=\sum_{k\in\N}z_{k+1}(y) \un_{\{z_k(y)\leq 0<z_{k+1}(y)\}}$
are $\mathcal{F}_x$-measurable.
$\hfill\Box$

We now prove the following lemma, which is useful in the proof of Lemma \ref{LemmaProbaStrongChangeOfSign}
between equations \eqref{eqProbaBCR} and \eqref{eqProbaBCR2},
in particular to show the independence used in \eqref{eqProbaBCR2}:

\begin{lemma}\label{LemmaMusurabiliteFX+}
For $k\geq 1$, the random variables $e[T_i(X_k)]$, $i\in\{-1,1\}$ are measurable with respect to $\mathcal{F}_{X_{k}+}$,
where $\mathcal{F}_{X_{k}+}=\{A\in\mathcal{F}_\infty,\ \forall x\geq 0,\ A\cap\{X_{k}<x\}\in\mathcal{F}_x\}$.
\end{lemma}

\noindent{\bf Proof:}
We use the same notation as in the previous proof.
Let $k\geq 1$ and $0<y< x$.
We define $K(y):=\sum_{\ell\in\N}  \ell\un_{\{z_\ell(y)\leq 0<z_{\ell+1}(y)\}}$,
so $x_i(W,y)=z_{K(y)+i}(y)$ for every $i\in\Z$ such that
$x_i(W,y)\in[x_0(W,x),x_1(W,x)]$,
and $K(y)$ is $\mathcal{F}_x$-measurable.
For $i\in\Z$ (with $z_j(y):=z_0(y)$ for $j<0$),
\begin{equation}\label{defdehi}
    h_i(y)
 :=
    \big|Z\big(z_{K(y)+i}(y)\big)-Z\big(z_{K(y)+i+1}(y)\big)\big|
=
    \sum_{k\in\N}\un_{\{K(y)=k\}}\big|Z\big(z_{k+i}(y)\big)-Z\big(z_{k+i+1}(y)\big)\big|
\end{equation}
is also $\mathcal{F}_x$-measurable (for every $0<y<x$).
And $h_i(y)=H(T_i(y))$ if the support of the slope $T_i(y)$ is included in $[x_0(W,x),x_1(W,x)]$, since
in this case, $Z\big(z_{K(y)+i}(y)\big)=Z(x_i(W,y))=W(x_i(W,y))$
and $Z\big(z_{K(y)+i+1}(y)\big)=Z(x_{i+1}(W,y))=W(x_{i+1}(W,y))$.

%\noindent
%We start with $e(T_1(X_k))$. Since $e(T_1(X_k))=H(T_1(X_k))-X_k$ and $X_k$ is $\mathcal{F}_{X_k}$ and
%then $(\mathcal{F}_{X_{k}+})$-measurable, we just have to prove that $H(T_1(X_k))$ is $(\mathcal{F}_{X_{k}+})$-measurable.

\noindent
We first prove that  $H(T_1(X_k))$ is ($\mathcal{F}_{X_k+}$)-measurable.
Let $a\in\R$; we have to prove that $\{H(T_1(X_k))\leq a\}\in(\mathcal{F}_{X_{k}+})$, which means that $\{H(T_1(X_k))\leq a\}\cap\{X_{k}< x\}\in\mathcal{F}_x$ for every $x\geq 0$. This is obvious for $0\leq x<c$ since $X_k\geq c$ a.s. We now fix $x\geq c$ and define for $p>1/c$
($h_1(u)$ is defined in \eqref{defdehi} for $0<u<x$, and we set $h_1(u):=0$ if $u\leq 0$)
$$
    D_p(x)
:=
    \sum_{i=1}^\infty h_1(x-i/p)\un_{\{0< x-i/p\}}\un_{\{x-i/p\leq X_k\}}
    \un_{\{X_k< x-(i-1)/p\}},
$$
which is $\mathcal{F}_x$-measurable. Moreover, on $\{X_k< x\}$, there exists a unique (random) $j=j(p)\geq 1$ such that
$x-j/p\leq X_k< x-(j-1)/p\leq x$, and then $x-j/p>0$ since $X_k\geq c>1/p$.
We have
\begin{equation}\label{InegInclusionSlopesyp}
[x_{-1}(W, x-j/p), x_{2}(W, x-j/p)]\subset
[x_{-1}(W, X_k), x_{2}(W, X_k)]
%\subset
%[x_{0}(W, X_{k+1}), x_{1}(W, X_{k+1})]
\subset
[x_{0}(W, x), x_{1}(W, x)].
\end{equation}
Indeed, the last inclusion comes from the fact that $X_k$ is a change of sign of $b(.)$, and $x>X_k$,
so $e(T_0(X_k))=0$ and $x_0(W,X_k)$ and $x_1(W,X_k)$ are not $x$-extrema
%(as says Cheliotis, the slopes $T_{-1}(X_k)$, $T_{0}(X_k)$ and $T_{1}(X_k)$ are glued together for $x>X_k$).

\noindent
Let $y_p:=(x-j(p)/p)\un_{\{X_k<x\}}$. So on $\{X_k<x\}$, $D_p(x)=h_1(y_p)=H(T_1(y_p))$
%$h_1(x-i/p)=H(T_1(x-i/p))$
(see the comments after \eqref{defdehi} since
the support of slope $T_1(y_p)$ is included in $[x_0(W,x),x_1(W,x)]$ by \eqref{InegInclusionSlopesyp}).
%And in particular,
Since $y_p\in(X_k-1/p,X_k]$,  $y_p\to_{p\to+\infty} X_k$ on $\{X_k<x\}$,
and since $H(T_1(.))$ is left-continuous on $(0,+\infty)$ on $\mathcal{W}$ by Lemma \ref{LemmaLeftcontinuity},
$H(T_1(X_k))=\lim_{p\to+\infty}H(T_1(y_p))=\lim_{p\to+\infty}D_p(x)$ on $\{X_{k}< x\}$.
Hence,
$$
    \{H(T_1(X_k))\leq a\}\cap\{X_k< x\}
=
    \big\{\lim_{p\to+\infty}D_p(x)\leq a\big\}
    \cap\{X_k< x\}.
$$
Since $\lim_{p\to+\infty}D_p(x)$ is the limit of a sequence of $\mathcal{F}_x$-measurable r.v., it is also
$\mathcal{F}_x$-measurable, and then $\{\lim_{p\to+\infty}D_p(x)\leq a\}\in \mathcal{F}_x$.
Since $\{X_k< x\}\in \mathcal{F}_x$, we get
$\{H(T_1(X_k))\leq a\}\cap\{X_k<x\}\in \mathcal{F}_x$, and this is true for every $x\geq 0$.
So $\{H(T_1(X_k))\leq a\}\in \mathcal{F}_{X_k+}$ for every $a\in\R$.

\noindent
Hence $H(T_1(X_k))$ and then $e(T_1(X_k))$ are ($\mathcal{F}_{X_k+}$)-measurable.
Finally, we show similarly that $H(T_{-1}(X_k))$ and then $e(T_{-1}(X_k))$ are ($\mathcal{F}_{X_k+}$)-measurable.
%II do not know if it is necessary but so is $e(T_0(X_k))$ since it is $0$ by definition !
\hfill$\Box$

\bigskip

\noindent\textbf{Acknowledgment : }
I am grateful to an anonymous referee for  reading the paper very carefully.

%notations : $\a''$, $k_1$, $k$ $\gamma$, $F_i(N)$ (bof); cf 25/08/08. et $q'$,  $Q$, $E_{9,i}$, ... $E_{7,i}$;

%%%%%%%%%%%%%%%%%%%%%%%%%%%%%%%%%%%%%%%%%%%%%%%%%%%%%%%%%%%
%                                                         %
%                     BIBLIO                              %
%                                                         %
%%%%%%%%%%%%%%%%%%%%%%%%%%%%%%%%%%%%%%%%%%%%%%%%%%%%%%%%%%%


\begin{thebibliography}{99}

{\baselineskip=14pt

\bibitem{Andreoletti_PhD} Andreoletti, P.: Localisation et Concentration de la Marche de Sinai.
Ph.D. thesis, Universit\'e Aix-Marseille II, 2003, available at \texttt{http://tel.archives-ouvertes.fr/tel-00004116}.




\bibitem{Papier_4_A}
    Andreoletti, P.: Alternative proof for the localization of {S}inai's walk. {\it J. Stat. Phys.}
 {\bf 118} (2005), 883--933.


\bibitem{Andreoletti_Devulder_1}
    Andreoletti, P. and Devulder A.: Localization and number of visited valleys for a transient diffusion in random environment.
    Preprint ArXiv (2013), 	arXiv:1311.6332.


\bibitem{AurzadaSimon}
Aurzada, F. and Simon, T.: Persistence probabilities \& exponents. Preprint ArXiv (2012), arXiv:1203.6554.

\bibitem{Bovier_Faggionato}
    Bovier, A. and Faggionato, A.: Spectral analysis of {S}inai's walk for small eigenvalues.
    {\it Ann. Probab.} {\bf 36} (2008), 198--254.




\bibitem{Bray}
    Bray A. J., Majumdar S. N. and Schehr G.: Persistence and First-Passage Properties in Non-equilibrium Systems. 
    {\it Advances in Physics} {\bf 62} (2013), 225--361.
%Preprint ArXiv (2013),  arXiv:1304.1195.



\bibitem{Brox} Brox, Th.: A one-dimensional diffusion process in a {W}iener
medium. {\it Ann. Probab.} {\bf 14} (1986), 1206--1218.


\bibitem{Cheliotis} Cheliotis, D.: Diffusion in random environment and the renewal theorem.
              {\it Ann. Probab.}  {\bf 33} (2005), 1760--1781.

\bibitem{Cheliotis_Favorite}
    Cheliotis, D.: Localization of favorite points for diffusion in a random environment.
    {\it Stoch. Proc. Appl.} {\bf 118} (2008), 1159--1189.


\bibitem{CoccoMonasson} Cocco, S. and Monasson, R.: Reconstructing a random potential from its random walks.
              {\it Europhysics Letters}  {\bf 81} (2008), 20002.

\bibitem{DemboGao}
    Dembo A., Ding, J.   and Gao F.: Persistence of iterated partial sums. 
    {\it Ann. Inst. H. Poincaré Probab. Stat.} {\bf 49} (2013), 873--884. 

\bibitem{Devulder} Devulder, A.:
Some properties of the rate function of quenched large deviations for random walk in random environment. {\it Markov Process. Related Fields} {\bf 12} (2006), 27--42.

\bibitem{Devulder_SPL} Devulder, A.: The speed of a branching system of random walks in random environment. {\it Statist. Probab. Lett.} {\bf 77} (2007), 1712--1721.

\bibitem{ELS} Enriquez N., Lucas C. and Simenhaus F.:
The Arcsine law as the limit of the internal DLA cluster generated by Sinai's walk.
{\it Ann. Inst. H. Poincaré Probab. Stat.} {\bf 46} (2010), 991-1000.


\bibitem{Golosov}  Golosov, A. O.:
Localization of random walks in one-dimensional random environments.
{\it Comm. Math. Phys.} {\bf 92} (1984),
491--506.



\bibitem{LMF} Le Doussal P., Monthus C., Fisher D.: Random walkers in one-dimensional
random environments; Exact renormalization group analysis.
              {\it Phys. Rev. E}  {\bf 59} (1999), 4795--4840.

\bibitem{HuLocal} Hu Y.:
Tightness of localization and return time in random environment.
{\it Stoch. Proc. Appl.} {\bf 86} (2000), 81--101.


\bibitem{Hug} Hughes, B.D.: {\it Random Walks and Random Environment,
vol. II: Random Environments}. Oxford Science Publications,
Oxford, 1996.

\bibitem{KMT} Koml{\'o}s, J., Major, P. and Tusn{\'a}dy, G.: An
approximation of partial sums of independent {${\rm RV}$}'s
              and the sample {${\rm DF}$. {I}}.
              {\it Z. Wahrsch. Verw. Gebiete}  {\bf 32} (1975), 111--131.

\bibitem{NP} Neveu J. and Pitman J.: Renewal property of the extrema and tree property of the
              excursion of a one-dimensional {B}rownian motion.
              {\it S\'eminaire de {P}robabilit\'es XXIII, Lecture Notes in Math.}   {\bf 1372} (1989), 239--247,
              Springer, Berlin.


\bibitem{Revesz} R{\'e}v{\'e}sz, P.: {\it Random walk in random and non-random environments}, second edition.
World Scientific, Singapore, 2005.

\bibitem{RevuzYor}
    Revuz, D. and Yor, M.:
    {\it Continuous Martingales and {B}rownian Motion}, second
    edition.
    {Springer}, Berlin, {1994}.
    %Attention, voir les ref de 99 et non pas 94...



\bibitem{Schumacher}
    Schumacher, S.:
    Diffusions with random coefficients.
    {\it Contemp. Math.} \textbf{41} (1985),  351--356.

\bibitem{S2}
     Shi, Z.:
     Sinai's walk via stochastic calculus. {\it Panoramas et Synth\`eses} {\bf 12} (2001), 53--74,
     Soci\'et\'e math\'ematique de France.

\bibitem{Simon}
    Simon, T: The lower tail problem for homogeneous functionals of stable processes with no negative jumps.
    {\it ALEA Lat. Am. J. Probab. Math. Stat.} {\bf 3} (2007), 165--179.


\bibitem{S3}
    Sinai, Ya. G.:
    The limiting behavior of a one-dimensional random walk in a random
    medium. {\it Th. Probab. Appl.} {\bf 27}
    (1982), 256--268.

\bibitem{SinaiSRW}
    Sinai, Ya. G.:
    Distribution of some functionals of the integral of a random walk.
    {\it Theoret. and Math. Phys.} {\bf 90}
    (1992), 219--241.


\bibitem{S1}
    Solomon, F.:
    Random walks in a random environment.
    {\it Ann. Probab.} {\bf 3} (1975), 1--31.

\bibitem{tanaka}
    Tanaka, H.:  Localization of a diffusion process in a one-dimensional Brownian environment.
    {\it Comm. Pure Appl. Math.} {\bf 47} (1994), 755--766.

\bibitem{Vito}
    Vysotsky, V.:
    On the probability that integrated random walks stay positive.
    {\it Stoch. Proc. Appl.} {\bf 120} (2010), 1178--1193.

\bibitem{Z1}
    Zeitouni, O.:
    Lectures notes on random walks in random environment.
   In: {\it Lect. Notes Math.}  {\bf 1837}
    193--312,  Springer, Berlin 2004.

}
     %%% end of "baselineskip=14pt"

\bibitem{Zindy} Zindy, O.: Upper limits of {S}inai's walk in random scenery. {\it Stoch. Proc. Appl.} {\bf 118}
      (2008), 981--1003.


\end{thebibliography}
\end{document}